\documentclass[twoside,a4paper,reqno,intlimits]{amsart}%
\usepackage[%
   colorlinks=true,
   urlcolor=rltblue,       
   filecolor=rltgreen,     
   citecolor=rltgreen,     
   linkcolor=rltred,       
   bookmarks=true,
   unicode,
   plainpages=false,
   ]{hyperref}
   \usepackage{mathtools}

   \newcommand{\pot}{U}
\newcommand{\potV}{U}
\newcommand{\Grenze}{g}
\DeclareMathOperator{\spt}{spt}
\newcommand{\function}{\omega}

\newcommand{\multappomega}{\function^{N}}

\newcommand{\appAq}[1]{{#1}^{A,q}}
\newcommand{\MC}{a}
\newcommand{\para}{{\delta}}
\newcommand{\diver}{\operatorname{div}}
\newcommand{\dist}{\operatorname{dist}}
\newcommand{\tran}[1]{{#1}_{-\tau}}
\newcommand{\trap}[1]{{#1}_{\tau}}
\newcommand{\difp}[1]{d^+{#1}}
\newcommand{\trapm}[1]{{#1}_{\pm\tau}}
\newcommand{\difpm}[1]{d^\pm{#1}}
\newcommand{\difn}[1]{d^-{#1}}
\newcommand{\intO}{\int\limits_{\Omega}}
\newcommand{\otimess}{\overset{s}{\otimes}}
 \newcommand{\td}{\partial_{\tau}}

\newcommand{\bue}{\bu^{N}}
\newcommand{\buen}[1]{\bu^{#1}}
\newcommand{\bFn}[1]{\bF^{#1}}
\newcommand{\Sn}[1]{S^{#1}}
\newcommand{\bSn}[1]{\bS^{#1}}
\newcommand{\param}{\lambda}

\usepackage[active]{srcltx}
\usepackage{color,esint}
\usepackage{amsmath,amsthm}
\usepackage{amssymb}
\usepackage{amsfonts}
\usepackage{enumitem}
\setenumerate{label={\rm (\alph{*})}}
\usepackage{graphicx}

\usepackage{xcolor}
\definecolor{rltred}{rgb}{0.75,0,0}
\definecolor{rltgreen}{rgb}{0,0.5,0}
\definecolor{rltblue}{rgb}{0,0,0.75}

\definecolor{egreen}{rgb}{0,0.6,0}

 \hypersetup{%
   pdftitle={}, pdfsubject={}, pdfauthor={}, pdfkeywords={}, pdfcreator={pdfLaTeX},
   pdfproducer={pdfLaTeX}, }

\newtheorem{theorem}[equation]{Theorem}
\newtheorem{lemma}[equation]{Lemma}
\newtheorem{proposition}[equation]{Proposition}
\newtheorem{corollary}[equation]{Corollary}
\newtheorem{definition}[equation]{Definition}
\newtheorem{remark}[equation]{Remark}
\numberwithin{equation}{section}

\usepackage[rose,frak,operator]{paper_diening}

\begin{document}
\title[Natural second-order regularity for elliptic/parabolic
problems] {Natural second-order regularity for parabolic systems with
  operators having $(p,\delta)$-structure and depending only on the
  symmetric gradient}
\author{Luigi C.\ Berselli}
\address{Dipartimento di Matematica, Universit{\`a} di Pisa, Via F.~Buonarroti 1/c,
  I-56127 Pisa, ITALY.}  \email{luigi.carlo.berselli@unipi.it}

\author{Michael R\r u\v zi\v cka{}} 
\address{Institute of Applied Mathematics, Albert-Ludwigs-University Freiburg,
  Ernst-Zermelo-Str.~1, D-79104 Freiburg, GERMANY.}
\email{rose@mathematik.uni-freiburg.de}

\begin{abstract}
  In this paper we consider parabolic problems with stress tensor
  depending only on the symmetric gradient. By developing a new
  approximation method (which allows to use energy-type methods
  typical for linear problems) we provide an approach to obtain global 
  regularity results valid for general potential operators with
  $(p,\delta)$-structure, for all $p>1$ and for all $\delta>0$. In
  this way we prove ``natural'' second order spatial regularity --up to
  the boundary-- in the case of homogeneous Dirichlet boundary
  conditions. The regularity results, are presented with full details
  for the parabolic setting in the case $p>2$. However, the same method
  also yields regularity in the elliptic case and for
  $1<p\leq 2$, thus proving in a different way results already known.
\end{abstract}
\date{\small \today}
\maketitle

\section{Introduction}
In this paper we consider an initial boundary value problem for
general nonlinear parabolic systems
\begin{align}
  \label{eq:pfluid}
  \begin{aligned}
    \frac {\partial \bfu}{\partial t}-\divo \bfS (\bfD\bfu) &= \bff
\qquad&&\text{in }
    I\times \Omega,
    \\
    \bu &= \bfzero &&\text{on } I\times\partial \Omega\,,
    \\
    \bu(0)&=\bu_0&&\text{in }\Omega\,,
  \end{aligned}
\end{align}
where the operator $\bS$ depends only on the symmetric
gradient~$\bD\bu= \frac 12({(\nabla\bu)^{\top}\!\!+\nabla\bu})$ and has
$(p,\delta)$-structure (cf.~Definition \ref{def:ass_S}). Here
$I:=(0,T)$ for some $T>0$ is a finite time interval, and
$\Omega\subset\setR^{3}$ is a sufficiently smooth, bounded domain.
The paradigmatic example for the operator in \eqref{eq:pfluid} is
given via
\begin{equation}
  \label{eq:example}
  \mathbf{S}(\bD\bu):=(\delta+|\bD\bu|)^{p-2}\bD\bu\qquad
  \delta\geq0,\ 1<p<\infty. 
\end{equation}
In this paper we only treat the case $p>2$. However, the method of
proof, based on an $(A,q)$-approximation (cf.~Section
\ref{sec:approx}) works for every $p\in (1,\infty)$. We focus to the
case $p >2$, since our main result in the case $p\in (1,2]$ has been
already proved in a different way (cf.~\cite{br-plasticity}) and the
method of the present paper simplifies a lot for these exponents. Note that
the elliptic problem corresponding to \eqref{eq:pfluid} can be treated
in the same way with much shorter proofs. Moreover, all our result
possess corresponding analogues in $d$-dimensional domains $\Omega
\subset \setR^d$, $d\ge 2$. For simplicity we only treat the case
$d=3$.

Our main goal is to prove a result of ``natural'' second-order spatial 
regularity for weak solutions. This corresponds to proving,
under appropriate (minimal) assumptions on the data, that weak
solutions satisfy
\begin{equation*}
\int\limits_I\int\limits_ \Omega (\delta+|\bD\bu|)^{p-2}|\nabla\bD\bu|^{2}\,d\bx
  \,ds\leq C,
\end{equation*}
which can be also equivalently  re-written as
$ \bF(\bD\bu) \in L^{2}(I;W^{1,2}(\Omega)) $ with 
\begin{equation}
  \label{eq:F}
  \bF(\bD\bu):=(\delta+|\bD\bu|)^{\frac{p-2}{2}}|\bD\bu|.
\end{equation}
We say ``natural'' as opposed to some recent results  proving
$\bS\in L^{2}(I;W^{1,2}(\Omega))$, which is equivalent to proving that  
\begin{equation*}
\int\limits_{I}\int\limits_{\Omega}\big
|\nabla\big((\delta+|\bD\bu|)^{p-2}\bD\bu\big)\big |^{2}\,d\bx
  \,ds\leq C,
\end{equation*}
which is called ``optimal'' second-order spatial regularity.  The two
notions of regularity are rather different in the spirit: the optimal
regularity is linked with nonlinear versions of the singular integral
theory, while the natural regularity is based on energy methods. This
yields estimates in quasi-norms, which are
of crucial relevance especially for the numerical analysis of the
problem, and in particular to study optimal convergence rates of spatial
discretizations (cf.~Barrett and Liu~\cite{baliu}).

The problem has a long history and many result concern mainly the
problem: a) in the scalar or elliptic case; b) with operators $\bS$
depending on the full gradient; c) the interior regularity. We refer
to the classical results by DiBenedetto~\cite{dibene}, Gilbarg and
Trudinger~\cite{gilbarg-trudinger}, Lady{\v{z}}henskaja et
al.~\cite{LU1968,la-sol-ur-parabolic}, Liebermann~\cite{Lie1988},
Uhlenbeck~\cite{Uhl1977}, Ural'ceva~\cite{Ura1968}, just to cite a
few; or the ones linked more to applications Bensoussan
and Frehse~\cite{BF2002}, Ne\v{c}as~\cite{necas-83}, and Fuchs and
Seregin~\cite{Fuc-Ser}. Even if the studies started in the sixties, we observe that the field is
still extremely active and very recent results are those
in~\cite{BM2020, BdVC2012, CM2019, CM2020}.



Our treatment of the case of systems with dependence only on the
symmetric gradient and up-to-the boundary is new, to the best of the
author's knowledge. We extend the so called $A$-approximation
technique from \cite {mnr3} such that it allows a treatment of all exponents $p
\in (1,\infty)$.
Here, we focus on the regularity of the quantity
in~\eqref{eq:F}. Thus, this work can be seen as a natural extension of
previous results we have done in the case $p \in (1,2]$ for the steady
problem in~ \cite{br-plasticity} and for the unsteady
continuous/discrete in~\cite{br-parabolic}. Note that our approach
allows to treat the full range of exponents $p \in (1,\infty)$, as in the
scalar case, even if we give full details only in the case $p>2$, as
the case $p \in (1,2)$ is already treated in a different way. Notice 
that the results in \cite{CM2020} hold only for $p>\frac{3}{2}$, which
has been improved in \cite{balci2021pointwise}, reaching
$p>4-2\sqrt{2}$. The limitation on $p>3/2$ was also present
 in prior results of ``natural'' regularity in the symmetric gradient 
case \cite{hugo-thin-nonflat}, but it has then later removed completely
in \cite{br-reg-shearthin} to the case $p>1$.

The techniques employed for $p>2$ are rather different from those
previously used in the case for $p<2$, where calculations can be more
easily justified by approximation of the system by means of adding the
term $-\epsilon\Delta\bu$ (and then showing that estimates for a
system with leading linear part could be made independent of
$\epsilon>0$). Anyway, the technique we use can be also employed in
the case $p \in (1,2]$ to prove in an alternative way the regularity results
already known. This requires some technical adjustments which are
left for a further investigation, since the technicalities are complex
enough already in the case $p>2$. The introduction of a different
regularization of the problem is due to the fact that for 
$p>2$ the  perturbation with the heat equation is not
enough to justify the computations; hence, we developed a new
(multiple) approximation technique, by
a sequence of operators, such that the last is an affine one, which
allows to use standard energy techniques leading to $W^{2,2}$-results.

\subsection{Sketch of the proof of the main result}\label{sec:ske}
To prove the main regularity result (cf.~Theorem \ref{thm:MT}) we
proceed as follows: (a) we introduce a proper multiple approximation
of the operator $\bS$; (b) we prove interior and tangential estimates
for second order derivatives by difference quotient methods; (c) we
use the equations point-wise to recover the remaining derivative; (d)
make again use of the point-wise equations and integration by parts
in the full domain to obtain estimates independent of the
approximation parameters; (e) and finally we pass to the limit with the multiple
approximation parameters.

For the reader's convenience, we explain here the main ideas in the
case that the operator $\mathbf{S}$ is given by
\eqref{eq:example} and that instead of \eqref{eq:pfluid} its
steady counterpart is treated. Most of the calculations are elementary, but
involved, and use  
various well-established techniques
from the regularity theory of partial differential equations. Since they
are linked in a quite intricate and delicate way and one has to be
careful in tracking the dependence on various parameters, we sketch the proof
now and then develop a full theory in the next sections.

A fundamental step in the approximation of general operators by
ones with linear growth dates back to \cite{mnr3}, where generalized Newtonian
fluids are treated. The results proved there are obtained by using for
$A\geq1$ the
following approximation\footnote{The precise form of the approximation
  in \cite{mnr3} is slightly different, since there the potential of
  the stress tensor was depending on $|\bD\bu|^{2}$, instead of
  $|\bD\bu|$ here.} $\bS^A$ defined via
\begin{equation*}
  \mathbf{S}^{A}(\bP)=\left\{
    \begin{aligned}
   &   (\delta+|\bP^{\sym}|)^{p-2}\bP^{\sym}&\quad \text{if }|\bP|\leq A\,,
      \\
     & c_{2}\bP^{\sym}+c_{1}&\quad \text{if }|\bP|> A\,,
    \end{aligned}
  \right.
\end{equation*}
with appropriately chosen constants $c_{i}=c_{i}(A,\delta,p)$ to
ensure an appropriate regularity of the stress tensor $\bS^A$. Hence,
the tensor $\bS^A$ grows
linearly for large $\bP$. This can be also restated by writing that
\begin{equation*}
  \mathbf{S}^{A}(\bP):=\frac{(\function^{A})'(|\bP^{\sym}|)
  }{|\bP^{\sym}|}\bP^{\sym}, 
\end{equation*}
where $\function^{A}:\setR^{\geq0}\to\setR^{\geq0}$ is a regular N-function such that
$(\function^{A})'(0)=0$, $(\function^{A})'(t)=(\delta+t)^{p-2}t$ for $t\leq A$ and
$(\function^{A})'(t)=c_{2}t+c_{1}$ for $t> A$.
\begin{remark}
  \label{rem:simple}
  In Section 2 we will show that -roughly speaking- once the results
  is established for this explicit example, then it can be extended to
  a rather wide class of nonlinear operators.
\end{remark}

To obtain results for the original problem we first consider 
the approximate problem 
\begin{align}
  \label{eq:pfluid-steady}
  \begin{aligned}
-\divo \mathbf{S}^{A} (\bfD\bfu^{A}) &= \bff
\qquad&&\text{in }
     \Omega,
    \\
    \bu^{A} &= \bfzero &&\text{on } \partial \Omega\,.
  \end{aligned}
\end{align}
For regular enough $\bff$ one can directly prove the existence of weak
solutions satisfying 
\begin{equation*}
\int\limits_{\Omega}|\bF^{A}(\bD\bu^{A})|^{2}\,d\bx\leq C,
\end{equation*}
where
\begin{equation*}
  \bF^{A}(\bP):=\sqrt{\frac{(\function^{A})'(|\bP^{\sym}|)}{|\bP^{\sym}|}}\bP^{\sym}.
\end{equation*}
Note that $|\bF^{A}(\bP)|^{2}\sim
\delta^{2}+|\bP^{\sym}|^{2}$, with constants depending on $A$. The special role of the quantity
\begin{equation*}
a^{A}(t):=\frac{(\function^{A})'(t)}{t}\,,
\end{equation*}
is evident from the definitions of $\mathbf{S}^{A}$ and $\bF^{A}$.

The estimates for the second order spatial derivatives are obtained by
using the difference quotient technique in the interior and along tangential
directions (after appropriate localization of the equations). Once
this step is done, one gets that the equations are satisfied almost
everywhere. Thus, the equations can be used point-wise to determine (by
ellipticity) estimations in the direction normal to the boundary. 
The outcome of this procedure, which is typical for second order
elliptic equations, leads to the estimates (cf.~Proposition
\ref{prop:JMAA2017-1}, Proposition \ref{prop:JMAA2017-2})
\begin{equation*}
  \begin{aligned}
 \delta^{p-2}
\int\limits_{\Omega_{0}}|\nabla\bD\bu^{A}|^{2}\,d\bx\leq\int\limits_{\Omega_{0}}|\nabla\bF^{A}(\bD\bu^{A})|^{2}\,d\bx 
  &\leq C_{1}
\qquad \forall \,\Omega_{0}\subset\subset \Omega,
  \\
   \delta^{p-2} \int\limits_{\Omega}|\nabla\bD\bu^{A}|^{2}\,d\bx\leq\int\limits_{\Omega}|\nabla\bF^{A}(\bD\bu^{A})|^{2}\,d\bx
   &\leq C_{2}(A),
\end{aligned}
\end{equation*}
where the constant $C_{1}$ is independent of $A$. In addition, one
gets that also tangential derivatives are regular up to the boundary
with a bound independent of $A$. Note that the linear growth of the
operator $\bS^A$ results in an $L^2$-setting, which allows us to use
the classical Korn inequality and to handle the dependence of the
operator on the symmetric gradient (instead of on the full gradient)
in the equations. An important feature of this step is that the proved
regularity is sufficient to justify the following step and to remove
the dependence on $A$ in the estimates in the direction normal to the
boundary. 

This is achieved by testing the equations locally near the boundary by
second order derivatives in the normal direction, and adapting a method
introduced by Seregin and Shilkin~\cite{SS00} for $1<p<2$
(cf.~\cite{br-plasticity,br-parabolic}). This results in the estimate
(cf.~Proposition~\ref{prop:main}, Proposition \ref{thm:estimate_for_ue})
\begin{equation*}
 \delta^{p-2}
 \int\limits_{\Omega}|\nabla\bD\bu^{A}|^{2}\,d\bx\leq\int\limits_{\Omega}|\nabla\bF^{A}(\bD\bu^{A})|^{2}\,d\bx \leq C_{3},
\end{equation*}
for some $C_{3}$ which is independent of $A$.

The final step is the passage to the limit $A\to\infty$. By uniform
boundedness it directly follows that $\bF^{A}(\bD\bu^{A})$ has a weak
limit $\widehat{\bF}\in W^{1,2}(\Omega)$ and by using also the uniform bound on
second order derivatives, it follows that $\bD\bu^{A}\to\bD\bu$ almost
everywhere. Combining these two information, the definition of $\bF^{A}$,
and the lower semi-continuity of the norm it follows that
\begin{gather*}
\widehat{\bF}=   \lim_{A\to\infty}\bF^{A}(\bD\bu^{A})=\bF(\bD\bu)\qquad\text{
      weakly in }W^{1,2}(\Omega) \text{ and a.e. in }\Omega,
    \\
    \int\limits_{\Omega}|\nabla\bF(\bD\bu)|^{2}\,d\bx\leq C_{3}.
\end{gather*}
It remains to prove that $\bu$ is the unique solution of the
steady version the original
problem~\eqref{eq:pfluid}. From the construction of $\bS^A$ follows 
$\mathbf{S}^{A}(\bP)\to\mathbf{S}(\bP)$ for every
$\bP\in \setR^{3\times 3}$. This fact, coupled with the almost everywhere
convergence of $\bD\bu^{A}$, implies that 
\begin{equation*}
  \lim_{A\to\infty}\mathbf{S}^{A}(\bD\bu^{A}(\bx))\to\mathbf{S}(\bD\bu(\bx))
  \qquad
  a.e.\ \bx\in \Omega, 
\end{equation*}
which is nevertheless \textit{not} enough to infer directly that
\begin{equation*}
\lim_{A\to\infty}
\int\limits_{\Omega}\mathbf{S}^{A}(\bD\bu^{A})\cdot\bD\bw\,d\bx=\int\limits_{\Omega}\mathbf{S}(\bD\bu)\cdot
\bD\bw\,d\bx\qquad \forall 
\,\bw\in C^{\infty}_{0}(\Omega),
\end{equation*}
and to pass to the limit in the weak formulation. To this end we need
-for instance- additionally an uniform bound on $\mathbf{S}^{A}(\bD\bu^{A})$ in
$L^{q}(\Omega)$ for some $q>1$. This implies that
$\mathbf{S}^{A}(\bD\bu^{A})\weakto \widehat{\mathbf{S}}$ in
$L^{q}(\Omega)$, and that the limit can be identified as
$\widehat{\mathbf{S}}=\mathbf{S}(\bD\bu)$, by a classical result.

Observe that from the definition of $\bS^A$ it follows
(cf.~Proposition \ref{prop:SA-ham}, Lem\-ma~\ref{lem:UAm}, Lemma
\ref{lem:UA}) that
\begin{equation*}
|  \mathbf{S}^{A}(\bD\bu^{A})|\leq \left\{
  \begin{aligned}
    &c\,(\delta^{p-1}+|\bD\bu^{A}|^{p-1})&\quad p>2,
    \\
    &c\,\delta^{p-2}|\bD\bu^{A}|&\quad 1<p\le 2,
  \end{aligned}
\right.
\end{equation*}
while the proved estimate $ \bF(\bD\bu^{A})\in W^{1,2}(\Omega)$, which
is uniformly with respect to $A$, implies by Sobolev embedding (in three-dimensions)
that $\|\bF(\bD\bu^{A})\|_{6}\leq C $. Using the properties  of $\bF^{A}$,
it follows that (cf.~Proposition \ref{prop:SA-ham},
Lem\-ma~\ref{lem:UAm}, Lemma~\ref{lem:UA}) 
\begin{equation*}
\begin{aligned}
    \|\bD\bu^{A}\|_{{6}}&\leq C &\quad p>2,
    \\
    \|\bD\bu^{A}\|_{{3p}}&\leq C&\quad 1<p\le 2.
  \end{aligned}
\end{equation*}
Hence we get that $\mathbf{S}^{A}(\bD\bu^{A})$ is bounded uniformly in
$L^{6/(p-1)}(\Omega)$ for $p>2$ and in $L^{3p}(\Omega)$ for
$1<p \le 2$, which implies that the above argument to pass to the
limit in the weak formulation works only for $1<p<7$.

To remove the restriction $p<7$ in the regularity 
result\footnote{The restriction depends on the space dimension and it
  is more stringent in the time-evolution case, due to different
  parabolic embedding results.}  we introduce and perform a \textit{multiple approximation} of the operator,
which is roughly speaking the following: for given decreasing sequences
$p>q_{1}>q_{2}>\dots>q_{N}=:2$ and $A_{N}>A_{N-1}>\dots >
A_1\ge 1$ we set
\begin{equation*}
  \mathbf{S}^{N}(\bP):=\left\{
    \begin{aligned}
   &   (\delta+|\bP^{\sym}|)^{p-2}\bP^{\sym}&\quad \text{if }|\bP|\leq A_{1}\,,
      \\
     & c_{2,q_{1}}|\bP^{\sym}|^{q_{1}-2}\bP^{\sym}+c_{1,q_{1}}&\quad \text{if
     }A_{1}<|\bP|\leq A_{2}\,,
     \\
     &\qquad \vdots&\vdots\qquad\quad
     \\
     & c_{2,q_{N-1}}|\bP^{\sym}|^{q_{N-1}-2}\bP^{\sym}+c_{1,q_{N-1}}&\quad \text{if
     }A_{N-1}<|\bP|\leq A_{N}\,,
     \\
          & c_{2,q_{N}}\bP^{\sym}+c_{1,q_{N}}&\quad \text{if
     }A_{N}<|\bP|\,,
    \end{aligned}
  \right.
\end{equation*}
where the various constants $c_{i,m}$ are chosen such that the
operator $\bS^N$ belongs to the class $C^{1}$. If the exponents $q_{n}$ are chosen
such that 
\begin{equation*}
  \frac {3q_{n}}{q_{n-1}}>1\qquad n=1,\dots,N\,,
\end{equation*}
it is possible to perform the
limiting process step by step, sending to infinity $A_{N}$ (with
$A_{n}$ for $n\leq N-1$ fixed), then taking the limit $A_{N-1}\to\infty$ with the
previous ones fixed, and so on.
This procedure requires to prove the precise dependence of the lower
and the upper bounds of the multiple approximation with respect to 
the parameters\footnote{Moreover, some care has to be taken in the
choice of the $A_{n}$ to ensure monotonicity of the resulting
potentials.} $A_{n}$. 

\bigskip

\noindent \textbf{Plan of the paper.}  The analysis of the approximate
operators is the content of  Section~\ref{sec:stress-tensor} of the
paper, where the procedure is carried out with
full details for general operators, derived from a potential and
having $(p,\delta)$-structure. Moreover, for the derivation of the
estimates for second derivatives, one also has to handle precisely the 
behavior of the related operators $\bF^{n}$. In particular, we will see
that a peculiar role is played by handling tensors derived from a 
potential $\pot$ satisfying $\pot'(t)/t\sim\pot''(t)$, which we call
\textit{balanced}. This allows us to reduce many of the estimations to
computable explicit cases, cf.~Remark~\ref{rem:simple}.

Next, in Section~\ref{sec:regularity} the existence and regularity for
solution of the approximate problems is treated in detail. 
Particular care is given to the full justification of
the calculations: the results are rather natural from a formal point
of view, while the rigorous treatment of all integrals needs certain
approximations and the application of difference quotients, in order to be sure
that we do not work with infinite quantities. First, some $A_{n}$ dependent
estimates are proved, in order to justify manipulating the system
\eqref{eq:pfluid-steady} point-wise and then to derive uniform
estimates by (improved) generalized energy methods. The limiting
process is carried out in the more technical parabolic case, using
space-time compactness results and convergences (at the price of a more
restrictive choice of the parameters $q_{n}$).

\section{Nonlinear operators and N-functions}
\label{sec:stress-tensor}
The goal of this section is to define an approximation, which
possesses nice properties, for operators appearing
in~\eqref{eq:pfluid}. The approximation is inspired by~\cite{mnr3},
while the proof of its properties is close to
\cite{dr-nafsa}. However, our notions are defined slightly different,
which simplifies and shortens the argumentation.
\subsection{Notation}
We use $c, C$ to denote generic constants, which may change from line
to line, but are not depending on the crucial quantities. Moreover, we
write $f\sim g$ if and only if there exists constants $c,C>0$ such
that $c\, f \le g\le C\, f$.

For a bounded, sufficiently smooth domain $\Omega\subset \setR^3$ we
use the customary Lebesgue spaces $(L^p(\Omega), \norm{\,.\,}_p)$,
$p \in [1,\infty]$, and Sobolev spaces
$(W^{k,p}(\Omega), \norm{\,.\,}_{k,p})$, $p \in [1,\infty]$,
$k \in \setN$. We use the notation $(f,g)=\int_\Omega f g\, d\bx$, whenever
the right-hand side is well defined. We do not distinguish between scalar, vector-valued or
tensor-valued function spaces in the notation if there is no danger of
confusion. However, we denote scalar functions by roman letters,
vector-valued functions by small boldfaced letters and tensor-valued
functions by capital boldfaced letters.  If the norms are considered
on a set $M$ different from $\Omega$, this is indicated in the
respective norms as $\norm{\,.\,}_{p,M}, \norm{\,.\,}_{k,p,M}$. We
equip $W^{1,p}_0(\Omega)$ (based on the \Poincare{} lemma) with the
gradient norm $\norm{\nabla \,.\,}_p$.  We denote by $\abs{M}$ the
$3$-dimensional Lebesgue measure of a measurable set $M$.
As usual the gradient of a vector field
$\bv :\Omega\subset \setR^3 \to \setR^3$ is denoted as
$\nabla \bv = (\partial _i v^j)_{i,j=1,2,3}=(\partial _i
\bv)_{i=1,2,3}$, while its symmetric part is denoted as $\bD\bv:=
\frac 12 \big (\nabla \bv + \nabla \bv ^\top\big )$.
The derivative of functions defined on tensors, i.e., 
$\potV:\setR^{3\times 3}\to \setR$ is denoted as
$\partial \potV = (\partial _{ij} \potV)_{i,j=1,2,3}$ where $\partial _{ij}$
are the partial derivatives with respect to the canonical basis of
$\setR^{3\times 3}$.

\subsection{N-functions}
We start with a discussion of some non-trivial properties of
\mbox{N-functions}  that we need in the sequel. For a detailed discussion of
Orlicz spaces and N-functions we refer to \cite{krasno,Mu,ren-rao,dr-nafsa}.
%
\begin{definition}[N-function and regular N-function]\label{def:N}
  A function $\varphi:\setR^{\geq0}\to\setR^{\geq 0}$ is called {\rm
    N-function} 
  if $\varphi$ is continuous,
  convex, strictly positive for $t>0$, and satisfies\footnote{In the
    following we use the convention that $\frac {\varphi'(0)}{0}:=0$.}
  \begin{equation*}
    \lim_{t\to 0^+}\frac{\varphi(t)}{t}=0\,,\qquad\qquad 
    \lim_{t\to \infty}\frac{\varphi(t)}{t}= \infty\,. 
  \end{equation*}
  If $\varphi$ additionally belongs to $C^1(\setR^{\ge 0})\cap
  C^2(\setR^{> 0})$ and satisfies $\varphi''(t)>0$ for all $t>0$, we call
  $\varphi$ a {\rm  regular N-function}.  
\end{definition}
The use of regular N-functions 
is sufficient for our purposes. Thus, in the rest of the paper we
restrict ourselves to this case. For a treatment in the general
situation we refer to the above mentioned literature. Note that for a
regular N-function we have $\varphi (0)=\varphi'(0)=0$. Moreover,
$\varphi'$ is increasing and $\lim _{t\to \infty} \varphi'(t)=\infty$.

The following notion plays an important role in the sequel.
\begin{definition}[$\Delta_2$-condition]\label{def:Delta2}
  A non-decreasing function $\varphi:\setR^{\geq0}\to\setR^{\geq0}$ is said to satisfy the
  {\rm $\Delta_2$-condition} if for some constant $K\geq 2$ it holds
\begin{equation}\label{delta_2}
   \varphi(2t)\leq K\varphi(t)\qquad \forall\,t\geq0\,.
\end{equation}
We write $\varphi\in \Delta_{2}$ if $\varphi$
satisfies the $\Delta_{2}$-condition. The $\Delta_2$-constant of
$\varphi$, denoted by $\Delta_2(\varphi)$, is the smallest constant
$K\geq 2$
satisfying \eqref{delta_2}. 
\end{definition}
We have the following results.
\begin{lemma}\label{lem:basic}
  For a regular N-function $\varphi$ the following properties are
  satisfied: 
  \begin{itemize}
  \item [{\rm (i)}] For all $t \ge 0$ there holds
    \begin{align*}
      \varphi(t) \le \varphi'(t)t \le \varphi(2t)\,.
    \end{align*}
  \item [{\rm (ii)}] If $\varphi \in \Delta_2$, then we have for all $t \ge 0$ 
    \begin{equation*}
      \varphi(t) \le \varphi'(t)t \le \Delta_2(\varphi)\,\varphi(t)\,.
    \end{equation*}
  \item [{\rm (iii)}] It holds that $\varphi \in \Delta{_2}$ if and only if $\varphi' \in
    \Delta_2$. In this situation we have $\Delta_2(\varphi)\leq 2 \Delta_2(\varphi') \leq (\Delta_2(\varphi))^{2}$.
  \end{itemize}
\end{lemma}
\begin{proof}
  Assertion (i) is contained in~\cite[Lemma~5.1]{dr-nafsa}.
  Assertion (ii) follows from (i).   Assertion (iii) is proved in~\cite[Lemma~5.2]{dr-nafsa}.  
\end{proof}

For a regular N-function $\varphi$ we define the \textit{complementary function}
$\varphi^*$ by
\begin{equation*}
  \varphi^*(t):=\int\limits_0^t(\varphi')^{-1}(s)\,ds\,.
\end{equation*}
It is easily seen from this definition, using elementary properties of
inverse functions (cf.~proof of \cite[Lemma 6.4]{dr-nafsa}), that
$\varphi^*$ is again a regular N-function.  We have the following
versions of Young inequality.
\begin{lemma}[Young type inequalities]\label{lem:young}
 Let the regular N-function $\varphi$ be such that $\varphi, \varphi^* \in
 \Delta_2$. Then, for all $t,u\geq 0$ there holds
\begin{align*}
  t u&\leq \vep\,\varphi(t)+(\Delta_2(\varphi^*))^M\,\varphi^*(u)\,,
       \\
  t u&\leq \vep\,\varphi^*(t)+(\Delta_2(\varphi))^M\,\varphi(u)\,,
       \\
  t \varphi'(u)&\leq \vep\,\varphi(t)+\Delta_2(\varphi)\,(\Delta_2(\varphi^*))^M\,\varphi^*(u)\,,
       \\
  \varphi'(t) u&\leq \vep\,\varphi^*(t)+(\Delta_2(\varphi))^N\,\varphi(u)
\end{align*}
for all $\vep \in (0,1)$, $M \in \setN$ such that $\vep^{-1}
\le 2^M$, and $N \in \setN$ such that $\Delta_2(\varphi)\,\vep^{-1}
\le 2^N$.
\end{lemma}
\begin{proof}
  The first two inequalities follow immediately from the classical
  Young inequality 
  \begin{equation*}
      t u\leq \varphi(t)+\varphi^*(u)\,,
  \end{equation*}
  $\varphi, \varphi^* \in \Delta_2$, and
  $\psi(\vep \,t)\le \vep\, \psi(t)$, valid for all convex functions
  $\psi$, $t\ge 0$ and $\vep\in (0,1)$. The last two inequalities
  follow from the first ones and the equivalence
  \begin{align}
    \label{eq:*'}
    (\Delta_2(\varphi^*))^{-1}\varphi(t)\le \varphi^*(\varphi'(t))\le \Delta_2(\varphi)\,\varphi(t)  \,,
  \end{align}
  valid for all $t \ge 0$ (cf.~\cite[(5.17)]{dr-nafsa}).
\end{proof} 
%
In the study of nonlinear problems like \eqref{eq:pfluid} and of N-functions the
property~\eqref{eq:phi_pp} below plays a fundamental role. To keep the presentation
shorter we call functions satisfying it  ``balanced function''.
\begin{definition}[Balanced function]
  \label{ass:phi}
  We call a regular N-function $\varphi$ {\em balanced}, if  there exist
  constants $\gamma_1\in (0,1]$ and  $\gamma_2 \ge 1$ such that for all $t> 0$
  there holds
  \begin{align}
    \label{eq:phi_pp}
    \gamma_1\,\varphi'(t)\le t\,\varphi''(t)\le \gamma_2\,\varphi'(t)
    \,.
  \end{align}
  The constants $\gamma_1$ and $\gamma_2$ are called {\em characteristics}
  of the balanced  N-function $\varphi$, and  will be denoted as $(\gamma_1,\gamma_2)$. 
\end{definition}
This property transmits itself to complementary functions. 
\begin{lemma}
  \label{lem:phi*}
  Let $\varphi$ be a balanced N-function with characteristics
  $(\gamma_1,\gamma_2)$. 
  Then, the complementary  N-function $\varphi^\ast$ is a
  balanced {N-function} with characteristics
  $(\gamma_2^{-1},\gamma_1^{-1})$. 
\end{lemma}
\begin{proof}
  The assertion is proved in \cite[Lemma 6.4]{dr-nafsa}. The proof
  uses only the condition~\eqref{eq:phi_pp}, and the formula for the
  derivative of the inverse function applied to 
  $(\varphi^\ast)'(t) = (\varphi')^{-1}(t)$.
\end{proof}
Balanced N-functions always satisfy the $\Delta_2$-condition (cf.~\cite{BDK12}).
\begin{lemma}\label{lem:d2}
   For a balanced N-function $\varphi$ we have that 
   $\varphi, \varphi^\ast \in \Delta_2$.  In
  particular, for all $t \ge  0$ there holds
  \begin{align*}
    \begin{aligned}
      \varphi(2t)&\le 2^{\gamma_2+1}\,\varphi(t)\,,
      \\
      \varphi^*(2t)&\le 2^{\frac 1{\gamma_1}+1}\,\varphi^*(t)\,,
    \end{aligned}
  \end{align*}
  i.e., the $\Delta_2$-constants of $\varphi$ and $\varphi^*$ possess an
  upper bound depending only on
  $\gamma_1$ and $\gamma_2$.
\end{lemma}
\begin{proof}
  From condition \eqref{eq:phi_pp} it follows for all $t>0$ that
\begin{equation*}
\frac{d}{dt}\log(\varphi'(t))=  \frac{\varphi''(t)}{\varphi'(t)}\leq \gamma_2\frac{1}{t} \,,
\end{equation*}
which implies by integration with respect to $t$ over $(s,2s)$, $s>0$, and using the
exponential function that
\begin{align*}
  \frac{\varphi'(2s)}{\varphi'(s)} \le 2^{\gamma_2}\,.
\end{align*}
A further integration with respect to $s$ over $(0,t)$, $t>0$,  proves, for all $t>0$, that 
\begin{equation*}
  \varphi(2t)\leq2^{\gamma_2+1}\varphi(t)\,,
\end{equation*}
showing the assertion for $\varphi$. The assertion for $\varphi^*$ follows
analogously by using Lemma \ref{lem:phi*}. 
\end{proof}
\begin{corollary}\label{cor:sim}
  For a balanced N-function  $\varphi$ we have
  \begin{align*}
    \varphi(t)\sim\varphi'(t)\,t \sim  \varphi''(t)\,t^{2}\qquad
    \text {for all } t>0\,, 
  \end{align*}
  with constants of equivalence depending only on the characteristics
  of $\varphi$. 
\end{corollary}
\begin{proof}
  This follows immediately from Lemma \ref{lem:basic} and Lemma
  \ref{lem:d2} since $\varphi$ is balanced.
\end{proof}
\begin{lemma}\label{lem:trans}
  Let  $\varphi$ a balanced N-function with characteristics
  $(\gamma_1,\gamma_2)$. Let $\potV \in C^1(\setR^{\ge 0})\cap C^2(\setR^{>0})$ with
  $\potV(0)=\potV'(0)=0$ satisfy for some $c_{0},c_{1}>0$ and for all $t > 0$
  \begin{align}
    \label{eq:equi1}
    c_0 \,\varphi ''(t) \le \potV''(t) \le c_1\, \varphi''(t)\,.
  \end{align}
  Then, also $\potV$ is a balanced N-function with
  characteristics $(\gamma_2 \frac {c_0}{ c_1}, \gamma_1\frac{c_1}{c_0})$,
  which satisfies for all $t\ge 0$ 
  \begin{equation}
    \label{eq:equi2}
    \begin{gathered}
      c_0 \,\varphi '(t) \le \potV'(t) \le c_1\, \varphi'(t)\,,
      \\
      c_0 \,\varphi (t) \le \potV(t) \le c_1\, \varphi(t)\,.
    \end{gathered}
  \end{equation}
\end{lemma}
\begin{proof}
  The inequalities in \eqref{eq:equi2} follow from \eqref{eq:equi1} by
  integration using that $\potV'(0)=\varphi'(0)=0$. From \eqref{eq:equi2}
  and \eqref{eq:phi_pp} it follows that $\potV$ is a balanced N-function with
  characteristics as indicated in the assertion.
\end{proof}

It turns out that the function $\MC_\varphi: \setR^{\ge 0} \to
\setR^{\ge 0}$, defined for regular N-functions $\varphi$ via 
\begin{align}\label{eq:mot}
  \MC_{\varphi}(t):=\frac{\varphi'(t)}{t} \,, 
\end{align}
plays an important role in the investigation of problem
\eqref{eq:pfluid}. 
\begin{lemma}
  \label{lem:giusti1}
  Let $\varphi $ be a regular N--function such that
  $\varphi, \varphi^* \in \Delta_2$.  Then, for all
  $\bfP, \bfQ \in \setR^{3 \times 3}$ there holds
  \begin{align*}
    \MC_\varphi(\abs{\bfP} + \abs{\bP-\bfQ})\sim \int\limits_0^1 \MC_{\varphi}( \abs{\theta\,\bfP
    +(1-\theta)\bfQ} )\, d\theta\,, 
  \end{align*}
  with constants of equivalence depending only on $\Delta_2(\varphi)$ and
  $\Delta_2(\varphi^*)$. 
\end{lemma}
\begin{proof}
  This follows immediately from \cite[Lemma 6.6]{dr-nafsa} by using
  Lemma \ref{lem:basic}, the convexity of $\varphi$, 
  $\varphi \in \Delta_2$, and
  $2^{-1} (\abs{\bfP} + \abs{\bfQ}) \le \abs{\bfP} + \abs{\bP-\bfQ}
  \le 2(\abs{\bfP} + \abs{\bfQ})$.
\end{proof}

It is convenient to introduce for all $p\in(1,\infty)$ and all
$ \para\in [0,\infty)$ the function
$\function_{p,\para}:\setR^{\ge 0} \to \setR^{\ge 0}$ via
\begin{equation*}
 \function(t)=\function_{p,\para}(t):=\int\limits_{0}^{t}(\para+s)^{p-2}s\,ds
 \qquad\forall\,t \ge 0\,, 
\end{equation*}
which is precisely the N-function associated with the definition of the tensor
$\mathbf{S}$ from~\eqref{eq:example}. If $p$ and $\delta$ are fixed
we often simply write
$\function(t):=\function_{p,\para}(t)$. Nevertheless, we will track
the possible dependence of constants in terms of these two parameters.
Clearly, $\function _{p,\para}$ is a regular N-function for all
$p\in(1,\infty)$ and all $ \para\in [0,\infty)$. The advantage is that
we have exact control of all relevant constants for these functions.
We have the following basic properties.
\begin{lemma}\label{lem:function}
  For any $\delta\in[0,\infty)$ and for any $p\in(1,\infty)$ there
  holds 
  \begin{equation}\label{eq:E}
    \begin{aligned}
     \function_{p,\para}(t)&\leq (\function_{p,\para})'(t)\,t\;\leq
      2^{p+1}\function_{p,\para}(t)\qquad \forall\,t\geq 0\,,
      \\
      \min\{1,p-1\}  \,    (\function_{p,\para})'(t)&\leq (\function_{p,\para})''(t)\,t\leq
      \max\{1,p-1\}\,(\function_{p,\para})'(t) \quad \forall\,t >0\,.\hspace*{-5mm}
    \end{aligned}
  \end{equation}
  In particular, $\function_{p,\para}$, $p\in (1,\infty)$, $\delta\ge 0$, are balanced
  N-functions with characteristics $( \min\{1,p-1\} , \max\{1,p-1\})$ and
  $\Delta_2$-constants depending only on $p$. Moreover, by the previous results also
  $(\function_{p,\para})^*$ are balanced N-functions with
  characteristics $(\min\{1,(p-1)^{-1}\}, \max\{1,(p-1)^{-1}\})$ and
  $\Delta_2$-constants depending only on $p$.
  \end{lemma}
  \begin{proof}
  The first assertion in \eqref{eq:E} follows from Lemma \ref{lem:basic} (ii) and
  \cite[Lemma~5.3]{dr-nafsa}, since $\Delta_2((\function_{p,0})')=2^{p-1}$. The
  second assertion \eqref{eq:E} follows from direct computations.
\end{proof}

\subsection{Nonlinear operators with $(p,\delta)$-structure}
In this section we collect the main results on nonlinear operators derived from a 
potential and having \mbox{$(p,\delta)$-struct-} ure. 
\begin{definition}[Operator derived from a potential]
  \label{def:potential}
  We say that an operator  \linebreak ${\bS:\,\setR^{3\times3}\to\setR^{3\times
    3}_{\sym}}$ is {\em derived from a potential} 
  $\pot: \setR^{\ge 0} \to \setR^{\ge 0}$, if $\bS(\bfzero)=\bfzero$
  and for all $\bP\in\setR^{3\times3}\setminus \set{\bfzero}$ there
  holds\footnote{Here we use the notation \eqref{eq:mot} also for a
    more general function $\pot$ (not necessarily a balanced or even a
    regular N-function).  } 
  \begin{align*}
    \bS(\bP)=\partial
    U(|\bP^{\sym}|)=\frac{\pot'(|\bP^{\sym}|)}{\abs{\bP^{\sym}}}\,\bP^{\sym}
    =\MC_{U}(\abs{\bP^{\sym}})\,\bP^{\sym}
  \end{align*}
  for some $\pot\in  C^{1}(\setR^{\geq0})\cap C^{2}(\setR^{>0})$
  satisfying  $\pot(0)=\pot'(0)=0$.
\end{definition}
\begin{remark}
  For ease of notation, in many cases we will also write
  $\bS=\partial \pot$ for an operator derived from the potential $\pot$ and
  note that from its definition it follows that
  $\bS(\bP)=\bS(\bP^{\sym})$, for all $\bP \in \setR^{3\times3}$.

  Note also that we consider the operator $\bS$ with domain
  $\setR^{3\times3}$, since we study the problem~\eqref{eq:pfluid} in
  the setting of three space-dimensions. Clearly, the same definition
  and results below can be applied to a general operator defined on
  $\setR^{d\times d}$, with $d\geq2$.
\end{remark}
\begin{remark}
  Note that in investigations of the regularity of solution of
  \eqref{eq:pfluid}, or its steady analogues, for operators derived from
  a potential $\pot$, the lower and upper bounds of the quantity
  \begin{align*}
    \frac {(\MC_\pot)'(t)\, t}{\MC_\pot(t)}= \frac{\pot''(t)\, t}{\pot'(t)}-1\,, 
  \end{align*}
  play an important role (cf.~the discussion in
  \cite{BM2020, CM2011, CM2019,CM2020}. 
  If $\pot$ is a balanced N-function
  these bounds are closely related to the characteristics
  $(\gamma_{1},\gamma_{2})$ of $\pot$. In fact, we have
  \begin{align*}
     \gamma_1-1\le \frac {(\MC_\pot)'(t)\,t}{\MC_\pot(t)}\le \gamma_2-1\,.
  \end{align*}
\end{remark}

\begin{definition}[Operator with  ${\varphi}$-structure]
  \label{def:ass_S}
  Let the operator
  ${\bS\colon \setR^{3 \times 3} \to \setR^{3 \times 3}_{\sym} }$,
  belonging to
  $C^0(\setR^{3 \times 3};\setR^{3 \times 3}_{\sym} )\cap C^1(\setR^{3
    \times 3}\setminus \{\bfzero\}; \setR^{3 \times 3}_{\sym} ) $,
  satisfy ${\bS(\bP) = \bS\big (\bP^{\sym} \big )}$ and
  $\bS(\mathbf 0)=\mathbf 0$. 
  We say that $\bS$ has {\em ${\varphi}$-structure} if there exist 
  a regular N-function $\varphi$ and constants
  $\gamma_3 \in (0,1]$, $\gamma_4 >1$ such that the inequalities
   \begin{subequations}
     \label{eq:ass_S}
     \begin{align}
       \sum\limits_{i,j,k,l=1}^3 \partial_{kl} S_{ij} (\bP) Q_{ij}Q_{kl}
       &\ge \gamma_3 \, \MC_{\varphi}(\abs{\bfP^{\sym}})\, |\bP^{\sym} |^2\,,\label{1.4b} 
       \\
       \big |\partial_{kl} S_{ij}({\bP})\big |
       &\le \gamma_4  \, \MC_{\varphi}(\abs{\bfP^{\sym}})\,,  \label{1.5b}
     \end{align}
   \end{subequations}
   are satisfied for all $\bP,\bQ \in \setR^{3\times 3} $ with $\bP^{\sym} \neq \bfzero$ and
   all $i,j,k,l=1,2, 3$.  The constants $\gamma_3$, $\gamma_4$, and $\Delta_2(\varphi)$
   are called the {\em characteristics} of $\bfS$ and will be denoted
   by $(\gamma_3,\gamma_4, \Delta_2(\varphi))$. 

   In the special case $\varphi= \function_{p,\delta} $ with $p \in (1, \infty)$ and
   $\para\in [0,\infty)$ we say that $\bS$ has {\em $(p,\para)$-structure} and call
   $(\gamma_3,\gamma_4,p)$ its characteristics.
\end{definition}

Closely related to an operator with $\varphi$-structure is 
the function
${\bF_\varphi\colon\setR^{3 \times 3} \to \setR^{3 \times 3}_{\sym}}$ defined via
\begin{equation}\label{def:F}
  \bF_\varphi(\bP):=\sqrt{\MC_{\varphi}(\abs{\bfP^{\sym}})}\,\bP^{\sym}
  = \frac{\sqrt{\varphi'(\abs{\bfP^{\sym}})\abs{\bP^{\sym}}}}{\abs{\bfP^{\sym}}}\,\bP^{\sym}\,,
\end{equation}
where the second representation holds only for $\bP^{\sym} \neq \bfzero$. However,
this form is convenient since it shows that $\bF_\varphi$ is derived from
the potential 
\begin{align}
  \psi(t):=\int\limits_0^t\sqrt{\varphi'(s)s}\, ds\,.\label{eq:psi}
\end{align}
It is easily seen that $\psi \in C^1(\setR^{\ge 0})\cap
C^2(\setR^{>0})$. 
In the special case of an operator $\bS$ with
$(p,\para)$-structure  we have with $\function=\function_{p,\para}$ 
\begin{equation}
  \label{def:F1}
  \bF(\bP):=\bF_{\function}(\bP)=\sqrt{\MC_{\function}(\abs{\bfP^{\sym}})}\,\bP^{\sym}
  = \big (\para +\abs{\bP^{\sym}}\big )^{\frac {p-2}2}\bP^{\sym}\,,
\end{equation}
which is consistent with the notation used in the previous literature, as
explained in the introduction, cf.~\eqref{eq:F}.

To derive a very important result for operators with
$\varphi$-structure we need the following result, which explains also
the link (and the choice of a similar name) between the characteristics
of a balanced N-function $\varphi$, and the characteristics of an
operator derived from a potential $\varphi$.
\begin{proposition}\label{prop:pot-phi}
  Let $\varphi$ be a balanced N-function with
  characteristics $(\gamma_1,\gamma_2)$.  Let $\bT=\partial\varphi$ be
  derived from the potential $\varphi$.  Then, $\bT$ has $\varphi$-structure with
  characteristics depending only on $\gamma_1$ and $\gamma_2$. 
\end{proposition}
\begin{proof}
  It follows from Lemma \ref{lem:d2} that the $\Delta_2$-constant of
  $\varphi$ depends only on $\gamma{_2}$. We have for all $\bP \in \setR^{3\times 3}$ with $\bP^{\sym} \neq {\bfzero}$
 \begin{align}\label{eq:pSphi}
  \partial_{kl} T_{ij} (\bP)& = \frac {\varphi'(|\bP^{\sym} |)}{|\bP^{\sym} |}
  \left(\delta^{\sym}_{ij,kl}-\frac{P^{\sym}
      _{ij}P^{\sym} _{kl}}{|\bP^{\sym} |^2}\right) + \varphi''(|\bP^{\sym} |) 
  \,\frac{P^{\sym} _{ij}P^{\sym} _{kl}}{|\bP^{\sym} |^2}\,, 
 \end{align}
 where $\delta^{\sym}_{ij,kl}:= \frac 12 (\delta_{ik}\delta_{j l}+\delta_{i
   l}\delta_{j k})$.
   Using \eqref{eq:phi_pp} we obtain from this, for all $j,k,l,m$,
  \begin{align*}
    \bigabs{\partial_{kl} {T_{ij}(\bfP)}} \leq 2\,
    \frac{\varphi'(\abs{\bfP^{\sym}})}{\abs{\bfP^{\sym}}} + \varphi''(\abs{\bfP^{\sym}})
    \leq (2+\gamma_2)\, \MC_{\varphi}(\abs{\bfP^{\sym}})\,,
  \end{align*}
  which proves \eqref{1.5b}.  From \eqref{eq:pSphi},
  \eqref{eq:phi_pp}, and $\gamma_1\le 1$ we obtain for
  $\bfP, \bfQ \in \setR^{3 \times 3}$ with $\bfP^{\sym} \not= \bfzero$
  \begin{align*}
    &\sum\limits_{i,j,k,l=1}^3 \partial_{kl} T_{ij} (\bP)
      Q_{ij}Q_{kl}
    \\
    &= \frac {\varphi'(|\bP^{\sym} |)}{|\bP^{\sym} |}
      \left(\abs{\bQ^{\sym}}^2-\frac{\abs{\bP^{\sym}\cdot
      \bQ^{\sym}}^2}{|\bP^{\sym} |^2}\right) + \varphi''(|\bP^{\sym} |)  
      \,\frac{\abs{\bP^{\sym}\cdot \bQ^{\sym}}^2 }{|\bP^{\sym} |^2}
    \\
    &\geq \gamma_1\, \frac {\varphi'(|\bP^{\sym} |)}{|\bP^{\sym} |}
      \left(\abs{\bQ^{\sym}}^2-\frac{\abs{\bP^{\sym}\cdot
      \bQ^{\sym}}^2}{|\bP^{\sym} |^2}\right) + \gamma_1\, \frac
      {\varphi'(|\bP^{\sym} |)}{|\bP^{\sym} |}  
      \,\frac{\abs{\bP^{\sym}\cdot \bQ^{\sym}}^2 }{|\bP^{\sym} |^2}
    \\
    &= \gamma_1\, \MC_ {\varphi}(|\bP^{\sym} |)\,\abs{\bQ^{\sym}}^2\,,
  \end{align*}
  which proves \eqref{1.4b}.
\end{proof}
We can now formulate the following crucial result for our
investigations (cf.~\cite[Section~6]{dr-nafsa}).

\begin{proposition}\label{prop:hammer-phi}
  Let
  $\varphi$ be a balanced N-function with characteristics
  $(\gamma{_1},\gamma{_2})$. Let $\bS$ have
  $\varphi$-structure with characteristics $(\gamma_3,\gamma{_4},
  \Delta{_2(\varphi)})$ and let
  $\bF_\varphi$ be defined in \eqref{def:F}. Then, we have for all
  $\bP,\bQ \in \setR^{3\times 3} $ that
  \begin{align}
    \big(\bS(\bP)-\bS(\bQ)\big)\cdot(\bP-\bQ)
    &\sim \MC_{\varphi}(\abs{\bfP^{\sym}} +
      \abs{\bfP^{\sym}-\bQ^{\sym}})\,|\bP^{\sym}-\bQ^{\sym}|^{2}  \label{eq:ham-ap}
    \\
    &      \sim |\bF_\varphi(\bP)-\bF_\varphi(\bQ)|^{2}\,, \label{eq:ham-bp}
    \\
    |   \bS(\bP)-\bS(\bQ)|&\sim \MC_{\varphi}(\abs{\bfP^{\sym}}+\abs{\bfP^{\sym}-\bfQ^{\sym}} 
       )\, |\bP^{\sym}-\bQ^{\sym}|\,, \label{eq:ham-cp}
\end{align}
where the constants of equivalence depend only on 
$\gamma_{1}, \gamma{_2}, \gamma{_3}$, and $\gamma_{4}$.
\end{proposition}
\begin{proof}
  First of all note that, due to Lemma \ref{lem:d2} and Lemma
  \ref{lem:phi*}, $\varphi$ and $\varphi^*$ satisfy the
  $\Delta_2$-condition with $\Delta{_2}$-constants depending only on
  $\gamma{_1}$ and $\gamma{_2}$.
  Using \eqref{eq:ass_S} and Lemma \ref{lem:giusti1} 
  we get that for all $\bP, \bQ \in \setR^{3 \times 3}$ with
  $\bP^{\sym} \neq \bfzero$
\begin{align*}
  &\big (\bS (\bP) - \bS (\bQ)\big ) \cdot(\bP -\bQ)\\
  & = \int\limits^1_0 \sum\limits_{i,j,k,l=1}^3 \partial_{kl} S_{ij} \big ( \theta \bfP
  +(1-\theta)\bfQ\big ) (P
    -Q)_{ij} (P -Q)_{kl} \, d\theta
  \\
  &\sim  \int\limits_0^1 \MC_{\varphi}( \abs{\theta \bfP^{\sym}
  +(1-\theta)\bfQ^{\sym}} )\, d\theta \;| \bP^{\sym} - \bQ ^{\sym}|^2 
\\
  & \sim \MC_{\varphi}( \abs{\bfP^{\sym}}
  +\abs{\bP^{\sym}-\bfQ^{\sym}} )\,| \bP^{\sym} - \bQ ^{\sym}|^2\,,
\end{align*}
which proves \eqref{eq:ham-ap} with constants of equivalence depending only on 
$\gamma_{1},\gamma_2,\gamma{_3}$, and $\gamma_{4}$. From \eqref{eq:ham-ap} we immediately
obtain, also using that $\bS$ is symmetric,  
\begin{align*}
  \MC_{\varphi}( \abs{\bfP^{\sym}}
  +\abs{\bP^{\sym}-\bfQ^{\sym}} )\,| \bP^{\sym} - \bQ ^{\sym}|^2
  &\le c \, \big (\bS (\bP) - \bS (\bQ)\big ) \cdot
    (\bP -\bQ)
  \\
  &\le c\, \abs  {\bfS(\bfP)-\bfS(\bfQ)} \abs {\bP^{\sym}-\bQ^{\sym}}\,,
\end{align*}
with constants depending only on 
$\gamma_{1},\gamma_2,\gamma{_3}$, and $\gamma_{4}$. This proves one inequality in \eqref{eq:ham-cp}. The other follows from 
\begin{align*}
  \abs {\bfS(\bfP)-\bfS(\bfQ)}
  &= \bigg( \sum_{i,j=1}^3\Big (\sum_{k,l=1}^3 \int\limits^1_0 \partial _{kl} S_{ij}
  \big ( \theta \bfP
  +(1-\theta)\bfQ\big ) (P -Q)_{kl} \, d\theta\Big )^2\bigg)^{\frac
    12}
  \\
  &\le c\, \int\limits_0^1 \MC_{\varphi}( \abs{\theta \bfP^{\sym}
  +(1-\theta)\bfQ^{\sym}} )\, d\theta \;| \bP^{\sym} - \bQ ^{\sym}|
\\
  & \le c\, \MC_{\varphi}( \abs{\bfP^{\sym}} +\abs{\bP^{\sym}-\bfQ^{\sym}} )\, | \bP^{\sym} - \bQ ^{\sym}|\,,
\end{align*}
with constants depending only on $\gamma_{1}, \gamma_2, \gamma{_3}$,
and $\gamma_{4}$. Here, we used again \eqref{1.5b}, the symmetry of
$\partial _{kl}S_{ij}$ with respect to $k,l$, and Lemma
\ref{lem:giusti1}.

To show \eqref{eq:ham-bp} we use that $\bF_\varphi$ defined in
\eqref{def:F} possesses $\psi$-structure, where $\psi$ is defined in
\eqref{eq:psi}. We have using \eqref{eq:phi_pp}, for all $t>0$, that
\begin{align*}
  \psi''(t)\, t &=\frac{\big(\varphi''(t) \,t +\varphi'(t)\big )\,t}{ 2\,
                  \sqrt{\varphi'(t)\, t}}
                  \sim 
                  \frac{\varphi'(t)\,t}{ 
                  \sqrt{\varphi'(t)\, t}} =  \psi'(t) \,.
\end{align*}
This shows that
$\psi$ is a balanced N-function with characteristics $(\frac{1+\gamma _1}{2}, \frac{1+\gamma _2}{2})$. Thus, Proposition
\ref{prop:pot-phi} yields that $\bF_\varphi$ has $\psi$-structure with
characteristics depending only on $\gamma_1$ and $\gamma_2$. The already
proven equivalence \eqref{eq:ham-cp} reads in this case as
\begin{align}
  \label{eq:Fequi}
  |\bF_\varphi(\bP)-\bF_\varphi(\bQ)|^{2}
  &\sim \big (\MC_{\psi}(\abs{\bfP^{\sym}} +
    \abs{\bP^{\sym}-\bfQ^{\sym}})\big)^2  \,|\bP^{\sym}-\bQ^{\sym}|^2\,,
\end{align}
with constants of equivalence depending only on 
$\gamma_{1}$ and $\gamma_{2}$. From the definition of $\psi$ we get, for all $t>0$,
\begin{align*}
  \big (\MC_{\psi}(t)\big)^2
  &= \bigg  (\frac{\sqrt{\varphi'(t)\,t}}{t}\bigg)^2 = \MC_{\varphi}(t)\,,
\end{align*}
which together with \eqref{eq:Fequi} yields \eqref{eq:ham-bp}. This
finishes the proof. 
\end{proof}

Let us finish this section by proving a useful result for the operator
occurring in~\eqref{eq:pfluid}. 
%
\begin{proposition}
\label{prop:potential-equivalence}
Let the operator $\bT=\partial\potV$, derived from a potential
$\potV$, have $\varphi$-structure, with characteristics
$(\gamma_3,\gamma_4,\Delta_2(\varphi))$. If $\varphi$ is a balanced
N-function with characteristics $(\gamma_1,\gamma_2)$, then $\potV$ is
a balanced N-function satisfying for all $t>0$
  \begin{equation}\label{eq:simU}
     \frac{ \gamma_3 }{\gamma_2}\varphi''(t)\le \potV''(t)\le
     \frac{ \gamma_4 }{\gamma_1}\varphi''(t)\,.
  \end{equation}
  The characteristics of $\potV$ is equal to $\big(\frac{\gamma_3}{\gamma_4}
    \,\frac{ \gamma_1^2}{\gamma_2}, \frac{\gamma_4}{\gamma_3}
    \,\frac{ \gamma_2^2}{\gamma_1}\big )$.
\end{proposition}
\begin{proof} 
  For $\bP =\frac{t}{\sqrt 3}\bI\bd$, $t>0$,
  $\bQ=\frac{1}{\sqrt 3}\bI\bd$ we get $\abs{\bP} =t$, $\abs{\bQ}=1$.
  Thus, \eqref{1.4b}, \eqref{1.5b}, and the definition of $\MC_\varphi$
  yield
  \begin{align*}
   \gamma_3 \frac{\varphi'(t)}{t}\le \sum\limits_{i,j,k,l=1}^3 \partial_{kl} T_{ij} (\bP)
   Q_{ij}Q_{kl} =\potV ''(t)\le \gamma_4 \frac{\varphi'(t)}{t}\,.
 \end{align*}
 This implies \eqref{eq:simU}, since $\varphi$ is balanced. The remaining
 assertions follow from Lemma \ref{lem:trans}.
\end{proof}

\begin{remark}\label{rem:pot}
  Proposition \ref{prop:potential-equivalence} states that $\potV$ is a
  balanced N-function with characteristics depending only on the
  characteristics of $\bS$ and on the characteristics of
  $\varphi$. Consequently, Lemma~\ref{lem:d2} and Lemma \ref{lem:phi*}
  yield that $\potV $ and $\potV ^*$ satisfy the $\Delta_{2}$-condition,
  with $\Delta_2$-constants depending only on the characteristics of
  $\bS$ and the characteristics of $\varphi$.
\end{remark}

\subsection{Approximations of a nonlinear operator}\label{sec:approx}

We now define the {$(A,q)$-ap\-}proximation and prove the relevant properties,
needed in the sequel. Note that the $(A,q)$-approximation in the special
case $p\ge 2$ and $q=2$ was introduced in a slightly different form
in~\cite{mnr3} (in that reference the potential depends on
$|\bP^{\sym}|^{2}$). 
The idea behind is that the operator induced by the $(A,q)$-approximation
for $q=2$
has linear growth at infinity (cf.~\cite[Lemma~2.22]{mnr3}) and
consequently, one can work on the level of this $(A,2)$-approximation
within the standard Hilbertian theory.
\begin{definition}[$(A,q)$-approximation of a
  scalar real function]\label{def:A-approx}
  Given a function $\potV \in  C^{1}(\setR^{\geq0})\cap C^{2}(\setR^{>0})$
  satisfying  $\potV(0)=\potV'(0)=0$  we
  define for $A\ge 1$ and $q\ge  2$ the {\em $(A,q)$-approximation}
  $\appAq{\potV}  \in   C^{1}(\setR^{\geq0})\cap C^{2}(\setR^{>0}) $ via
  \begin{equation*}
    \appAq{\potV}(t):=\left\{
      \begin{aligned}
        &\potV(t)\qquad &t\leq A\,,
        \\
        &\alpha_{2,q}\,t^{q}+\alpha_{1,q}\,t +\alpha_{0,q}\qquad &t> A\,.
      \end{aligned}
    \right.
  \end{equation*}
  Consequently, the constants $\alpha_{i,q}=\alpha_{i,q}(\potV)$, $i=0,1,2$, are given by 
  \begin{equation*}
    \begin{aligned}
      \alpha_{2,q}&=\frac{1}{q(q-1)}\,\frac{\potV''(A)}{A^{q-2}}\,,
      \\
      \alpha_{1,q}&=\potV'(A)-\frac{1}{q-1}\,\potV''(A)\,A\,,
      \\
      \alpha_{0,q}&=\potV(A)-\potV'(A)\,A+\frac{1}{q}\potV''(A)\,A^{2}\,.
    \end{aligned}
  \end{equation*}
\end{definition}

\begin{remark}\label{rem:VA}
  If  $\varphi$ is a regular N-function and $q=2$, the definition of the
  $(A,2)$-approximation $\varphi^{A,2}$ and the properties of $\varphi$ immediately imply
  that there exists a constant $c(A,\varphi)$ such that for all $t\ge 0$ there holds
  \begin{align*}
    a_{\varphi^{A,2}}(t) = \frac{(\varphi^{A,2})'(t)}{t} \le c(A,\varphi)\,.
  \end{align*}
\end{remark}
Next, we define the $(A,q)$-approximation of an operator derived
from a potential.

\begin{definition}[$(A,q)$-approximation of an operator derived from a
  potential]\label{def:SA}
  Let the operator $\bS=\partial\pot$ be derived from the potential
  $\pot$.  Then, we define for given $A\ge 1$ and $q\ge 2$ the {\em
    $(A,q)$-approximation} $\appAq{\bS}:=\partial \appAq{\pot}$ of
  $\bS$ as the operator derived from the potential $\appAq{\pot}$,
  i.e., $\appAq\bS$ satisfies $\appAq\bS(\bfzero)=\bfzero$ and for all
  ${\bP\in \setR^{3\times 3}\setminus\set{\bfzero}}$ there holds
  \begin{equation*}
    \appAq\bS(\bP):=\partial\appAq{\pot}(|\bP^{\sym}|)
    =\frac{(\appAq\pot)'(|\bP^{\sym}|)}{\abs{{\bP^{\sym}}}}\,\bP^{\sym}
    = \MC_{\appAq{\pot}}(\abs{\bP^{\sym}})\,\bP^{\sym}
    \,.
\end{equation*}
\end{definition}
As explained in the introduction, for an operator with
$(p,\para)$-structure, for large $p$, we need also multiple approximations,
which we define now.
\begin{definition}[Multiple approximation of an operator]\label{def:mult}
  Let the operator $\bS$ have $(p,\delta)$-structure for some
  $p \in (2, \infty)$ and $\para\in [0,\infty)$ and let $\bS$ be derived from the
  potential $\pot$.  For given $N\in \setN$ and $q_n \in [2,p]$,
  $n=0,\ldots, N$ with $q_0=p$, $q_{N}=2$ and $q_n>q_{n+1}$,
  $n=0, \ldots ,N-1$, and $A_n\ge 1$, $n=1,\ldots, N$ with
  $A_{n+1}\ge A_n+1$, $n=1, \ldots ,N-1$, we set
\begin{equation*}
\pot^{0}:=\pot,\quad\bS^0:=\bS,\quad\function^0:=\function
_{p,\para},\quad\bF^{0} :=\bF_{\function^0},\quad a^{0}:=a_{\function^0},
\end{equation*}
and then recursively
\begin{equation*}
  \begin{aligned}
    \pot^{n}:=(\pot^{n-1})^{A_{n},q_{n}},\quad \!\bS^n:=\partial
    U^n,\quad \!  \function^{n}:= (\function^{n-1})^{A_{n},q_{n}},
    \quad \!  \bF^n:=\bF_{\function^n}, \quad \!\MC^n:=
    \MC_{\function^n},
  \end{aligned}
  \end{equation*}
  for $n=1,\ldots, N$.  We call $\pot^n$, $\bS^n$, $\function^n$,
  $\bF^n$, and $\MC^n$, $n=1,\ldots, N$, {\em multiple approximation}
  of $\pot$, $\bS$, $\function_{p,\delta}$, $\bF$, and $\MC$,
  respectively.
\end{definition}
\begin{remark}
  \label{rem:single}
  As we will see later on (for the parabolic problem in three-space
  dimensions) strictly speaking the multiple approximation is not
  needed for $p\in (1,\frac {13}3)$.  Since in the definition of a multiple
  approximation the case ${N=1}$ is included, also a single
  $(A,q)$-approximation is a special case of a multiple
  approximation. To unify the presentation we also call the
  $(A,2)$-approximation for $p\in (1,\frac {13}3)$ multiple approximation. In this
  case we have $\pot^1=\pot^{A,2}$, $\bS^1=\bS^{A,2}=\partial
  U^1$, 
  $\function^1= (\function^{0})^{A,2}=(\function _{p,\para})^{A,2}$,
  $\bF^1=\bF_{\function^1}$, and $\MC^1= \MC_{\function^1}$.
\end{remark}

In the following we derive various properties of multiple
approximations for an operator $\bS$ which is derived from  a potential $\pot$ and has
$(p,\delta)$-structure. In particular, we need to carefully track any
possible dependence of the various constants and on the parameters
$A_{n}$, $n=1,\dots,N$. We start with a single approximation, showing
in particular independence of the characteristics of $\appAq\varphi$ on $A\geq1$.

\begin{lemma}
  \label{lem:eq12}
  Let $\varphi$ be a balanced N-function with characteristics
  $(\gamma_1,\gamma_2)$. Then, for
  all $A\ge 1$ and $q\ge 2$ the $(A,q)$-approximation $\appAq\varphi$
  is a balanced N-function with characteristics
  $\big(\gamma_1, \max \set{\gamma_2,q-1}\big )$. 
\end{lemma}
\begin{proof} 
By construction we have 
\begin{gather*}
  \appAq\varphi \in    C^{1}(\setR^{\geq0})\cap C^{2}(\setR^{>0}),\qquad
  \appAq\varphi (0)=(\appAq\varphi )'(0)=0,
  \\
  (\appAq\varphi )''(t)>0\qquad \textrm{ for }t >0\,.
\end{gather*}
For $t\leq A$ we have $\appAq\varphi (t)=\varphi(t)$, while
$\frac{\appAq\varphi (t)}{t}=\alpha_{2,q}\,t^{q-1}+\alpha_{1,q}+\frac{\alpha_{0,q}}{t}$,
for $t>A$, which implies that $\appAq\varphi $ is a regular N-function, since
\begin{equation*}
\lim_{t\to0^{+}}  \frac{\appAq\varphi (t)}{t}=0\qquad \text{and}\qquad \lim_{t\to{+}\infty}
\frac{\appAq\varphi (t)}{t}=\infty\,,
\end{equation*}
where we used in the first limit that $\varphi$ is an N-function. 
From $\appAq\varphi (t)=\varphi(t)$, for $t\leq A$, we get, for all $t
\in (0,A]$, that 
\begin{equation}\label{eq:equiUA}
  \gamma_1   (\appAq\varphi  )'(t)\le (\appAq\varphi
  )''(t)\,t \le \gamma_2  \,  (\appAq\varphi  )'(t)\,,
\end{equation}
since $\varphi $ is balanced. On the other hand, for $t\geq A$ we have 
\begin{equation*}
  \frac{(\appAq\varphi )'(t)}{(\appAq\varphi )''(t)\,t}=  \frac{q\alpha_{2,q} \,t^{q-1}+\alpha_{1,q}}
{q(q-1)\alpha_{2,q}\,t^{q-1}}=:g^{A}(t).
\end{equation*}
Observe that for fixed $A\ge 1$ there holds
$\lim_{t\to\infty}g^{A}(t)=\frac 1{q-1}$, while 
\begin{equation*}
  g^{A}(A)=
  \frac{\varphi'(A)}{\varphi''(A)\,A}\in\left[\gamma_{2}^{-1},\gamma_{1}^{-1}\right]\,. 
\end{equation*}
From $(g^{A})'(t)=-\frac{1}{t^{q-2}}\frac{\alpha_{1,q}}{q\alpha_{2,q}}$
it follows that the sign of the derivative depends only on
$\alpha_{1,q}=\varphi'(A)-\varphi''(A)\,A$. Thus, $g_{A}(t)$ is
monotone. Distinguishing between $\alpha_{1,q}\ge 0$
and $\alpha_{1,q}\le 0$, using $\gamma_1\le 1$ and $\gamma_2\ge 1$, as
well as \eqref{eq:equiUA}, one easily sees that for all $t\ge 0 $ there
holds
\begin{align*}
\min \Bigset{\frac 1{\gamma_2},
  \frac 1{q-1}}\leq g^A(t)\leq  \frac 1{\gamma_1}\,, 
\end{align*}
implying the assertion. 
\end{proof}
\begin{corollary}
  \label{cor:eq12}
  For all $A\ge 1$ and $q\ge 2$ the $(A,q)$-approximation
  $\appAq{(\function_{p,\delta})}$ of $\function_{p,\delta}$ with
  $p \in (1,\infty)$ and $\delta\in [0,\infty)$ is a balanced
  N-function with characteristics
  \begin{equation*}
    \big(  \min \set{1,p-1}, \max  \set{1,p-1,q-1}  \big ).
  \end{equation*}
\end{corollary}
\begin{proof}
  This follows immediately from   Lemma \ref{lem:function} and
  Lem\-ma~\ref{lem:eq12}. 
\end{proof}

We have the following analogue of
Proposition~\ref{prop:potential-equivalence} for
$(A,q)$-approximations of $\potV$ and 
$\varphi$.
\begin{lemma}
\label{lem:pot-equiA}
Let $\varphi$ be a balanced N-function with characteristics
$(\gamma_1,\gamma_2)$. Let the operator $\bS=\partial\potV$
have $\varphi$-structure with
characteristics $(\gamma_3,\gamma_4,\Delta_2(\varphi))$. For $A\ge 1$
and $q\ge 2$ let $\appAq{\potV} $ and $\appAq{\varphi} $ be the
$(A,q)$-approximation of $\potV$ and $\varphi $, respectively. Then,
there holds for all $t>0$
\begin{equation*}
     \frac{ \gamma_3 }{\gamma_2}(\appAq\varphi)''(t)\le (\appAq\potV)''(t)\le
     \frac{ \gamma_4 }{\gamma_1}(\appAq\varphi)''(t)\,.
  \end{equation*}
\end{lemma}
\begin{proof}
  By definition we have $\appAq{\potV}(t)=\potV(t)$ and
  $\appAq{\varphi}(t) = \varphi (t)$ for $t\leq A$. Thus, the
  assertions for $t\le A$ follow from \eqref{eq:simU}.  For $t\ge A$
  we have $(\appAq{\potV})''(t)=\potV''(A) \frac{t^{q-2}}{A^{q-2}}$ and
  $(\appAq{\varphi})''(t)=\varphi''(A)
  \frac{t^{q-2}}{A^{q-2}}$. Thus, for $t\ge A$ the  assertion
  follows again from \eqref{eq:simU}. 
\end{proof}

The properties of the function $\appAq{\potV}$ proved in the previous
lemmas allow us to show that the operator $\appAq{\bS}$ has
$\appAq{\varphi}$-structure.

\begin{proposition}\label{prop:SA-struc}
  Let $\varphi$ be a
  balanced N-function with characteristics $(\gamma_1,\gamma_2)$. Let the operator
  $\bS=\partial\potV$ 
  have $\varphi$-structure with characteristics
  $(\gamma_3,\gamma_4,\Delta_2(\varphi))$. For
  $A\ge 1$ and $q\ge 2$ let $\appAq{\potV} $ and $\appAq{\varphi} $ be
  the $(A,q)$-approximation of $\potV$ and $\varphi $,
  respectively. Then, the operator
  $\appAq{\bS}:=\partial \appAq{\potV}$ has both
  $\appAq{\potV}$-structure and $\appAq{\varphi}$-structure, with
  characteristics depending only on
  $\gamma_1,\gamma_2, \gamma_3, \gamma_4$, and $q$.
\end{proposition}

\begin{proof}
  The operator $\appAq{\bS}$ is derived from the potential $\appAq{\potV}$
  which, according to Proposition \ref{prop:potential-equivalence} and
  Lemma~\ref{lem:eq12}, is a balanced N-function with 
  characteristics
  $\big(\frac{\gamma_3}{\gamma_4} \,\frac{ \gamma_1^2}{\gamma_2},
  \max\set{q-1,\frac{\gamma_4}{\gamma_3} \,\frac{
      \gamma_2^2}{\gamma_1}}\big )$. Thus, the proof of Proposition
  \ref{prop:pot-phi} shows that $\appAq{\bS}$ has
  $\appAq{\potV}$-structure with characteristics
  $$
  \big(\frac{\gamma_3}{\gamma_4} \,\frac{ \gamma_1^2}{\gamma_2},
  2+\max\bigset{q-1, \frac{\gamma_4}{\gamma_3} \,\frac{
      \gamma_2^2}{\gamma_1}}, \Delta_2(\appAq{\potV})\big )\,,
  $$
  where $\Delta_2(\appAq{\potV})$ depends
  only on $\max\set{q-1,\frac{\gamma_4}{\gamma_3} \,\frac{
      \gamma_2^2}{\gamma_1}}$, according to Lemma \ref{lem:d2}. Now Lemma \ref{lem:pot-equiA} yields
  that the operator $\appAq{\bS}$ has $\appAq{\varphi}$-structure
  with characteristics
  $$
  \big(\frac{\gamma_3^2}{\gamma_4} \,\frac{ \gamma_1^2}{\gamma_2^2}, \frac
  {\gamma_4}{\gamma_1}
  \big ( 2 + \max\bigset{q-1, \frac{\gamma_4}{\gamma_3} \,\frac{
      \gamma_2^2}{\gamma_1}}\big) , \Delta_2(\appAq{\varphi})\big).
  $$
  Lemma \ref{lem:eq12}  and Lemma
  \ref{lem:d2} yield that $\Delta_2(\appAq{\varphi})$ depends only
  on $\max \set{\gamma_2,q-1} $.  This finishes the proof.
\end{proof}
We have the following crucial result (cf.~Proposition~\ref{prop:hammer-phi}).
\begin{proposition}
\label{prop:SA-ham}
Let $\varphi$ be a balanced N-function with characteristics
$(\gamma_1,\gamma_2)$. Let the operator $\bS=\partial\potV$
have $\varphi$-structure with
characteristics $(\gamma_3,\gamma_4,\Delta_2(\varphi))$.  For ${A\ge 1}$
and $q\ge 2$ let $\appAq{\varphi} $ and $\appAq{\bS} $ be the
$(A,q)$-approximation of $\varphi$ and $\bS $, respectively, and let
$\bF_{\appAq{\varphi}}$ be defined in \eqref{def:F}. Then, we have for
all $\bP,\bQ \in \setR^{3\times 3} $ that
  \begin{align*}
    (\appAq{\bS}(\bP)-\appAq{\bS}(\bQ))\cdot(\bP-\bQ)
    &\sim \MC_{\appAq{\varphi}}(\abs{\bfP^{\sym}} +
      \abs{\bfP^{\sym}-\bfQ^{\sym}})\,|\bP^{\sym}-\bQ^{\sym}|^{2}  \,,
    \\
    &      \sim |{\bF}_{\appAq{\varphi}}(\bP)-{\bF}_{\appAq{\varphi}}(\bQ)|^{2}, 
    \\
    |   \appAq{\bS}(\bP)-\appAq{\bS}(\bQ)|&\sim {\MC}_{\appAq{\varphi}}(\abs{\bfP^{\sym}} +
        \abs{\bfP^{\sym}-\bfQ^{\sym}})\,|\bP^{\sym}-\bQ^{\sym}|\,, 
\end{align*}
where the constants of equivalence depend only on $\gamma_1, \gamma_2,
\gamma_{3},\gamma_{4}$, and $q$.
\end{proposition}
\begin{proof}
  This is a direct consequence of Proposition \ref{prop:SA-struc} 
  and Proposition \ref{prop:hammer-phi}.
\end{proof}
\begin{remark}
  For the limiting processes, it is of fundamental relevance that in Proposition~\ref{prop:SA-ham}
  the constants do not depend on $A\geq1$.
\end{remark}

Based on Proposition \ref{prop:SA-ham} we can prove the validity of
equivalent expressions for $\nabla\bF_\varphi(\bD\bu)$ which play a
crucial role in the proof of regularity for the problem
\eqref{eq:pfluid} (cf.~\cite{br-plasticity,br-parabolic,SS00}). To
this end, we define for a sufficiently smooth operator $\bS:\,\setR^{3\times3}\to\setR^{3\times
    3}_{\sym}$ the functions
$\mathbb P_i\colon \setR^{3\times 3} \to \setR$, $i=1,2,3$ via
\begin{equation*}
  \mathbb{P}_i (\bP):=\partial_i\bS(\bP)\cdot\partial_i\bP =
  \sum_{j,k,l,m=1}^3\partial_{jk}S_{lm}(\bP) \,\partial_i P_{jk}\,\partial_i P_{lm }\,,
\end{equation*}
and emphasize  that there is no summation over the index $i$.

If $\bS^n$, $n\in \set{1,\ldots , N}$, is a multiple approximation of
$\bS$ we define analogously $\mathbb P_i^n \colon \setR^{3\times 3}
\to \setR$, $i=1,2,3$, for $n\in \set{1,\ldots , N}$, via
\begin{equation*}
  \mathbb{P}^n_i (\bP):=\partial_i\bS^n(\bP)\cdot\partial_i\bP =
  \sum_{j,k,l,m=1}^3\partial_{jk}S^n_{lm }(\bP) \,\partial_i P_{jk}\,\partial_i P_{lm }\,.
\end{equation*}

\begin{proposition}\label{prop:pFA}
  Let the operator $\bS$ have $(p,\delta)$-structure for some
  $p \in (1, \infty)$ and $\para\in [0,\infty)$, with characteristics
  $(\gamma_3,\gamma_4,p)$ and let $\bS$ be derived from a potential $\pot$ with
  characteristics $(\gamma_1,\gamma_2)$. For given $N\in \setN$ let
  $\bS^n$, $n \in \set{1,\ldots , N}$ be a multiple approximation of
  $\bS$. If for a vector field $\bv \colon \Omega\subset \setR^3 \to \setR^3$ there
  holds $ \bF^n(\bD\bv) \in W^{1,2}(\Omega)$, then we have for
  $i=1,2,3$ and a.e.~in $\Omega$ the following equivalences
  \begin{align}\label{eq:pn}
    \begin{aligned}
      |\partial_i \bF^{n}(\bD\bv)|^{2} &\sim\MC^n(|\bD\bv|)\,
      |\partial _i\bD\bv|^{2}
      \\
      &\sim \mathbb P^{n}_i(\bD\bv)\,,
      \\
      |\partial_i \bS^{n}(\bD\bv)|^{2} &\sim\MC^n(|\bD\bv|)\, \mathbb
      P^n _i(\bD\bv)\,.
    \end{aligned}
  \end{align}
where the constants of equivalence depend only on
$\gamma_1,\gamma_2,\gamma_3,\gamma_4, p$, and $q_{n}$, for \linebreak ${n=1,\ldots,N}$. 
\end{proposition}
\begin{proof}
  For $h>0$ and $i=1,2,3$ let
  $\Delta^+_i \bv (\bx):= \bv (\bx +h\,\be_i)-\bv (\bx)$ and
  $d^+_i\bv (\bx):= h^{-1}\Delta^+_i\bv (\bx)$ be the classical increments
  and difference
  quotients in direction $\be_i$ of the canonical basis. The standard
  theory of difference quotients (cf.~\cite{breit-pde}) and
  $\bF^n(\bD\bv) \in W^{1,2}(\Omega)$ yield that
  $d^+_i \bF^n(\bD\bv) \to \partial_i \bF^{n}(\bD\bv)$ a.e.~and in
  $L^2_{\loc}(\Omega)$ as $h\to 0$ and
  \begin{align*}
    \int\limits_{\Omega_h} \abs{d^+_i\bF^n(\bD\bv) }^2\, dx \le c\,,
  \end{align*}
  where the constant $c$ is independent of $h$, and where we used the
  notation 
  \begin{equation*}
    \Omega_{h}:=\left\{\bx\in \Omega \fdg d(\bx,\partial\Omega)>2h\right\}.
  \end{equation*}
 Thus, Proposition~\ref{prop:SA-ham}, Lemma~\ref{lem:UA} for $p\le 2$, and
  Lemma~\ref{lem:UAm} for $p>2$ yield that
  \begin{align*}
    \int\limits_{\Omega_h}\abs{d^+_i \bD\bv }^2\, dx \le c\,,
  \end{align*}
  with $c$ independent of $h$ (even if it may depend on $\delta$ and
  $A_{n}$).  Consequently, we obtain that
  $\bD\bv \in W^{1,2}(\Omega)$, and
  $d^+_i \bD\bv \to \partial _i \bD\bv$,
  $\Delta^+_i \bD\bv \to \bfzero $ a.e.~and in $L^2_{\loc}(\Omega)$, as
  $h\to 0$. Proposition \ref{prop:SA-ham} yields
  \begin{align*}
    \abs{d^+_i\bF^n(\bD\bv)(\bx) }^2
    &\sim \MC^n(\abs{\bD\bv(\bx)} +\abs{\Delta^+_i\bD\bv(\bx)})\,
      \abs{d^+_i\bD\bv (\bx)}^2\,, 
  \end{align*}
which implies, using the above proved convergences, \eqref{eq:pn}$_1$ as $h \to 0$.

  Proposition \ref{prop:SA-struc} shows that $\bS^n$ has
  $\function^n$-structure, which implies
  \begin{equation*}
    \mathbb P^n_i (\bD\bv )=\sum_{j,k,l,m=1}^3\partial_{jk}S_{lm}^{n}(\bD\bv) \,\partial_i
    D_{jk}\bv\,\partial_i D_{lm}\bv \sim  \MC^n(|\bD\bv|)
    |\partial_{i}\bD\bv|^2 \,,
  \end{equation*}
  showing \eqref{eq:pn}$_2$. To prove \eqref{eq:pn}$_3$ we use the
  definition of $\mathbb P_i^n$ and \eqref{eq:pn}$_{1,2}$ to get
  \begin{equation*}
    \MC^n(|\bD\bv|)\, \big (\mathbb{P}_i^{n}(\bD\bv)\big)^2 \le
    \MC^n(|\bD\bv|)\,
    |\partial_{i}\bS^{n}(\bD\bv)|^{2}|\partial_{i}\bD\bv|^{2} \sim 
    {\mathbb P^n_{i}(\bD\bv)}\,|\partial_{i}\bS^{n}(\bD\bv)|^{2}\,, 
  \end{equation*}
which implies 
  \begin{equation*}
     \MC^n(|\bD\bv|)\, \mathbb{P}_i^{n}(\bD\bv)    \le c\, |\partial_{i}\bS^{n}(\bD\bv)|^2\,.
   \end{equation*}
   On the other hand, the fact that $\bS^n$ has
   $\function^n$-structure and \eqref{eq:pn}$_{1,2}$ imply that
   $\mathbb{P}_{i}(\bP)=0$  if and only if
   $\bP=\bfzero$. Consequently, we obtain  
  \begin{equation*}
    \begin{aligned}
      |\partial_{i}\bS^{n}(\bD\bv)|^{2}&\leq
      \sum_{k,l=1}^3|\partial_{kl}
      \bS^{n}(\bD\bv)\,\partial_{i}D_{kl}\bv|^{2} \leq c\, \left
        (\MC^n(|\bD\bv|)\right)^{2}|\partial_{i}\bD\bv|^{2}
      \\
      & \leq c\, \MC^n(|\bD\bv|)\,|\partial_{i}\bF^n(\bD\bv)|^{2} \leq
      c\, \MC^n(|\bD\bv|)\,\mathbb P^n_{i}(\bD\bv)\,.
    \end{aligned}
  \end{equation*}
  Note that all constants just depend on the quantities indicated in the
  formulation of the assertion. This yields  the reverse estimate, proving
  \eqref{eq:pn}$_3$. 
\end{proof}

To derive a priori estimates and to perform the limiting process we
need estimates, 
which do not involve the parameters $A_{n}$, for $n=1,\dots,N$. We restrict ourselves to
the case that $\varphi =\function_{p,\para}=\function$ with
$p \in (1,\infty)$ $\delta\in[0,\infty)$ and distinguish between the cases $p\le 2$ and $p> 2$ for the sake of a simpler presentation.

\subsection{Some estimates specific to the case $p> 2$}
In this section we prove some results, which are specific of the
case $p>2$. In particular, in the case $p\geq\frac{13}{3}$ we need
multiple approximations, which makes the results more technical.

\begin{lemma}\label{lem:UAm}
  For given $p>2$, $\para >0$, and $ N \in \setN$ let $\MC^n$, $n\in \set{1,\ldots, N}$,
  be a multiple
    approximation of $\MC^0$.  Then, there exist $\widehat A_n
    =\widehat{A}_n(\para,p, q_1,\ldots, q_n )\ge 1$ such that for all
    $A_n\ge 
    \max\set{\para, \widehat A_n}$ the
  function ${\MC^n}$ is non-decreasing and satisfies for all $t \ge 0
  $  
  \begin{gather}
    \label{eq:UA}
    \begin{gathered}
      \frac 1{(p-1) 2^{q_1-2}} \,\para^{p-2}\le \frac 1{p-1} \,
      \frac{\para^{p-q_n}}{2^{q_1-2}}\, \MC_{\function_{q_n,\para}}
      (t) \le {\MC^n} (t)\,,
      \\
      {\MC^n} (t) \le \frac {p-1}{q_n-1}\, 2^{p-2}\,
      (A_{n-1})^{p-q_{n-1}}\,\MC_{\function_{q_{n-1},\para}}(t)\,,
      \\
      {\MC^n} (t)  \le     \frac {p-1}{q_n-1}\, 2^{p-2}\, \MC^0(t)\,.
      \end{gathered}
  \end{gather}
\end{lemma}
\begin{proof}
  For ease of presentation we show the assertion just in the first two
  cases, i.e., $n=1,2$. The remaining cases follow in the same way and
  are left to the interested reader. 
  
  The case $\mathbf{ (n=1):}$ For simplicity we set $A:=A_1$ and $q:=q_1$. For
  $t \le A$ we have $\MC^1 (t) =\MC^0(t)=(\para +
  t)^{p-2}$. Thus, $q_0=p>q\ge 2$ implies
  \begin{align*}
    \para^{p-2}\le \para ^{p-q}\, (\para +t)^{q-2}\le (\para
    +t)^{p-2}\,,
  \end{align*}
  and $\frac 1{(p-1) 2^{q-2}}  \le 1 \le \frac {p-1}{q-1}
  (1+\para)^{p-2}$, which proves \eqref{eq:UA} for $n=1$ and $t\le A$. Moreover,
  $(\para +t )^{p-2} $ is an increasing function in $t$.

  For $t \ge A$ we have 
  \begin{align}\label{def:gA}
    \MC^1(t)
    &= q\, \alpha_{2,q}\, t^{q-2}\big (1 + \frac
      {\alpha_{1,q}}{ q\, \alpha_{2,q}\, t^{q-1}}\big )
      =:  q\, \alpha_{2,q}\, t^{q-2}g^A(t)\,,
  \end{align}
  where $\alpha_{i,q}=\alpha_{i,q}(\function^0)$, $i=1,2$. The
 expressions for the coefficients $\alpha_{i,q}$, $i=1,2$, imply
  $g^A(A)=(q-1)\frac{(\function^0)'(A)}{(\function^0)''(A)A}$,
  $\lim_{t\to \infty}g^A(t)=1$, and
  \begin{align}\label{eq:gA'}
  (g^A)'(t)=\frac{\alpha_{1,q}(1-q)}{q\, \alpha_{2,q}\,
    t^{q}}&= (q-1)^2\, \frac{A^{q-1}}{t^q} \Big (\frac 1{q-1} -
  \frac{(\function^0)'(A)}{(\function^0)''(A)A } \Big )\,.
  \end{align}
  From the properties of $\function ^0$ it follows that
  $\lim_{t\to \infty} \frac{(\function^0)'(t)}{(\function^0)''(t)t}=
  \frac 1{q-1}$ and that
  $\frac{(\function^0)'(t)}{(\function^0)''(t)t}$ is strictly monotone
  increasing. Thus, for $A \ge \widehat A(p,q,\para)$ the expression
  in the parenthesis in \eqref{eq:gA'} is non-negative and thus $g^A$
  is a non-decreasing function. Consequently, we get that $\MC^1$ is a
  non-decreasing function, since $q\, \alpha_{2,q}\, t^{q-2}$ is
  increasing for $q>2$ (non-decreasing for $q=2$). Using these
  properties and that $\function^0 $ is balanced we obtain that for
  $t \ge A \ge \widehat A$ there holds
  \begin{align}\label{eq:gA}
    \frac {q-1}{p-1} \le g^A(t) \le 1\,.
  \end{align}
  It remains to estimate the factor in front of $g^A$ in
  \eqref{def:gA}. From the expression for  $\alpha_{2,q}$ we get that
  $ q\, \alpha_{2,q}\, t^{q-2}=\frac{(\function^0)''(A)}{q-1} \big (\frac
  tA \big )^{q-2}$. Thus, using \eqref{eq:E}, $t\ge A\ge \para\ge 0$,
  and $2\le q<p$ 
  we obtain
  \begin{align}      \label{eq:1u}
    q\, \alpha_{2,q}\, t^{q-2}
    &\le \frac{p-1}{q-1} \frac {(\para
      +A)^{p-2}}{A^{q-2}} \, t^{q-2} = \frac{p-1}{q-1} \Big
      (\frac{\para}A +1\Big )^{p-2}\, A^{p-q}\,
      t^{q-2} \notag 
    \\
    &\le \frac{p-1}{q-1} \, 2^{p-2}\, A^{p-q}\, (\para + t)^{q-2} \notag
      =
      \frac{p-1}{q-1} \, 2^{p-2}\,\Big (\frac {\para +t}A\Big )^{q-p}
      \, (\para +t)^{p-2}\notag
    \\
    &     \le \frac{p-1}{q-1} \,2^{p-2}\, (\para +t)^{p-2}\,.
  \end{align}
  For $A\ge \delta$ we get that $t\ge 2^{-1}\,
  (\delta+t)$ for all $t \ge A$. Using this, the definition of
  $\function^0$, \eqref{eq:E}, $2\le q<p$ and $t\ge A\ge \para\ge 0$ we obtain
  \begin{align}
    q\, \alpha_{2,q}\, t^{q-2}
    &\ge \frac{1}{q-1} \frac {(\para
      +A)^{p-2}}{A^{q-2}} \, t^{q-2}
      \ge  \frac{1}{q-1} (\para +A )^{p-q} \frac{1}{2^{q-2}}\, (\para +t)^{q-2}\notag
    \\
    &\ge \frac{\para ^{p-q}}{q-1} \frac{1}{2^{q-2}} \,
      \MC_{q,\para}(t)\ge  \frac{\para ^{p-2}}{p-1} \frac{1}{2^{q-2}}
      \,.
      \label{eq:1l}
  \end{align}
  The inequalities \eqref{eq:gA}--\eqref{eq:1l} imply \eqref{eq:UA}  
  for $n=1$ and $t\ge A \ge \max\{\para,\widehat A\}$. This completes the proof for
  $n=1$.\\[-3mm]

  $\mathbf{ (n=2):}$ For simplicity we set $B:=A_2\geq A_{1}=:A$ and
  $r:=q_2$, $q:=q_1$. For $t \le B$ we have
 ${\MC}^2 (t) =\MC^1(t)$. Thus, $p>q>r\ge 2$ implies
  \begin{align*}
    \para^{p-2}\le \para ^{p-r}\, (\para +t)^{r-2}\le \para^{p-q}\,(\para
    +t)^{q-2}\,,
  \end{align*}
  which together with \eqref{eq:UA} for $n=1$ shows \eqref{eq:UA}$_1$
  for $n=2$ and $t\le B$. The estimate \eqref{eq:UA}$_{2,3}$ for $n=2$
  and $A\le t\le B$ follows from ${\MC}^2 (t) =\MC^1(t)$,
  \eqref{def:gA}, \eqref{eq:gA}, the estimates in \eqref{eq:1u} and
  $r<q$; while for $t \le A$ it follows from $\MC^2(t)=\MC^0(t)$,
  $\para\le A$, $2\le r<q<p$, and
  \begin{align*}
    (\para +t)^{p-2}=    (\para +t)^{p-q}\,    (\para +t)^{q-2}\le
    2^{p-q}   \, (\para +t)^{q-2}\,. 
  \end{align*}

  For $t \ge B$ we have 
  \begin{align}\label{def:gA2}
    {\MC}^2(t)
    &= r\, \alpha_{2,r}\, t^{r-2}\big (1 + \frac
      {\alpha_{1,r}}{ q\, \alpha_{2,r}\, t^{r-1}}\big )
      =:  r\, \alpha_{2,r}\, t^{r-2}h^B(t)\,,
  \end{align}
    where $\alpha_{i,r}=\alpha_{i,r}(\function_1)$, $i=1,2$.
  The expressions of the coefficients $\alpha_{i,r}$, $1=1,2$,
  imply $h^B(B)=(r-1)\frac{(\function^{1})'(B)}{(\function^{1})''(B)B}$, 
  $\lim_{t\to \infty}h^B(t)=1$, and
  \begin{align}\label{eq:hB'}
  (h^B)'(t)=\frac{\alpha_{1,r}(1-r)}{r\, \alpha_{2,r}\,
    t^{r}}&= (r-1)^2\, \frac{B^{r-1}}{t^r} \Big (\frac 1{r-1} -
  \frac{(\function^{1})'(B)}{(\function^{1})''(B)B } \Big )\,.
  \end{align}
  In the proof of Lemma \ref{lem:eq12} we showed that
  $\lim_{t\to \infty} \frac{(\function ^1)'(t)}{(\function^1)''(t)\,t}=
  \frac 1{q-1}$ and that $\frac{(\function ^1)'(t)}{(\function^1)''(t)\,t}$
  is strictly monotone. Thus, for $B \ge \widehat B(p,q,r,\para)$ the expression in
  the parenthesis in \eqref{eq:hB'} is non-negative and thus $h^B$ is a non-decreasing
  function. Consequently, we
  get that $\MC^2$ is a non-decreasing function, since
  $r\, \alpha_{2,r}\, t^{r-2}$ is increasing for $r>2$. Using these properties
  and Lemma \ref{lem:eq12} for $\function^1 $ we obtain that for
  $t \ge B \ge \widehat B$ there holds 
  \begin{align}\label{eq:gA2}
    \frac {r-1}{p-1} \le h^B(t) \le 1\,.
  \end{align}
  It remains to estimate the factor in front of $h^B$ in
  \eqref{def:gA2}. From the expressions for  $\alpha_{2,r}$ and $(\function^1)''(B)$ we get that
  $ r\, \alpha_{2,r}\, t^{r-2}=\frac{(\function^1)''(B)}{r-1} \big
  (\frac tB \big )^{r-2}=\frac{(\function^0)''(A)}{r-1} \big (\frac BA
  \big )^{q-2}\big (\frac tB \big )^{r-2}$. Thus, using \eqref{eq:E},
  the definition of $\function ^1$, $2\le r<q<p$ and
  $0\le \para\le A\le B\le t$ we obtain
  \begin{align*}
    r\, \alpha_{2,r}\, t^{r-2}
    &\le \frac{p-1}{r-1}  \Big (\frac{\para}A
      +1\Big )^{p-2}\, \Big (\frac tB\Big )^{r-q}\, A^{p-q}\, t^{q-2}\notag
    \\
   &\le \frac{p-1}{r-1} \, 2^{p-2}\, A^{p-q}\,(\para + t)^{q-2} =
     \frac{p-1}{r-1} \, 2^{p-2}\,\Big (\frac {\para +t}A\Big )^{q-p}
     \, (\para +t)^{p-2}\notag
    \\
    &     \le \frac{p-1}{r-1} \,2^{p-2}\, (\para +t)^{p-2}\,.
  \end{align*}
  For $B\ge A\ge \delta$ we get that $t\ge 2^{-1}\,
  (\delta+t)$ for all $t \ge B$. Using this, the definition of
  $\function^0$, \eqref{eq:E}, $2\le r<q<p$ and $t\ge B\ge
  A\ge \para\ge 0$ we obtain
  \begin{align*}
     r\, \alpha_{2,r}\, t^{r-2}
    &\ge \frac{1}{r-1} \frac {(\para
                     +A)^{p-2}}{A^{q-2}} \,B^{q-r}\, t^{r-2}
   \ge  \frac{1}{r-1} (\para +A )^{p-q} \,B^{q-r}\,\frac{1}{2^{r-2}}\, (\para +t)^{r-2}\notag
    \\
    &\ge \frac{\para ^{p-r}}{r-1} \frac{1}{2^{r-2}} \,
      \MC_{r,\para}(t) \ge \frac{\para ^{p-2}}{r-1} \frac{1}{2^{r-2}}
      \,.
  \end{align*}
  The last two inequalities and \eqref{eq:gA2} imply \eqref{eq:UA}  
  for $n=2$ and $t\ge B \ge \max\{A,\widehat B\}$. This completes the proof for
  $n=2$.
\end{proof}
In the proof of regularity  we will need
mainly the following corollary.
 \begin{corollary}\label{cor:UAm}
   Let the operator $\bS=\partial \pot$, derived from the potential $\pot$, have
   $(p,\para)$-structure with $p >2$ and $\delta> 0$, and
   characteristics $(\gamma_3,\gamma_4,p)$. For
   $ N \in \setN$ let $\function^n$, $\bF^n$, $\bS^n$, $n\in \set{1,\ldots, N}$, be a multiple
   approximation of $\function^0$, $\bF^0$, $\bS^0$.  Then, for
   all $A_n\ge \max\set{\para, \widehat A_n}$ with $\widehat A_n$ from Lemma
   \ref{lem:UAm} and all $t \ge 0 $ there holds
  \begin{gather}
      \frac {\para^{p-2}}{(p-1) 2^{q_1-1}}\,t^2\le \,
      \frac{\para^{p-q_n}}{ (p-1)2^{q_1-2}}\, {\function_{q_n,\para}}
      (t) \le {\function^n} (t)  \le     \frac {p-1}{q_n-1}\,
      2^{p-2}\, \function^0(t)\,, \notag
      \\
      {\function^n} (t) \le \frac {p-1}{q_n-1}\, 2^{p-2}\,
      (A_{n-1})^{p-q_{n-1}}\,{\function_{q_{n-1},\para}}(t)\,,    \label{eq:UA3m}
      \\
    (\function^n)^*(t)\le \frac {(p-1) 2^{q_1-3}} {\para^{p-2}} \,
    t^2\,.\notag 
  \end{gather}
  Moreover, for all $\bP \in \setR^{3\times 3}$  there holds 
  \begin{align}
    \label{eq:UA4m}
    \begin{aligned}
      \abs{\bF^n(\bP)}^2&\sim \function^n(\abs{\bP^{\sym}})\,,
      \\
      c\, \para ^{p-q_n}\,\abs{\bF_{\function _{q_n,\para}}(\bP)}^2
      &\le \abs{\bF^n(\bP)}^2
      \,, 
      \\
    \abs{\bS^n(\bP)} &\le c\, (A_{n-1})^{p-q_{n-1}}\,
(    \function_{q_{n-1},\para})' (\abs{\bP^{\sym}}) \,
    \end{aligned}
  \end{align}
  with constants $c$ depending only on $\gamma_3, \gamma_4$, $q_n$, $q_1$ and $p$. 
\end{corollary}
 \begin{proof}
   The proof of the estimates \eqref{eq:UA3m}$_{1,2}$ is a direct
   application of the previous lemma, the definition in
   \eqref{eq:mot}, $\function ^n(0) =\function _{q_n,\delta}(0) =\function _{q_{n-1},\delta}(0)=0$ and
   integration.

   Estimate \eqref{eq:UA3m}$_{3}$ follows from \eqref{eq:UA3m}$_{1,2}$
   and the equivalent definition of the complementary function since
  \begin{align*}
   (\function^{n})^*(t)&:= \sup _{s \ge 0} s\,t -\function^{n}(s)
   \\
                       &\le \sup _{s \ge 0} s\,t -      \frac {\para^{p-2}}{(p-1) 2^{q_1-1}}\,s^2
   \\
                     &= \frac {(p-1) 2^{q_1-3}} {\para^{p-2}}\,  t^2\,.
                        \end{align*}
The estimates \eqref{eq:UA4m} follow immediately from the definition
of multiple approximations of $\bS,\,\bF,\,\function$, $a$;
Proposition \ref{prop:SA-ham} and \eqref{eq:UA3m}.
\end{proof}

Let us finish this section with a more technical estimate needed in
the proof of regularity in the case $p>2$.
\begin{lemma}
  \label{lem:ast}
  For given $p>2$, $\para \ge 0$ and $ N \in \setN$ let $\MC^n$,
  $\function^n$ and $n\in \set{1,\ldots, N}$, be a multiple
  approximation of $\MC^0$ and $\function^0$ with $A_n\ge \max\set{\para,1}$, respectively.  Then,
  there exists a constant $c =c( p, q_1,\ldots, q_n )$
  such that for all $s,t \ge 0$ there holds 
  \begin{gather}
    \label{eq:ast}
    {\MC^n} (t)\, s^2 \le c\,
    \big(\delta^{p}+\function^n(s)+\function^n(t) \big )\,.    
  \end{gather}
\end{lemma}
\begin{proof}
  The assertion follows essentially from Young inequality and the
  expressions for the coefficients of the approximations. However, we
  have to distinguish several cases.  For ease of presentation we show
  the assertion just in the first two cases, i.e., $n=1,2$. The
  remaining cases follow in the same way and are left to the
  interested reader.

  The case $\mathbf{ (n=1):}$ For simplicity we set $A:=A_1$ and
  $q:=q_1$. For $t \le A$ we have
  $\MC^1 (t) =\MC^0(t)=(\para + t)^{p-2}$. For $s,t \le A$ Young
  inequality with $\frac p2$ and $\frac{p}{p-2}$, $\para \ge 0$ and
  $p>2$ yield
  \begin{align*}
    (\para +t)^{p-2}\, s^2
    &\le c \,\big ( (\para +t)^{p} + s^p\big ) \le  c\, \big (\para ^p
      +(\para +t)^{p-2}t^2+(\para +s)^{p-2} s^2 \big )\,, 
  \end{align*}
  which implies \eqref{eq:ast} for $s,t \le A$, since Corollary \ref{cor:sim} implies
  \begin{align}\label{eq:tA}
  (\para +t)^{p-2}t^2 \sim \function^0(t)= \function^1(t)\,,
  \end{align}
  valid for $t \le A$.
  
  Next, assume that $s,t \ge A$. Since $\function^1$ is
  balanced with characteristics $(1,p-1)$, the definition of $\function^1$ implies for $t\ge A$ that there holds
  \begin{align}\label{eq:stA}
    \MC^1(t)= \frac {(\function^1)'(t)}{t} \sim 
    (\function^1)''(t)= \frac {(\function^0)''(A)}{A^{q-2}}\, t^{q-2}\,,
  \end{align}
  with constants of equivalence depending only on $p$. This and
  Young inequality with $\frac q2$ and $\frac{q}{q-2}$ imply
  \begin{align*}
    \MC^1(t)\, s^2
    &\le c(p,q) \, \frac {(\function^0)''(A)}{A^{q-2}}\, \big (
      t^{q-2}\, t^2 + s^{q-2}\, s^2\big ) = c\, \big (
      (\function^1)''(t)\, t^2 +(\function^1)''(s)\, s^2\big )\,, 
  \end{align*}
 which yields \eqref{eq:ast} for $s,t \ge A$, since Corollary
  \ref{cor:sim} shows $(\function^1)''(t) \, t^2 \sim \function^1(t) $.

  Next, assume $s\le A\le t$. Using \eqref{eq:stA}, Young inequality
  with $\frac{p}{p-2}$ and $\frac p2$, $\para \ge 0$, \eqref{eq:tA},
  \eqref{eq:stA}, and again $(\function^1)''(t) \, t^2 \sim \function^1(t) $ we
  obtain
  \begin{align}\label{eq:sAt}
    \begin{aligned}
      \MC^1(t)\, s^2 &\le c(p ) \, \bigg (\Big ( \frac
      {(\function^0)''(A)}{A^{q-2}}\Big )^{\frac p{p-2}}\, t^{\frac
        {q-2}{p-2}\,p} + s^{p-2}\, s^2\bigg )
      \\
      &\le c \, \bigg (\frac {(\function^0)''(A)}{A^{q-2}}\,
      t^{q-2}\,t^2\, \Big ( \frac {(\function^0)''(A)}{A^{q-2}}\Big
      )^{\frac p{p-2}-1}\, t^{\frac {q-2}{p-2}\,p-q} + \function^{1}(s)\bigg )
      \\
      &\le c \, \bigg (\function^1(t)\, \Big ( \frac
      {(\function^0)''(A)}{A^{q-2}}\Big )^{\frac 2{p-2}}\, t^{2\frac
        {q-p}{p-2}}+ \function^1(s)\bigg )\,.
    \end{aligned}
  \end{align}
  Using $(\function^0)''(A)\sim (\para +A)^{p-2}$,
  $\max\{1,\delta\}\le A$ and $q<p$ we get
  \begin{align*}
    \Big ( \frac {(\function^0)''(A)}{A^{q-2}}\Big )^{\frac 2{p-2}}\,
    t^{2\frac {q-p}{p-2}}
    & \le c\, \frac {(\para +A)^2} {A^2} \,  \frac {t^{2\frac
      {q-p}{p-2}}}{A^{2\frac {q-2}{p-2}-2}} \le c\,
      \Big (\frac tA\Big )^{2\frac {q-p}{p-2}}\le c\,,
  \end{align*}
   which together with the last estimate implies \eqref{eq:ast} for $s
   \le A\le t$.

   Finally, for $t\le A\le s$ we get
   \begin{align*}
    \MC^1(t)\, s^2
    &= (\para +t)^{p-2} \, s^2= (\para +t)^{p-2} \, \Big ( \frac {(\function^0)''(A)}{A^{q-2}}\Big
      )^{-\frac 2q}\, \Big ( \frac {(\function^0)''(A)}{A^{q-2}}\Big
      )^{\frac 2q}\, s^2\,.
   \end{align*}
   We use Young inequality with $\frac q{q-2}$ and $\frac q2$,
   \eqref{eq:stA} and again $(\function^1)''(t) \, t^2 \sim
   \function^1(t) $  to
   arrive at
   \begin{align*}
    \MC^1(t)\, s^2
    &\le c(p,q) \, \bigg ( (\para + t)^{q\frac {p-2}{q-2}}\,\Big ( \frac {(\function^0)''(A)}{A^{q-2}}\Big
      )^{\frac {-2}{q-2}} + \frac {(\function^0)''(A)}{A^{q-2}}\,
      s^{q-2}\,s^2 \bigg      )\,.
   \end{align*}
   From $(\function^0)''(A)\sim (\para +A)^{p-2}$, $p>q\ge 2$, $t\le
   A$ and  $\para \ge 0$ we obtain 
   \begin{align*}
     (\para + t)^{q\frac {p-2}{q-2}}\,\Big ( \frac {(\function^0)''(A)}{A^{q-2}}\Big
      )^{\frac {-2}{q-2}}  
     &\le c\, (\para + t)^p \,  (\para +
       A)^{q\frac {p-2}{q-2}-p+2 - 2\frac{p-2}{q-2}}
     \\
     &\le c\, (\para + t)^p \le c\, (\para ^p +
     (\para +  t)^{p-2}\,t^2)\,,
   \end{align*}
   which together with \eqref{eq:tA}, \eqref{eq:stA} and the last estimate yields \eqref{eq:ast} for $t
   \le A\le s$. This finishes the proof for $n=1$. 

   \medskip

  The case $\mathbf{ (n=2):}$ For simplicity we additionally set $B:=A_2$ and $r:=q_2$. For
   $s, t \le B$ we have ${\MC}^{2} (t) =\MC^{1}(t)$. Thus, the assertion
   \eqref{eq:ast} for $n=2$ is already proved in the case $n=1$
   above. In the case $s,t \ge B$ we proceed exactly as in the case
   $s,t\ge A$ for $n=1$ just replacing $\function^1$, $\function^0$,
   $q$ and $A$ by $\function^2$, $\function^1$, $r$ and $B$,
   respectively.

   In the case $s\le B\le t$ we have to distinguish between $s\le A$
   and $A\le s\le B$. In the former case we use the analogue of
   \eqref{eq:stA} for $\MC^2$ and proceed as in \eqref{eq:sAt} to
   arrive at
   \begin{align*}
    \begin{aligned}
      \MC^2(t)\, s^2 &\le c(p) \, \bigg (\function^2(t)\, \Big ( \frac
      {(\function^1)''(B)}{B^{r-2}}\Big )^{\frac 2{p-2}}\, t^{2\frac
        {r-p}{p-2}}+ \function^2(s)\bigg )\,.
    \end{aligned}
  \end{align*}
  Using $(\function^1)''(B) = \frac{(\function^0)''(A)}{A^{q-2}} B^{q-2}
  \le (p-1) \frac{(\para +A)^{p-2}}{A^{q-2}} B^{q-2}$, $1\le A\le B$
  and $r<q<p$ we get the estimate
  \begin{align*}
    \Big ( \frac {(\function^1)''(B)}{B^{r-2}}\Big )^{\frac 2{p-2}}\,
    t^{2\frac {r-p}{p-2}}
    & \le c\, \frac {(\para +A)^2} {A^2} \, \Big (\frac AB\Big
      )^{2\frac {p-q}{p-2}}\Big (\frac tB\Big )^{2\frac {r-p}{p-2}}\le
      c\,, 
  \end{align*}
   which together with the last estimate implies \eqref{eq:ast} for $s
   \le A \le B\le t$. For $A\le s\le B\le t $ we have 
   \begin{align*}
    \MC^2(t)\, s^2
    &\le \frac{(\function^1)''(B)}{B^{r-2}}\,t^{r-2} \, \Big ( \frac {(\function^0)''(A)}{A^{q-2}}\Big
      )^{-\frac 2q}\, \Big ( \frac {(\function^0)''(A)}{A^{q-2}}\Big
      )^{\frac 2q}\, s^2\,.
   \end{align*}
   Using Young inequality with $\frac {q}{q-2}$ and $\frac q2$, the
   analogue of \eqref{eq:stA} for $\MC^2$ and $\MC^1$,
  and  $(\function^2)''(t) \, t^2 \sim \function^2(t) $ we get 
   \begin{align*}
     \MC^2(t)\, s^2
     &\le c(p,q) \, \bigg ( t^{q\frac {r-2}{q-2}}\,\Big ( \frac {(\function^0)''(A)}{A^{q-2}}\Big
       )^{\frac {-2}{q-2}}\,\Big ( \frac {(\function^1)''(B)}{B^{r-2}}\Big
       )^{\frac {q}{q-2}} + \frac {(\function^0)''(A)}{A^{q-2}}\,
       s^{q-2}\,s^2 \bigg      )
     \\
     & \le c\, \bigg ( \function^2(t)\,  t^{2\frac {r-q}{q-2}}\, \frac
       {A^{2}}{((\function^0)''(A))^{\frac {2}{q-2}}} \,\Big ( \frac {(\function^1)''(B)}{B^{r-2}}\Big
       )^{\frac {2}{q-2}} + \function_2(s)\bigg      )\,.
   \end{align*}
   From $(\function^1)''(B) = \frac{(\function^0)''(A)}{A^{q-2}}
   B^{q-2}$, $q>r\ge 2$ and  $t\ge
   B$ we obtain 
   \begin{align*}
     t^{2\frac {r-q}{q-2}}\, \frac   {A^{2}}{((\function^0)''(A))^{\frac {2}{q-2}}}
     \, \Big ( \frac {(\function^1)''(B)}{B^{r-2}}\Big  )^{\frac {2}{q-2}}
     = \Big ( \frac {t}{B}\Big  )^{2\frac {r-q}{q-2}} \le 1\,,
   \end{align*}
   which together with the last estimate yields \eqref{eq:ast} for $A
   \le s\le B\le t$.  

   In the case $t\le B\le s$ we have to distinguish between $t\le A$
   and $A\le t\le B$. In the former case we have
   \begin{align*}
    \MC^2(t)\, s^2
    &= (\para +t)^{p-2} \, \Big ( \frac {B^{r-2}}{(\function^1)''(B)}\Big
      )^{\frac 2r}\, \Big ( \frac {(\function^1)''(B)}{B^{r-2}}\Big
      )^{\frac 2r}\, s^2\,,
   \end{align*}
   which by Young inequality with $\frac {r}{r-2}$ and $\frac r2$,
   and the analogue of \eqref{eq:stA} for $\MC_2$ yields
   \begin{align*}
     \MC^2(t)\, s^2
     &\le c(p,r) \,\bigg ( (\para +  t)^{r\frac {p-2}{r-2}}\,\Big ( \frac {B^{r-2}}{(\function^1)''(B)}\Big
       )^{\frac {2}{r-2}} + \function^2(s)\bigg      )
     \\
     &\le c \,\bigg ( \big (\para ^p + \MC^0(t)\big )\, (\para +
       t)^{2\frac {p-r}{r-2}}\,
       \frac {B^{2}}{((\function^1)''(B))^{\frac {2}{r-2}}} + \function^2(s)\bigg      )\,.
   \end{align*}
   From
   $(\function^1)''(B) = \frac{(\function^0)''(A)}{A^{q-2}} B^{q-2}$,
   $q>r\ge 2$, $\para \ge 0$ and $B\ge A\ge t$ we obtain
   \begin{align*}
     (\para +
     t)^{2\frac {p-r}{r-2}}\, \frac
     {B^{2}}{((\function^1)''(B))^{\frac {2}{r-2}}}
     \le  \, \Big ( \frac {\para +t}{ \para +A}\Big  )^{2\frac {p-r}{r-2}}
     \,\Big ( \frac {A}{B}\Big  )^{2\frac {q-r}{r-2}} \le 1\,,
   \end{align*}
   which together with $\MC^0(t) =\MC^2(t)$, $t\le A$, and the last
   estimate yields \eqref{eq:ast} for $ t\le A\le B\le s$. For $A\le t\le B\le s $ we have 
   \begin{align*}
    \MC^2(t)\, s^2
    &\le \frac{(\function^0)''(A)}{A^{q-2}}\,t^{q-2} \, \Big ( \frac {B^{r-2}}{(\function^1)''(B)}\Big
      )^{\frac 2r}\, \Big ( \frac {(\function^1)''(B)}{B^{r-2}}\Big
      )^{\frac 2r}\, s^2\,,
   \end{align*}
   which, by Young inequality with $\frac {r}{r-2}$ and $\frac r2$
   and \eqref{eq:stA} for $\MC^2$ and $\MC^1$, yields
   \begin{align*}
     \MC^2(t)\, s^2
     &\le c(p,r) \,\bigg (  t^{r\frac {q-2}{r-2}}\,
       \frac{((\function^0)''(A))^{\frac r{r-2}}}{A^{r\frac {q-2}{r-2}}}\,
       \frac {B^{2}}{((\function^1)''(B))^{\frac {2}{r-2}}} +
       \function^2(s)\bigg      )
     \\
     &\le c \,\bigg (  \function^2(t) \, t^{2\frac {q-r}{r-2}}\,
       \frac{((\function^0)''(A))^{\frac 2{r-2}}}{A^{2\frac {q-2}{r-2}}}\,
       \frac {B^{2}}{((\function^1)''(B))^{\frac {2}{r-2}}}+
       \function^2(s)\bigg      )       \,.
   \end{align*}
   Using $(\function^0)''(A)\sim (\para + A)^{p-2}$ and
   $(\function^1)''(B) = \frac{(\function^0)''(A)}{A^{q-2}} B^{q-2}$,
   $q>r\ge 2$ and $B\ge t$ we obtain
   \begin{align*}
     t^{2\frac {q-r}{r-2}}\,
       \frac{((\function^0)''(A))^{\frac 2{r-2}}}{A^{2\frac {q-2}{r-2}}}\,
     \frac {B^{2}}{((\function^1)''(B))^{\frac {2}{r-2}}}
     \le c\, \Big ( \frac {t}{B}\Big  )^{2\frac {q-r}{r-2}} \le 1\,,
   \end{align*}
   which together with the last
   estimate yields \eqref{eq:ast} for $ A\le t\le B\le s$. This
   finishes the proof of the case $ n=2$.
\end{proof}
\subsection{Some estimates specific of the case $1<p\leq2$}
For completeness we deduce estimates for the case $p\in(1,2]$, which are the
counterpart of those proved in the previous section and which can be
used to prove the regularity results also in the case $p\in(1,2]$. Note
that in this case it is enough to use a single approximation with
$q=2$. 
\begin{lemma}\label{lem:UA}
  For $p\in (1,2]$ and $\delta> 0$ let
  $\function =\function_{p,\para}$ and $\MC =\MC_\function$. For
  $A\ge 1$ and $q=2$ we set $\function^{A}:=\function^{A,2}$ and
  $\MC^A:= \MC_{\function^{A,2}}$. Then, the
  function $\MC^{A}$ is non-increasing and satisfies for all $t\ge 0$
  \begin{align}
    \label{eq:wAl1}
    \begin{gathered}
      (p-1)\, \MC (t) \le \MC^A(t) \le
      \delta^{p-2}\,,
      \\
      (p-1)\, (\delta+A)^{p-2} \le \MC^A(t) \,.
    \end{gathered}
  \end{align}
\end{lemma}
\begin{proof}
  The statement is clear for $t \le A$ using
  $\MC^A(t)=\MC(t)=(\para +t)^{p-2} $, $0\le \para$, $t \le A$, and $p
  \le 2$. Moreover, $(\para+t)^{p-2}$ is a non-increasing function in
  $t$.

  For $t\ge A$ we
  have
  $\MC^A(t)=\function ''(A) + \frac{\function'(A)-\function''(A)\,
    A}{t}$. Thus, we get that $\MC^A(A)=(\delta+A)^{p-2}$,
  $\lim_{t\to
    \infty}\MC^A(t)=(\delta+A)^{p-3}\,\big(\delta+(p-1)A\big)$ and
  $(\MC^A)'(t)=- \frac{\function'(A)-\function''(A)\, A}{t^2}\le 0$ in
  view of \eqref{eq:E}, and $p\le 2$, hence proving that $a^{A}$ is
  non-increasing also for $t>A$ (contrary to the case $p>2$ we do not have
  any restriction on the choice of $A$).  This yields
  \begin{align*}
    (\para +A)^{p-2}\ge \MC^A(t) \ge (\para
    +A)^{p-3}((p-1)A+ \para) \ge (p-1)\, (\para +A)^{p-2}\,,
  \end{align*}
  which immediately implies the assertions using
  $\delta^{p-2} \ge (\para +A)^{p-2}$ and
  $(\para +A)^{p-2}\ge (\para +t)^{p-2}$ in view of $t\ge A$, and
  $p\le 2$.
\end{proof}
As in the case $p>2$  we would need
mainly the following corollary in the proof of regularity.
\begin{corollary}\label{cor:UA}
  Let the operator $\bS$, derived from the potential $\pot$, have
  $(p,\para)$-structure with $p\in (1,2]$ and $\delta> 0$. Denote
  $\function=\function_{p,\para}$, $\bF= \bF_\function$ and for
  $A\ge 1$ set ${\function^A:=\function^{A,2}}$,
  $\bF^A:= \bF_{\function^A}$ and $\bS^A:= \bS^{A,2}$.  Then, there
  holds for all $t\ge 0$ that
  \begin{align}
    \label{eq:UA3}
    \begin{gathered}
      (p-1)\, \function (t) \le \function^A(t) \le
      \frac {\para^{p-2}}2\, t^2\,,
      \\
      \frac{p-1}2\, (\para + A)^{p-2}\,t^2 \le \function^A(t) \,,
      \\
    (\function^A)^*(t)\le (p-1) \,(\Delta_2(\function^*))^M \,\function^*(t)\,,
    \end{gathered}
  \end{align}
  where $M\in \setN_0$ is chosen such that $(p-1)^{-1}\le 2^M$. Moreover, 
  for all $\bP \in \setR^{3\times 3}$ there holds
  \begin{align}
    \label{eq:UA4}
    \begin{aligned}
      \abs{\bF^A(\bP)}^2&\sim \function^A(\abs{\bP^{\sym}})\,,
      \\
    c\, \abs{\bF(\bP)}^2 &\le  \abs{\bF^A(\bP)}^2\,, 
    \\
    \abs{\bS^A(\bP)} &\le c\, \para^{p-2} \abs{\bP^{\sym}} \,,
    \end{aligned}
  \end{align}
  with constants $c$ depending only on $\gamma_3, \gamma_4$, and $p$. 
\end{corollary}
\begin{proof}
  Assertions \eqref{eq:UA3}$_{1,2}$ follow from \eqref{eq:wAl1}$_{1,2}$, the
  definition of $\MC$, $\MC^A$,
  $\function(0)=\function^A(0)=0$  and integration. Using the first inequality in
  \eqref{eq:UA3}$_1$ we get for all $t \ge 0$ that
  \begin{align*}
    (\function^A)^*(t)&= \sup _{s \ge 0} \ s\,t -\function^A(s)
    \\
                      &\le (p-1)\,\sup _{s \ge 0} \ s\,\frac t{p-1} -\function(s)
    \\
                      &= (p-1) \,\function^*\Big (\frac t{p-1}\Big )\le (p-1)\,
                        (\Delta_2(\function^*))^M\, \function^*(t)\,,
  \end{align*}
  with $M\in \setN_0 $ as chosen above. This proves
  \eqref{eq:UA3}$_2$. The inequalities in \eqref{eq:UA4} follow from
  Proposition \ref{prop:SA-ham} with $\bQ=\bfzero$, the definition of
  $\MC^A$, the fact that $\function $, $\function ^A$ are balanced,
  \eqref{eq:UA3}$_1$, the equivalences
  for $\bF$ and $\bS$ in Proposition~\ref{prop:hammer-phi}, and Lemma~\ref{lem:UA}.
 \end{proof}
 \section{On the existence and uniqueness of regular solutions}
\label{sec:regularity}
In this section we prove our main result, namely Theorem~\ref{thm:MT},
i.e., the existence and uniqueness of regular solutions
of~\eqref{eq:pfluid}, solely based on appropriate assumptions on the
regularity (but not on the size) of the data. From now on we will
restrict to the case $p>2$, but with a few (but non completely
trivial) changes the same arguments can be applied also to the case
$p \in (1,2]$, where a single approximation would be enough. Even if
the theory of approximation gives a unified approach valid for all
$p$, we decided to focus on the case $p>2$ since many estimates should
be changed, starting already from the a priori estimates and we think
that explaining the steps that need to be changed in the case $p\le 2$
would fragment the presentation in such a way that the readability of
the paper would be much more difficult. Since the result in the case
$p\le 2$ is already contained in \cite{SS00} and
\cite{br-plasticity,br-parabolic}, using a different approximation, we
preferred to skip them. Nevertheless, they will be presented in a
forthcoming paper~\cite{br-volume-chung}.

 \begin{definition}[Regular solution]
  Let the operator $\bS$ in~\eqref{eq:pfluid}, derived from a
  potential $\pot $, have $(p,\delta)$-structure for some
  $p\in(1,\infty)$, and $\delta\in[0, \infty)$ fixed but arbitrary.
  Let $\Omega\subset\setR^3$ be a bounded domain with $C^{2,1}$
  boundary, and let $I=(0,T)$, $T\in (0,\infty)$, be a finite time
  interval. Then, we say that $\bu$ is a regular solution of \eqref{eq:pfluid} if
  $\bu \in L^p(I;W^{1,p}_0 (\Omega))$ satisfies for all
  $\psi \in C_0^\infty (0,T)$ and all $\bw\in W^{1,p}_0(\Omega)$
\begin{equation*}
\begin{aligned}
  \int\limits _0^T\Bighskp{\frac {\partial\bu(t)}{\partial
      t}}{\bw}\,\psi(t)   + 
  \hskp{\bS(\bD\bu(t))}{\bD\bw}\,\psi (t)\,dt 
=
\int\limits
  _0^T\hskp{\bff(t)}{\bfw}\, \psi(t)\, dt\,,
\end{aligned}
\end{equation*}
and  fulfils
  \begin{align*}
    \begin{split}
      { \bfu}&\in L^{\infty}(I;W^{1,2}_{0}(\Omega))\cap {W^{1,2}(I; L^2(\Omega))}, 
     \\
      \bF(\bD\bfu)&\in L^{\infty}(I;L^{2}(\Omega))\cap L^{2}(I;W^{1,2}(\Omega)) ,
    \end{split}
  \end{align*}
 \end{definition}
 \begin{remark}
   Note that we are focusing on the ``natural'' second order
   spatial regularity, especially we are proving that 
   $\bF(\bD\bu)\in L^{2}(I;W^{1,2}(\Omega))$. In the parabolic case it
   is possible, at the price of some more restrictive hypotheses on
   the data, also to prove that
   $\bF(\bD\bu)\in W^{1,2} (I; L^{2} (\Omega))$. This result can be
   obtained independently on what we prove later on and nevertheless
   implies also some simplifications of the argument concerning the
   treatment of the time derivative. The regularity of
   $\frac{\partial\bF(\bD\bu)}{\partial t}$  would be
   needed in case of time-discretization to prove optimal convergence
   results, as done in \cite{br-parabolic}.
 \end{remark}
\begin{remark}\label{rem:ini}
To formulate clearly the dependence on the data in the various
estimates  we introduce the quantity
\begin{equation*}
  |||\bu_0,\bff|||^2:=
    \int\limits_\Omega 
    |\bu_0|^2 + 
    \function(|\bD\bu_0|)
    \,d\bx+\int\limits_0^T\int\limits_\Omega
    |\bff|^2
    \,d\bx \,dt\,.
\end{equation*}
Using the equivalences
$\function _{p,\para} (t) + \para ^p \sim t^p +\para ^p$ and
$\function ^*(t) \sim (\para ^{p-1} +t)^{p'-2}t^2$, valid for all
$p \in (1,\infty)$, $t,\para \ge 0$ with constants of equivalence just
depending on $p$, together with Korn and Poincar\`e inequalities, one
easily checks that $|||\bu_0,\bff||| $ is finite if
$\bu_0 \in W^{1,p}_0(\Omega)$ and
$\bff \in L^{2}(I\times \Omega)$, provided that $p\ge \frac 65$.
\end{remark}
We can now state the main result of this paper.

 \begin{theorem}
  \label{thm:MT}
  Let the operator $\bS$ in~\eqref{eq:pfluid}, derived from a
  potential $\pot $, have $(p,\delta)$-structure for some
  $p\in(2,\infty)$, and $\delta\in(0, \infty)$ fixed but arbitrary.
  Let $\Omega\subset\setR^3$ be a bounded domain with $C^{2,1}$
  boundary, and let $I=(0,T)$, $T\in (0,\infty)$, be a finite time
  interval. Assume that $\bu_0 \in W^{1,p}_0(\Omega)$
  and $\bff \in L^{2}(I\times \Omega)$.

Then, the system~\eqref{eq:pfluid} has a unique regular
solution 
  with norms  depending  only on the characteristics of
  $\bS$, $\delta^{-1}$, $T$,  $\Omega$, and $|||\bu_{0},\bff|||$.
\end{theorem}
To prove Theorem \ref{thm:MT} we use an approximate problem, obtained by replacing the
operator $\bS=\partial \pot$ with $(p,\delta)$-structure by the last
item $\bSn{N}= \partial \pot ^N$ of a special multiple
approximation $\bS^n$, $n=1,\ldots, N$, of $\bS$, i.e., 
\begin{equation*}
  \bSn{N}(\bP)=\frac{(\function^{N})'(|\bP^{\sym}|)}{|\bP^{\sym}|}\bP^{\sym}
  =\MC^{N}(|\bP^{\sym}|)\bP^{\sym}\,,
\end{equation*}
which we define now.
%
\begin{definition}[Special multiple approximation]\label{def:spec}
  Let the operator $\bS=\partial \pot$, derived from a potential
  $\pot$, have $(p,\delta)$-structure for some $p>2$ and $\para >0$
  with characteristics $(\gamma_3,\gamma_4,p)$. We call $\bS^n, \pot ^n, \function
^n$, and $\MC^n$, $n=1,\ldots , N$, a {\emph special} multiple
approximation of $\bS, \pot , \function_{p,\delta}$, and
$\MC_{\function_{p,\delta}}$ if it is a multiple approximation generated by  $N:=\lceil
  \frac{p-2}{2} \rceil$,  exponents
\begin{equation*}
q_n:=p-2n\quad\text{for } n=1, \ldots,
  N-1\qquad
\text{ and }\qquad    q_{N}:=2\,,
\end{equation*}
and parameters $A_n$, $n=1, \ldots, N$ satisfying the conditions
in Definition \ref{def:mult} and in Corollary \ref{cor:UAm}.
\end{definition}

\begin{remark}\label{rem:bal}
  (i) Let $\bS^n$, $n=1,\ldots, N$, be a special multiple
  approximation as in Definition \ref{def:spec}. Lemma
  \ref{lem:function} and a successive application of
  Proposition~\ref{prop:SA-struc}, and Lemma \ref{lem:d2} yields that
  for each $n=1,\ldots, N$ the operator $\bS^n$ has
  $\function^n$-structure with characteristics depending only on the
  characteristics of $\bS$, i.e., on $\gamma_3,\gamma_4,p$, due to the
  special choice of $q_j$, $j=1,\ldots, N$. The special choice of
  $q_j$, $j=1,\ldots, N$, Lemma~\ref{lem:function} and Lemma
  \ref{lem:eq12} imply that the characteristics of $\function^n$,
  $n=1,\ldots, N$, depends only on $p$. Thus, the constants in
  Proposition \ref{prop:pFA} as well as in Proposition
  \ref{prop:SA-ham} and Corollary \ref {cor:UAm} applied to $\bS^n$,
  $n=1,\ldots, N$, depend only on the characteristics of~$\bS$.

  (ii) In view of (i), Lemma \ref{lem:UAm}, and Remark \ref{rem:VA} the
  operator $\bS^N$ has $(2,\para)$-structure with characteristics
  depending on $p$, $\gamma_3, \gamma_4$, $\function ^{N-1}$, and
  $A_N$. 
\end{remark}

In view of the previous remark  we can work in the $W^{1,2}$-setting, which is sufficient
to justify all forthcoming computations, which is the main reason for
the introduction of these approximations.
\subsection{The approximate problem and some global regularity in 
  time} 
We have the following result on existence and uniqueness of
time-regular solutions of the approximate problem.
\begin{proposition}
\label{thm:existence_perturbation}
Let the operator $\bS=\partial U$, derived from the potential $\pot$,
have $(p,\delta)$-structure for some $p\in(2,\infty)$ and
$\delta\in(0,\infty)$. Assume that $\bu_0 \in W^{1,p}_0(\Omega)$ and
$\bff \in L^{2}(I\times \Omega)$. Let $\bSn{N} $ be the last item of a
special \footnote{Observe that the results of this section, prior to
  the passage to the limit, are in fact valid for any sequence of
  parameter $q_{n}$ as described in the definition of the multiple
  approximation.} multiple approximation $\bS^n$, $n=1,\ldots, N$, of
$\bS$ as in Definition~\ref{def:spec}. Then, the approximate problem
\begin{equation}
  \label{eq:eq-e}
  \begin{aligned}
    \frac{\partial\bue}{\partial t}-\divo \bSn{N} (\bfD\bue)&=\bff
\qquad&&\text{in }I\times\Omega\,,
    \\
    \bue &= \bfzero &&\text{on } I\times\partial \Omega\,,
    \\
    \bue(0)&=\bu_0&&\text{in }\Omega\,,
  \end{aligned}
\end{equation}
possesses a unique strong solution
$\bue$
, i.e., $\bue\in W^{1,2}(I;L^{2}(\Omega))$ with
$\bF^N(\bD\bue) \in L^{\infty}(I;L^{2}(\Omega))$, which
satisfies for all
$\psi \in C_0^\infty (0,T)$ and all $\bw\in W^{1,2}_0(\Omega)$
\begin{equation}
  \label{eq:weak-eps}
\begin{aligned}
  \int\limits _0^T\Bighskp{\frac {\partial\bue(t)}{\partial
      t}}{\bfw}\,\psi(t)   + 
  \hskp{\bSn{N}(\bD\bue(t))}{\bD\bfw}\,\psi (t)\,dt 
=
\int\limits
  _0^T\hskp{\bff(t)}{\bfw}\, \psi(t)\, dt\,.
\end{aligned}
\end{equation}
In addition, the solution $\bue$ satisfies 
the estimate
\begin{gather}
  \label{eq:main-apriori-estimate2}
  \begin{aligned}
    &\esssup _{t \in
      I}\|\bue(t)\|_{2}^{2}+\|\bFn{N}(\bD\bue(t))\|_{2}^{2} +\para
    ^{p-2}\norm{\nabla\bue(t)}_2^{2} 
    \\
    &\quad  +\para ^{p-2}\int\limits_{0}^{T}\norm{\bD\bue(t)}_2^{2}\,dt
    +\int\limits_{0}^{T}\Bignorm{\frac{\partial\bue(t)}{\partial
        t}}_2^2 \,dt
 \leq C\,\big (1 +|||\bu_0,\bff|||^2\big)\,,
  \end{aligned}
 \end{gather}
with $C$ depending only on the characteristics of $\bS$,
${\para^{p-2}}$, and
$\Omega$.
\end{proposition}
\begin{proof}
  The proof is based on a standard Faedo-Galerkin approximation of
  \eqref{eq:eq-e}. The existence of solutions of the Galerkin approximations
  follows from the standard Carath\'eodory theory for systems of ordinary differential
  equations. As pointed out in Remark \ref{rem:bal}, the operator $\bSn{N}$ has
  $(2,\delta)$ structure, hence the system can be treated essentially
  as the heat equation. In particular, once the existence of the
  Galerkin solution $\bue_{k}$, $k \in \setN$, is obtained, passing
  to the limit as $k\to \infty$ can be done within the standard theory
  of evolutionary problems with monotone operators.
Since this is a standard procedure, we just derive the a priori
estimates necessary for it. 

The first a priori estimate, derived by using $\bue _k$ as test
function in the Galerkin approximation for $\bue_{k} $, is the
following one:
\begin{equation*}
  \begin{aligned}
    \frac 12 \frac d{dt} \|\bue_k\|_2^2
    +c\, \|\bFn{N}(\bD\bue_{k})\|_2^{2} 
    & \leq c_{\epsilon}\int\limits_{\Omega}(\multappomega)^{*}(|\bff|)
    \,d\bx+\epsilon\int\limits_{\Omega}\multappomega(|\bue_{k}|)\,d\bx
    \\
    &\leq c_{\epsilon}\int\limits_{\Omega}(\multappomega)^{*}(|\bff|)
    \,d\bx+\epsilon\,
    C\int\limits_{\Omega}\multappomega(|\bD\bue_{k}|)\,d\bx,
  \end{aligned}
\end{equation*}
where we used in the first line Proposition~\ref{prop:SA-ham} with
$\bQ=\bfzero$ together with Young inequality,  and in the second line
\begin{equation*}
  \int\limits_{\Omega} \multappomega(|\bue_{k}|)\,d\bx\leq
  C_{P}\int\limits_{\Omega}\multappomega(|\nabla\bue_{k}|)\,d\bx\leq
  C_{P}C_{K}\int\limits_{\Omega}\multappomega(|\bD\bue_{k}|)\,d\bx, 
\end{equation*}
which follows from the modular versions of Poincar\`e and Korn
inequalities in Orlicz spaces (cf.~\cite{talenti, john,bdr-phi-stokes}). Moreover, we absorb the last term on 
the left-hand side of the previous estimate using
$\int_{\Omega}\multappomega(|\bD\bue_{k}|)\,d\bx \sim
\|\bFn{N}(\bD\bue_{k})\|_2^{2}$ in view of 
Corollary~\ref{cor:UAm}. Note that all constants are independent of
$A_{n}$, $n=1,\ldots, N$, and depend only on the characteristics of
$\bS$ and on $\Omega$, due to Remark \ref{rem:bal} and \cite{talenti, john,bdr-phi-stokes}.
Moreover, from 
Corollary \ref{cor:UAm} and the
definition of $|||\bu_{0},\bff|||$, it follows that
  \begin{equation*}
    \int\limits_{0}^{T}\int\limits_{\Omega}(\multappomega)^{*}(|\bff|)\,d\bx\,
    ds\leq c( p) \, \para ^{p-2}
    \, |||\bu_{0},\bff|||^{2}<\infty\,. 
  \end{equation*}
Hence, after the limiting procedure $k \to \infty$ we arrive at 
\begin{equation*}
  \esssup _{t \in
    I}\norm{\bue(t)}_2^{2}+\int\limits_{0}^{T}\norm{\bFn{N}(\bD\bue(s))}_2^{2}\,ds
  +\para ^{p-2}\int\limits_{0}^{T}\norm{\bD\bue(s)}_2^{2}\,ds \leq c \, |||\bu_{0},\bff|||^{2}\,,
\end{equation*}
where we also used Corollary \ref{cor:UAm}.  Next, we take
$\frac{\partial\bue_{k}}{\partial t}$ as test function in the Galerkin
approximation, use the fact that $\bSn{N}=\partial U^{N}$ is derived
from the potential $\pot^N \sim \function^n$ in view
of Lemma~\ref{lem:pot-equiA}, Lemma \ref{lem:trans} and
Proposition~\ref{prop:potential-equivalence}, to arrive at 
  \begin{equation*}
    \int\limits_{0}^{t}\Bignorm{\frac{\partial \bue_k(s)}{\partial
        t}}^{2}_{2}\,ds+\int\limits_{\Omega}\multappomega(|\bD\bue_k(t)|)\,d\bx\leq
    c\int\limits_{\Omega}\multappomega(|\bD\bu_{0}^k|)\,d\bx+c\int\limits_{0}^{t}\|\bff(s)\|_{2}^{2}\,ds\,.
  \end{equation*}
   Since $p>2$, we use Corollary \ref{cor:UAm}
   to arrive at 
  \begin{equation*}
     \int\limits_{\Omega}\multappomega(|\bD\bu_{0}^k|)\,d\bx \le c(p)\int\limits_{\Omega}\function(|\bD\bu_{0}^k|)\,d\bx \,.
   \end{equation*}
  These properties together with 
  Corollary \ref{cor:UAm} imply, after the limiting procedure $k \to
  \infty$, that for a.e. $t\in[0,T]$  
  \begin{equation*}
    \begin{aligned}
      \norm{\bF^N(\bD\bue(t))}_2^{2}  +\para ^{p-2}\norm{\nabla\bue(t)}_2^{2}
     + \int\limits_{0}^{T}\Bignorm{\frac{\partial \bue(s)}{\partial
        t}}^{2}_{2}\,ds 
      &\leq c\,\big (  1 + |||\bu_{0},\bff|||^{2}\,\big )\,,
    \end{aligned}
  \end{equation*}
where we also used Corollary \ref{cor:UAm} and Korn inequality. 

 
The uniqueness of the solution $\bue$ follows in a
standard manner. 
\end{proof}
\begin{remark}\label{rem:weak}
  Note that by the fundamental theorem of calculus of variations
  the weak formulation \eqref{eq:weak-eps} is equivalent to
\begin{align}\label{eq:weak-eps1}
  \Bighskp{\frac {\partial\bue (t)}{\partial t}
  }{\bfw} +\hskp{\bSn{N}(\bD\bue(t))}{\bD\bfw} =
  \hskp{\bff(t)}{\bfw}\,,
\end{align}  
being satisfied for $\textrm{a.e.}~t \in I$  and all $ \bw \in W^{1,2}_0(\Omega)$.
\end{remark}

\bigskip

In order to prove existence and uniqueness of regular solutions to~\eqref{eq:pfluid}, by
taking the various limits $A_{n}\to \infty$, we need to prove further regularity for the solution
$\bue$, namely on the second order spatial derivatives, in such a way that
$\bD\bue $ 
converges almost everywhere.
The regularity in the spatial variables requires an ad hoc treatment (localization) for the
Dirichlet boundary value problem.  To do this we adapt the argument introduced in
\cite{mnr3} (treating the case $p>2$) and that 
in~\cite{br-plasticity,br-parabolic} (treating the case $p<2$).  We sketch
the relevant steps, pointing out the main new aspects which are present in the
time-dependent case.
\subsection{Description and properties of the boundary}
\label{sec:bdr} 
$\hphantom{}$
We assume that the boundary $\partial\Omega$ is of class $C^{2,1}$, that
is for each point $P\in\partial\Omega$ there are local coordinates such
that in these coordinates we have $P=0$ and $\partial\Omega$ is locally
described by a $C^{2,1}$-function, i.e.,~there exist
$R_P,\,R'_P \in (0,\infty),\,r_P\in (0,1)$ and a $C^{2,1}$-function
$\Grenze_{P}:B_{R_P}^{2}(0)\to B_{R'_P}^1(0)$ such that
\begin{itemize}
\item   [\rm (b1) ] $\bx\in \partial\Omega\cap (B_{R_P}^{2}(0)\times
  B_{R'_P}^1(0))\ \Longleftrightarrow \ x_3=\Grenze_{P}(x_1,x_2)$,
\item   [\rm (b2) ] $\Omega_{P}:=\{(x',x_{3})\fdg x'=(x_1,x_2)
 \in  B_{R_P}^{2}(0),\ \Grenze_{P}(x')<x_3<\Grenze_{P}(x')+R'_P\}\subset \Omega$, 
\item [\rm (b3) ] $\nabla \Grenze_{P}(0)=\bfzero,\text{ and }\forall\,x'=(x_1,x_2)^\top
  \in B_{R_P}^{2}(0)\quad |\nabla \Grenze_{P}(x')|<r_P$,
\end{itemize}
where $B_{r}^k(0)$ denotes the $k$-dimensional open ball with center
$0$ and radius ${r>0}$.  Note that $r_P $ can be made arbitrarily
small if we make $R_P$ small enough.  In the sequel we will also use,
for $0<\lambda<1$, the  scaled open sets $\lambda\,
\Omega_P\subset \Omega_P$, defined as follows
\begin{equation*}
  \lambda\, \Omega_P:=\{(x',x_{3})\fdg x'=(x_1,x_2)^\top
 \in
  B_{\lambda R_P}^{2}(0),\ \Grenze_{P}(x')<x_3<\Grenze_{P}(x')+\lambda R_P'\}.
\end{equation*}
To localize near  $\partial\Omega\cap \partial\Omega_P$, for $P\in\partial\Omega$, we fix smooth
functions $\xi_{P}:\setR^{3}\to\setR$ such that 
\begin{itemize}
\item [$\rm (\ell 1)$] $\chi_{\frac{1}{2}\Omega_P}(\bx)\leq\xi_P(\bx)\leq
  \chi_{\frac{3}{4}\Omega_P}(\bx)$,
\end{itemize}
where $\chi_{A}(\bx)$ is the indicator function of the measurable set
$A$. 
For the remaining interior estimate we  localize by a smooth function
${0\leq\xi_{0}\leq 1}$ with $\spt \xi_{0}\subset\Omega_{0}$,
where $\Omega_{0}\subset \Omega$ is an open set such that
$\dist(\partial\Omega_{0},\,\partial\Omega)>0$.  
Since the boundary $\partial\Omega $ is compact, we can use an appropriate
finite sub-covering which, together with the interior estimate, yields
the global estimate.

Let us introduce the tangential derivatives near the boundary. To
simplify the notation we fix $P\in \partial\Omega$, $h\in (0,\frac{R_P}{16})$,
and simply write $\xi:=\xi_P$, $\Grenze:=\Grenze_{P}$. We use the standard notation
$\bx =(x',x_3)^\top$ and denote by $\be^i,i=1,2,3$ the canonical
orthonormal basis in $\setR^3$. In the following lower-case Greek
letters take values $1,\, 2$. For a function $f$ with $\spt
f\subset\spt\xi$ we define for $\alpha=1,2$ tangential translations:
\begin{equation*}
\begin{aligned}
  \trap f(x',x_3) = f_{\tau _\alpha}(x',x_3)&:=f\big (x' +
  h\,\be^\alpha,x_3+\Grenze(x'+h\,\be^\alpha)-\Grenze(x')\big )\,,
\end{aligned}
\end{equation*}
tangential differences $\Delta^+ f:=\trap f-f$, and tangential
difference quotients
%
  $\difp f:= h^{-1}\Delta^+ f$.  
%
For
simplicity we denote $\nabla \Grenze:=(\partial_1 \Grenze,\partial_{2}\Grenze, 0)^\top$
and use the operations $\trap {(\cdot)}$, $\tran {(\cdot)}$,
$\Delta^+(\cdot) $, $\Delta^+(\cdot) $, $\difp {(\cdot)}$ and $\difn
{(\cdot)}$ also for vector-valued and tensor-valued functions,
intended as acting component-wise.

We will use the following properties of the difference
quotients, all proved in~\cite{hugo-petr-rose}. 
\begin{lemma}
  \label{lem:TD1} 
  Let $\bv\in W^{1,1}(\Omega)$ such that $\spt \bv
  \subset\spt\xi$. Then 
\begin{equation*}
\begin{aligned}
  \nabla\difpm \bv &=\difpm{\nabla \bv }+\trap{(\partial_3 \bv
    )}\otimes\difpm{\nabla \Grenze},
  \\
  \bD\difpm \bv &=\difpm{\bD \bv }+\trap{(\partial_3 \bv
    )}\otimess\difpm{\nabla \Grenze},
  \\
  \diver\difpm \bv &=\difpm\diver \bv +\trapm{(\partial_3 \bv
    )}\difpm{\nabla \Grenze}
  \\
  \nabla \bv _{\pm\tau} &= (\nabla \bv )_{\pm\tau} + \trapm{(\partial_3 \bv
    )}\difpm{\nabla \Grenze},
\end{aligned}
\end{equation*}
where $\otimess$ is defined component-wise also for scalar and tensor-valued functions.
\end{lemma}

As for the classical difference quotients, $L^{q}$-uniform (with
respect to $h>0$) bounds for $\difp f$ imply that $\partial _\tau f $
belongs to $L^q(\spt\xi)$.
\begin{lemma}
  \label{lem:Dominic}
 It holds\footnote{Note that $\partial_{\tau}f$ denotes a
    tangential derivative, and to avoid confusion with time
    derivatives, the latter will be always denoted as
    $\frac{\partial f}{\partial t}$.}  that, if
  $ f \in W^{1,1}(\Omega)$, then we have for $\alpha=1,2$
  \begin{align}
    \label{eq:1}
    \difp f \to \td f=\partial _{\tau_\alpha}f :=\partial_\alpha f +\partial_\alpha
    \Grenze\, \partial_3 f  \qquad \text{ as } h\to 0,
  \end{align} 
  almost everywhere in $\spt\xi$, (cf.~\cite{mnr3}).
  If we
  define, for $0<h<R_P$
  \begin{equation*}
    \Omega_{P,h}=\left\{\bx\in \Omega_P\fdg x'\in B^2_{R_P-h}(0)\right\},
  \end{equation*}
  and if $f\in W^{1,q}_\loc(\Omega)$, $1\le q<\infty$, then
  \begin{equation*}
    \int\limits_{\Omega_{P,h}}|d^+f|^q\,d\bx\leq c\int\limits_{\Omega_{P}}|\partial_\tau f|^q\,d\bx.
  \end{equation*}
 Moreover, if $d^{+}f\in
L^q(\Omega_{P,h_0})$, $1< q<\infty$, and if
  \begin{equation*}
    \exists\,c_1>0:\quad   \int\limits_{\Omega_{P,h_0}}|d^{+}f|^q\,d\bx\leq c_1\qquad
    \forall\,h_0\in(0,R_P)\text{ and } \forall\,h\in(0,h_0),
  \end{equation*}
  then $\partial_{\tau}f\in L^q(\Omega_P)$ and
  \begin{equation*}
    \int\limits_{\Omega_{P}}|\partial_{\tau}f|^q\,d\bx\leq c_1.
  \end{equation*}
\end{lemma}
\begin{remark}\label{rem:orlicz}
  All assertions of the previous lemma also hold in Orlicz spaces
  generated by N-functions $\phi \in \Delta_2$, as can be easily seen
  by adapting the proof carried out in \cite{evans-pde} to this
  situation. 
\end{remark}
The following variant of formula of integration by parts will be often used.
\begin{lemma}
  \label{lem:TD3}
  Let $\spt g\cup\spt f\subset\spt\xi=\spt\xi_P$ and $0<h<\frac{R_P}{16}$. Then
  \begin{equation*}
    \intO f\tran g \, d\bx =\intO\trap f g\, d\bx.
  \end{equation*}
  Consequently, $\intO f\difp g \, d\bx= \intO(\difn f )g\, d\bx$.
  Moreover, if in addition $f$ and $g$ are smooth enough and at least
  one vanishes on $\partial\Omega$, then 
\begin{equation*} \intO f\td g \, d\bx= -\intO(\td f )g\,
    d\bx.
  \end{equation*}
\end{lemma}
Also the following properties of the difference quotient will be used
in the sequel. 
\begin{lemma}
  \label{lem:TD2}
  Let $\spt g \subset\spt\xi$. Then
\begin{equation*}
  \trap{(\difn g )}=-\difp g ,\quad \tran{(\difp g )}=-\difn g , \quad 
  \difn  g_\tau = - \difp g .
\end{equation*}
\end{lemma}
\begin{lemma}
 \label{lem:TD4}
 Let $\spt g\cup\spt f\subset\spt\xi$. Then
  \begin{equation*}
  \difpm (f g) = f_{\pm\tau} \,\difpm g + (\difpm f )\, g.
\end{equation*}
\end{lemma}
\subsection{A first regularity result in space}
%
We start proving spatial regularity for the approximate problem. %
The estimates proved in this intermediate step are uniform with
respect to $A_{n}$, $n=1,\ldots, N$, only  (a) in the interior of $\Omega$ and (b) in
the case of tangential derivatives. On the contrary estimates depend
on $A_{n}$ in the normal direction. Nevertheless, this allows later on
to use the equations point-wise to prove in a different way estimates
independent of $A_{n}$, $n=1,\ldots, N$, in the normal direction. Thus,
we can pass to the limit with $A_{n}\to\infty$, to treat the original
problem in the non-degenerate case.

We observe that by using a translation method, the result below is
proved rigorously for the solutions we constructed.
\begin{proposition}  
  \label{prop:JMAA2017-1}
  Let the operator $\bS=\partial \pot$, derived from the potential $\pot$, have
  $(p,\delta)$-structure for some $p\in(2,\infty)$, and
  $\delta\in(0,\infty)$   with characteristics $(\gamma_3,\gamma_4,p)$. Let $\Omega\subset\setR^3$ be a bounded
  domain with $C^{2,1}$ boundary and assume that
  $\bu_0 \in W^{1,p}_0(\Omega)$ and
  $\bff \in L^{2}(I\times \Omega)$. Let $\bSn{N} $ be
  the last item of the special multiple approximation $\bS^n$,
  $n=1,\ldots, N$, of $\bS$ from Definition \ref{def:spec}.
  Then, the unique strong solution
  $\bue$ of the approximate problem~\eqref{eq:eq-e} satisfies for
  a.e.~$t\in I$
  \begin{align}
    \begin{aligned}
      \int\limits_0^t \int\limits_{\Omega} \xi_0^2
      \abs{\nabla \bFn{N}(\bD\bue(s))}^2
      +
      \delta^{p-2} \xi_{0}^2|\nabla^{2}\bue(s)|^{2}\,d\bx\,ds
     &\le c_{0} \,,
      \\[3mm]
      \int\limits_0^t
      \int\limits_{\Omega} \xi^2_P \abs{\td
        \bFn{N}(\bD\bue(s))}^2+
     \delta^{p-2}       \xi_{P}^2|\td\nabla\bue(s)|^{2}\,d\bx\,  ds
      &
      \le c_{P}\,,
        \end{aligned} \label{eq:est-eps}
  \end{align}
  where $c_{0}=c_{0} (\para
  ^{2-p}\|\bff\|_2^2\,,|||\bu_0,\bff|||,\norm{\xi_0}_{ 1,\infty},\gamma_3,\gamma_4,p)$, while
  the constant related to the neighborhood of $P$ is such that 
  $c_{P}=c_{P} (\para
  ^{2-p}\|\bff\|_2^2\,,|||\bu_0,\bff|||,\norm{\xi_P}_{1,\infty},\norm{\Grenze_{P}}_{C^{2,1}},$
  $\gamma_3,\gamma_4,p)$. 
  Here, $\xi_{0}(\bx)$ is a cut-off function with support in the interior of
  $\Omega$ and, for arbitrary $P\in \partial \Omega$, the tangential
  derivative 
  is defined locally in $\Omega_P$ via~\eqref{eq:1}. 
\end{proposition}
  Proposition~\ref{prop:JMAA2017-1} and Proposition~\ref{thm:existence_perturbation} imply 
  $\bue(t) \in W^{2,2}(\Omega)$ and ${\frac{\partial \bue}{\partial
    t}(t)\in L^2(\Omega)}$ for a.e.~$t \in I$.
  Hence, equations~\eqref{eq:eq-e} hold point-wise a.e.~in~$I\times \Omega$.

  We employ this to  deduce the
following result, by using the equations in a point-wise sense,
yielding however a critical dependence on the approximation of the operator.
\begin{proposition}  
  \label{prop:JMAA2017-2}
  Under the assumptions of Proposition \ref{prop:JMAA2017-1} there exists a constant
  $C_1>0$ such that, 
  provided  in the local description of the boundary there holds $r_P<C_1$ in $(b3)$,
  where $\xi_{P}$ is a cut-off function with support in
  $\Omega_P$, then there holds for a.e.~$t\in I$
  \begin{equation}
    \label{eq:est-eps-1}
    \begin{aligned}
      \int\limits_0^t \int\limits_{\Omega} 
      \xi^2_{P} \abs{\partial _3 \bFn{N}(\bD\bue(s))}^2
      +\delta^{2-p}\xi^2_{P}    \abs{\partial_3 \bD\bue(s)}^2\,d\bx\,ds
      \le C_{N}\,,
    \end{aligned}
  \end{equation}
  where $C_{N}=C_{N}(\delta^{2-p}, \para
    ^{2-p}\|\bff\|_2^2\,,|||\bu_0,\bff|||,\norm{\xi_P}_{1,\infty},\norm{\Grenze_{P}}_{C^{2,1}},\gamma{}_3,\gamma{}_4,p,
  A_N,\function_{N-1})$.
\end{proposition}
\begin{remark}
  We consider only the case $\delta>0$ and in the estimates of the two above propositions
  all dependencies on $\delta$ are traced in a precise and explicit way, showing how they
  deteriorate in the degenerate case. The degenerate problem could be treated by assuming
  more stringent assumptions on the regularity of the data (namely the regularity of the
  right-hand side $\bff$). The same phenomenon is well-known to happen even for the
  $p$-Laplace problem.  In that case sharpness of additional assumptions and links
  with the fractional regularity of the solution are proved and
  discussed in detail by Brasco and
  Santambrogio~\cite{BS2018} and the references therein.
\end{remark}

\begin{proof}[Proof of Proposition~\ref{prop:JMAA2017-1}]
Fix $P\in \partial \Omega$ and define in $\Omega_P$
  \begin{equation*}
    \bw:=\difn{(\xi^2\difp(\bue |_{\frac 12 {\Omega}_P}))},
  \end{equation*}
  where $\xi:=\xi_P$, $\Grenze:=\Grenze_{P}$, and
  $h\in(0,\frac{R_{P}}{16})$ and use the function $\bw$ extended by
  zero outside of $\Omega_{P}$ as a test function
  in~\eqref{eq:weak-eps1}. This yields, using the properties of the
  difference quotient in Lemma \ref{lem:TD1}, Lemma \ref{lem:TD3},
  Lemma \ref{lem:TD4}, for a.e.~$t \in I$
\allowdisplaybreaks
\enlargethispage{5mm}
\begin{align*}
     &\int\limits_{0}^{t}\int\limits_{\Omega}{\xi^2\difp{\frac{\partial\bue(s)}{\partial t}}}\cdot{ \difp\bue (s)}\, d\bx \,ds + \int\limits_{0}^{t}\intO
    \xi^2\difp{\bSn{N}(\bD\bue (s))}\cdot \difp \bD\bue  (s)\, d\bx \,ds
    \\
    &=-\int\limits_{0}^{t}\intO \xi^{2}\difp{\bSn{N}}(\bD\bue (s))\cdot\big(\partial _3 \bue (s)\big)_{\tau} \otimess\difp\nabla \Grenze\, d\bx\, ds
    \\
    &\quad -2\int\limits_{0}^{t}\intO\difp {\bSn{N}}(\bD\bue (s))\cdot \xi\nabla \xi \otimess\difp\bue\, d\bx\, ds
    \\
    &\quad +\int\limits_{0}^{t}\intO \bSn{N}(\trap{(\bD\bue)})\cdot \big(2 \xi
    \partial_3\xi \difp\bue  \big)\otimess\difp\nabla \Grenze \,
    d\bx \,ds 
    \\
    &\quad +\int\limits_{0}^{t}\intO \bSn{N}(\trap{(\bD\bue)})\cdot \big(\xi^2
    \difp\partial_3\bue  \big)\otimess\difp\nabla \Grenze \, d\bx \,ds 
    \\
    &\quad +\int\limits_{0}^{t}\intO\ff (s)\cdot\difn(\xi^2 \difp \bue (s))\, d\bx\, ds
      =:\sum_{j=1}^{5}\int\limits_{0}^{t}\mathcal{I}_j(s)\,ds\,.
\end{align*}
Proposition \ref{prop:SA-ham} yields for a.e.~$s \in I$ the following equivalence
\begin{equation*}
\intO \xi^2 \bigabs{ \difp{\bFn{N}}  (\bD\bue(s)) }^2 \, d\bx \sim\intO
    \xi^2\difp{\bSn{N}(\bD\bue (s))}\cdot \difp \bD\bue  (s)\, d\bx\,,
\end{equation*}
with constants depending only on the characteristics of $\bS$, due to
Remark \ref{rem:bal}. This equivalence provides the
``natural"{} quantity on the left-hand side.  We estimate the integrals
$\mathcal{I}_j$, $j=1,\ldots, 5$, similarly as in \cite{hugo-petr-rose}. Note that all
constants in the following can depend on the characteristics of $\bS$
and that other dependencies will be indicated.

We start estimating the first one as 
\begin{equation*}
  \begin{aligned}
    \mathcal{I}_{1}&\leq c \int\limits_\Omega
    \xi^{2}|\difp{\bD\bue}|\,\MC^{N}(|\bD\bue|+|\Delta^+\bD\bue|)
    |(\nabla \bue)_\tau|\,|\difp \nabla \Grenze|\, d\bx  
    \\
    &
    \leq  c\, \norm{\Grenze}_{C^{1,1}}\left(\int\limits_\Omega
      \xi^{2}\,\MC^{N}(|\bD\bue|+|\Delta^+\bD\bue|)|\difp{\bD\bue}|^{2}\, d\bx
    \right)^{1/2}\times
    \\
    &\qquad \times \left(\int\limits_{\Omega}
      \xi^{2}\MC^{N}(|\bD\bue|+|\Delta^+\bD\bue|)|(\nabla
      \bue)_\tau|^{2}\, d\bx\right)^{1/2} 
    \\
    &\leq     \epsilon\,\norm{\xi \,\difp{\bFn{N}}(\bD\bue)}_2^{2}+C\, \Big (\delta^{p} + 
    \int\limits_\Omega
    \multappomega(|\bD\bue|)\,d\bx\Big ),
  \end{aligned}
\end{equation*}
where we used Proposition \ref{prop:SA-ham}, H\"older and Young
inequalities, Lemma~\ref{lem:ast}, the convexity and
$\Delta_2$-condition of the balanced N-function $\function ^N$, the
substitution theorem and Korn inequality. The constant $C$ depends on
$\norm{   \Grenze}_{C^{1,1}}$ and $\epsilon^{-1}$. 
Note that in view of $\int\limits_{\Omega}\multappomega(|\bD\bue|)\,d\bx \sim
\|\bFn{N}(\bD\bue)\|_2^{2}$, estimate~\eqref{eq:main-apriori-estimate2},
and the substitution theorem the right-hand
side of the last estimate is finite. This comment also applies to the estimates of the other terms
$\mathcal I_j$, $j=1,\ldots, 5$.

The second term is estimated more or less in the same way 
\begin{align*}
    \mathcal{I}_{2}&\leq c \int\limits_\Omega
    \xi|\difp{\bD\bue}|\,\MC^{N}(|\bD\bue|+|\Delta^+\bD\bue|)|\nabla\xi| |\difp \bue
   | \, d\bx 
    \\
    &
    \leq c\,\norm{\nabla\xi}^{2}_{\infty}
\bigg(\int\limits_\Omega
      \xi^{2}\,\MC^{N}(|\bD\bue|+|\Delta^+\bD\bue|)|\difp{\bD\bue}|^{2}\,d\bx
    \bigg)^{1/2}\times
    \\
    &\qquad \times \bigg(\,\int\limits_{\Omega\cap \spt \xi} \MC^{N}(|\bD\bue|+|\Delta^+\bD\bue|)|\difp \bue|^{2}\, d\bx\bigg)^{1/2}
    \\
    &\leq
    \epsilon\,\norm{\xi\,\difp{\bFn{N}}(\bD\bue)}^2+C(\epsilon^{-1},
   \norm{\xi}_{1,\infty}) \,\Big (\delta^{p} + 
    \int\limits_\Omega
    \multappomega(|\bD\bue|)\,d\bx\Big ),
\end{align*}
where we additionally used Remark \ref{rem:orlicz}.

To estimate the integral $\mathcal{I}_3$ we use that due to
Proposition \ref{prop:SA-ham} there holds
$\abs{\bSn{N}(\bP)} \le c\,
(\multappomega)'(\abs{\bP^{\sym}})$. Using this,  
Young inequality, \eqref{eq:*'}, the substitution theorem,
Remark \ref{rem:orlicz} and Korn inequality we get
\begin{equation*}
  \begin{aligned}
    \abs{\mathcal{I}_{3}}&\leq c(\norm{\xi}_{1,\infty}, \|\Grenze\|_{C^{2,1}})\intO
    (\multappomega)^{*}(|\bSn{N}(\trap{(\bD\bue )})|) + \multappomega(| \difp\bue| )\,
    d\bx
    \\
    & \leq C(\norm{\xi}_{1,\infty}, \|\Grenze \|_{C^{2,1}})\intO \multappomega(| \bD\bue |
    )\, d\bx.
\end{aligned}
\end{equation*}

The integral $\mathcal I _4$ is estimated by using Lemma
\ref{lem:TD3}--Lemma \ref{lem:TD4} to obtain
\begin{equation*}
  \begin{aligned}
    |\mathcal{I}_{4}| 
    &=\Big|\intO
    \big(-\xi^2\difp{\bSn{N}}(\bD\bue)\difp\nabla\Grenze+\bSn{N}(\bD\bue)\difp(
    \nabla\Grenze) \,\difn(\xi^2 )\big)\otimess\partial_{3}\bue
    \\
    &\qquad\;\; +\bSn{N}(\bD\bue)(\xi^2)_{-\tau}\,\difn {\difp \nabla \Grenze}
    \otimess\partial_{3}\bue \, d\bx\,\Big|
    \\
    &\leq \epsilon\,\norm{\xi\,\difp{\bFn{N}}(\bD\bue)}_2^2+c(\vep^{-1}, \norm{\xi}_{1,\infty}, \norm{\Grenze}_{C^{2,1}})\, \Big (
    \delta^{p}+\int\limits_\Omega \multappomega(|\bD\bue|)\,d\bx\Big ),
  \end{aligned}
\end{equation*}
where the first term was treated as $\mathcal I_1$, while the other
two were treated as $\mathcal I_3$ .

On the other hand, the integral related to the right-hand side can be
estimated as follows 
\begin{align*}
    \mathcal{I}_{5}&\leq c(\vep ^{-1})\, \para ^{2-p}\norm{\bff}_2^{2} +
    \epsilon\, \para ^{p-2}\int\limits_{\Omega}|\difn(\xi^2 \difp \bue )|^2\,d\bx  
    \\
    & \leq c(\vep ^{-1})\, \para ^{2-p}\norm{\bff}_2^{2} +
     c(\norm{\xi}_{1,\infty}, \norm{\Grenze}_{C^{1,1}}) \, \para
     ^{p-2}\int\limits_{\Omega}|\bD\bue |^2\,d\bx
     \\
     &\quad+\epsilon\,
    c\, \para ^{p-2}\int\limits_{\Omega}\xi^{2}\,|\difp \bD \bue |^2\,d\bx
    \\
    &\leq c(\vep ^{-1})\, \para ^{2-p}\norm{\bff}_2^{2} +
     c(\norm{\xi}_{1,\infty}, \norm{\Grenze}_{C^{1,1}}) \, \para
     ^{p-2}\int\limits_{\Omega}|\bD\bue |^2\,d\bx
     \\
     &\quad+\epsilon\,
     c\, \int\limits_{\Omega}|\difp{\bFn{N}}(\bD\bue)|^2 \,d\bx,
\end{align*}
where we used standard properties of the difference quotient in $L^2$,
Korn inequality, the substitution theorem, as well as Proposition~\ref{prop:SA-ham}, and
Lemma~\ref{lem:UAm}, which yield
$\para^{p-2}|\difp \bD \bue |^2 \le c\, |\difp{\bFn{N}}(\bD\bue)|^2$.

Observing that
$d^+\frac{\partial\bue}{\partial t}=\frac{\partial d^+\bue}{\partial
  t}$, choosing $\vep>0$ sufficiently small, and using
\eqref{eq:main-apriori-estimate2} we proved that for a.e.~$t\in I$
\begin{align}
  \begin{aligned}\label{eq:Ft}
    &\frac{1}{2}\intO\xi^2|\difp\bue(t)|^2\,d\bx
    +c\int\limits_{0}^{t}\intO \xi^2 \bigabs{ \difp{\bFn{N}}
      (\bD\bue(s)) }^2 \, d\bx \,ds
    \\
    &\le \frac{1}{2}\intO\xi^2|\difp\bu_{0}|^2\,d\bx+
    c(\norm{\xi}_{1,\infty},\norm{\Grenze}_{C^{2,1}},{\para
      ^{2-p}}\|\bff\|^{2}_2\,,|||\bu_0,\bff|||, \gamma_3,\gamma_4,p)
    \\
    &\le C_0 (\norm{\xi}_{ 1,\infty},\norm{\Grenze}_{C^{2,1}},{\para
      ^{2-p}}\|\bff\|^{2}_2\,,|||\bu_0,\bff|||,
    \gamma_3,\gamma_4,p)\,,
  \end{aligned}
\end{align}
where we also used the assumption on the data. Since $C_0$ does not depend 
on $h>0$, it follows by Lemma~\ref{lem:Dominic}
that for a.e.~$t \in I$
\begin{align*}
  \int\limits_{0}^{t}\intO \xi^2 \bigabs{ \td \bFn{N} (\bD\bue(s)) }^2
  \, d\bx\, ds
  &\leq \int\limits_{0}^{t}\intO \xi^2 \bigabs{
    \difp{\bFn{N}}(\bD\bue(s)) }^2 \, d\bx\, ds \le C_0\,,
\end{align*}
proving the estimate for the first term in \eqref{eq:est-eps}$_2$.
Next, observe that Proposition~\ref{prop:SA-ham} and Lemma
\ref{lem:UAm} imply
\begin{equation*}
  \delta^{p-2}
  \int\limits_{0}^{t}\int\limits_{\Omega}\xi^{2}|\difp\bD\bue(s)|^{2}\,d\bx\,
  ds \leq \int\limits_{0}^{t}\intO \xi^2 \bigabs{
    \difp{\bFn{N}} (\bD\bue(s)) }^2 \, d\bx \, ds\leq C_0\,.
\end{equation*}
Now we proceed exactly as in the proof of
\cite[(3.12)--(3.14)]{br-reg-shearthin} for the special choice
$\phi(t)=t^2$ to get
\begin{align}\label{eq:l2}
  &\delta^{p-2}  \!\int\limits_{0}^{t}\int\limits_{\Omega}\xi^{2}|\difp
  \nabla\bue(s)|^{2}\,d\bx \, ds
  \\
  &\!\!\leq   \delta^{p-2}\!
    \int\limits_{0}^{t}\int\limits_{\Omega}\xi^{2}|\difp\bD\bue(s)|^{2}\,d\bx\,
    ds + c(\norm{\xi}_{1,\infty},\norm{\Grenze}_{C^{1,1}})\, \delta^{p-2}
    \int\limits_{0}^{t}\int\limits_{\Omega}|\bD\bue(s)|^{2}\,d\bx\,
    ds \notag
  \\
  &\!\!\le C_0 + c(\norm{\xi}_{1,\infty},\norm{\Grenze}_{C^{1,1}})\, \delta^{p-2}
    \int\limits_{0}^{t}\int\limits_{\Omega}|\bD\bue(s)|^{2}\,d\bx\,
    ds\,.\notag
\end{align}
This, the a priori estimate \eqref{eq:main-apriori-estimate2}, and
Lemma \ref{lem:Dominic} finally shows for a.e.~$t \in I$
\begin{equation*}
  \delta^{p-2}  \int\limits_{0}^{t}\int\limits_{\Omega}\xi^{2}|\td\nabla\bue|^{2}\,d\bx\leq
  C (\norm{\xi}_{1,\infty},\norm{\Grenze}_{C^{2,1}},{\para
    ^{2-p}}\|\bff\|^{2}\,,|||\bu_0,\bff|||, \gamma_3,\gamma_4,p)\,, 
\end{equation*}
proving the estimate for the second term in \eqref{eq:est-eps}$_2$.
%

\smallskip

The same argument used with a test function $\xi_{0}$ with compact
support in $\Omega$, and standard difference quotients can be used to
prove~$(\ref{eq:est-eps})_1$.
\end{proof}

\begin{corollary}\label{cor:Ftau}
  Under the assumptions of Proposition \ref{prop:JMAA2017-1} there holds
  a.e.~in $I\times \Omega$
  \begin{align*}
    \bigabs{\td\bF^N(\bD\bue)}^2 \sim \MC^N(\abs{\bD\bue}) \abs{\td\bD\bue}^2
  \end{align*}
  with constants depending only on the characteristics of $\bS$.
\end{corollary}
\begin{proof}
  Proposition \ref{prop:SA-ham} implies
  \begin{align*}
    \bigabs{\difp\bF^N(\bD\bue)}^2 \sim
    \MC^N(\abs{\bD\bue}+\abs{\Delta\bD\bue}) \abs{\difp\bD\bue}^2\,. 
  \end{align*}
  The estimates \eqref{eq:Ft}, \eqref{eq:l2} and Lemma
  \ref{lem:Dominic} yield that a.e.~in $I\times \Omega$ there holds
  $\difp\bF^N(\bD\bue) \to \td\bF^N(\bD\bue)$ and
  $\difp\bD\bue \to \td\bD\bue$ as $h\to 0$. These observations immediately imply
  the assertion.
\end{proof}

Now we prove the result on the regularity in the ``normal"{} direction
from~\eqref{eq:est-eps-1}, which is valid up to the boundary, but is 
dependent on the chosen multiple approximation.
\begin{proof}[Proof of Proposition~\ref{prop:JMAA2017-2}]
  Thanks to the previous results we can re-write the equations
in~\eqref{eq:eq-e} a.e.~in $I\times \Omega$ as follows
  \begin{equation*}
    \label{eq:linear_system}
    -\frac {\partial u^{N}_{i}}{\partial t}   + \sum_{k=1}^3\partial_{k 3}\Sn{N}_{i 3}(\bD\bue)\partial_3
    D_{k 3}\bue + \sum_{\alpha=1}^2\partial_{3\alpha}\Sn{N}_{i 3}(\bD\bue)\partial_3
    D_{ 3\alpha}\bue  =\mathfrak{f}_{i} \,,
  \end{equation*}
  where 
  \begin{equation*}
\mathfrak {f}_{i} :=-f_{i}-\sum_{\gamma,\sigma=1}^2
\partial_{\gamma \sigma}\Sn{N}_{i 3}(\bD\bue)\partial_3
  D_{\gamma \sigma}\bue- \sum_{k,l=1}^3\partial_{k l}\Sn{N}_{i \beta}(\bD\bue)\partial_\beta
  D_{k l}\bue\,,
\end{equation*}
for $i=1,2,3$.  We now proceed as in~\cite[Eq.~(3.3)]{br-plasticity}
and multiply these equations by $\partial _3 \widehat D_{i 3}\bue$,
where $\widehat D_{\alpha\beta}\bue =0$, for $\alpha,\beta=1,2$,
$\widehat D_{\alpha 3}\bue =\widehat D_{3\alpha}\bue
=2D_{\alpha3}\bue$, for $\alpha=1,2$,
$\widehat D_{33}\bue =D_{33}\bue $ and sum over $i=1,2,3$. Since
$\bS^N$ has $\function^N$-structure we get 
  \begin{equation*}
    -\sum_{i=1}^3\frac {\partial u^{N}_{i}}{\partial t}\partial_3\widehat{D}_{i3}\bue+
    \gamma\,\MC^{N}(|\bD\bue|)| {\boldsymbol { \mathfrak b}}|^2 \leq
     |\boldsymbol { \mathfrak f}|| {\boldsymbol {   \mathfrak b}}|\qquad\textrm {a.e.\ in
    }I\times \Omega\,,
  \end{equation*}
where $\mathfrak b_i :=\partial_3 D_{i3}\bue$ and where the constant
$\gamma$ just depends on the characteristics of $\bS$. 

  By straightforward manipulations (cf.~\cite[Sections 3.2~and~4.2]{br-reg-shearthin}) we
  obtain that a.e.~in $I\times \Omega_P$ it holds
  \begin{equation*}
    \begin{aligned}
      |\boldsymbol { \mathfrak f}| &\leq c \left(|\bff|
        \MC^{N}(|\bD\bue|\right )\left(|\partial_\tau\nabla
          \bue|+\|\nabla\Grenze\|_{\infty}|\nabla^2 \bue|\right)\,,
      \\
      |\boldsymbol { \mathfrak b}|&\geq2|\widetilde{\boldsymbol{\mathfrak
          b}}|-|\partial_{\tau}\nabla\bue|-\|\nabla
      \Grenze\|_{\infty}|\nabla ^{2}\bue|\,,
    \end{aligned}
  \end{equation*}
  for $\widetilde {\mathfrak b}_i:=\partial ^2_{33}u^{N}_{i}$,
  $i=1,2,3$. Consequently we get a.e.~in $I\times \Omega_P$
  \begin{align*}
    -&\sum_{i=1}^3\frac {\partial u^{N}_{i}}{\partial t}\partial_3\widehat{D}_{i3}\bue+ 
       2\gamma\,\MC^{N}(|\bD\bue|) |{\widetilde {\boldsymbol { \mathfrak b}}}|^2
    \\
    & \leq c \left[ \,|\bff|+\MC^{N}(|\bD\bue|)\left(|\partial_\tau\nabla
        \bue|+\|\nabla \Grenze\|_{\infty}|\nabla^2 \bue|\right)\right]  |{{\boldsymbol {
        \mathfrak b}}}|\,.
  \end{align*}
  We then add on both sides, for
  $\alpha=1,2$ and $i,k=1,2,3$, the term (which is finite a.e.) 
  \begin{equation*}
   2\gamma\, \MC^{N}(|\bD\bue|)\,|\partial_\alpha\partial_i
    u^{N}_{k}|^2\,,
  \end{equation*}
  use the estimate $|\boldsymbol{\mathfrak{b}}|\le |\nabla ^2 \bue|$
  and Young inequality, yielding
  \begin{equation*}
    \begin{aligned}
      &-\sum_{i=1}^3\frac{\partial u^{N}_{i}}{\partial
        t}\partial_3\widehat{D}_{i3}\bue+ 2\gamma\, \MC^{N}(|\bD\bue|)
      | {\nabla ^2\bue }|^2
      \\
      & \leq {\gamma}\MC^{N}(|\bD\bue|)|\nabla^2 \bue|^2
      \!+\!\frac{c\,|\bff|^2}{\MC^{N}\!(|\bD\bue|)} \!+\!c\,\MC^{N}\!(|\bD\bue|)
      \big(|\partial_\tau\nabla \bue|^2\!+\!\|\nabla
      \Grenze\|_{\infty}^2|\nabla^2 \bue|^2\big),
    \end{aligned}
  \end{equation*}
  where in the right-hand side we used also the definition of the tangential derivative
  (cf.~\eqref{eq:1}).  Next, we choose the sets $\Omega_{P}$  such that
  $\|\nabla \Grenze\|_{\infty}=\|\nabla \Grenze_{P}(x_1,x_2)\|_{{\infty},\Omega_{P}}$ is small enough,
  so that we can absorb the last term from the right-hand side. We finally arrive at the
  following pointwise inequality
  \begin{equation}
   \label{eq:pointwise-estimate}
   \begin{aligned}
      &-\sum_{i=1}^3\frac {\partial u^{N}_{i}}{\partial
        t}\partial_3\widehat{D}_{i3}\bue + \gamma\,\MC^{N}(|\bD\bue|)
      | {\nabla ^2\bue }|^2 
      \\
      &
      \leq c \left(\frac{|\bff|^2}{\MC^{N}(|\bD\bue|)}+\MC^{N}(|\bD\bue|)\,|\partial_\tau\nabla
          \bue|^2\right)\qquad \text{a.e. in }I\times \Omega_P\,.
    \end{aligned}
  \end{equation}
  We 
  multiply~\eqref{eq:pointwise-estimate} by $\xi^2$, and integrate for
  a.e.~$t \in I$  over the sub-domain
  \begin{equation*}
    (0,t)\times \Omega_{P,\epsilon}:=(0,t)\times\left\{\bx\in \Omega_{P}\fdg \Grenze_{P}+\epsilon<x_{3}<\Grenze_{P}+R_{P}'\right\}\,,
  \end{equation*}
  for $0<\epsilon<R_{P}'$.  This shows, using also Young inequality, that
  \begin{align*}
    &
      \gamma \int\limits_{0}^{t}\int\limits_{\Omega_{P,\epsilon}}
    \xi^{2} \MC^{N}(|\bD\bue|) | {\nabla ^2\bue }|^2d\bx \,ds
    \\
    & \leq \int\limits_{0}^{t}\int\limits_{\Omega_{P,\epsilon}}
    \!\!\!\! c\,\xi^{2} \!
    \left(\frac{|\bff|^2}{\MC^{N}(|\bD\bue|)}+\MC^{N}(|\bD\bue|)\,|\partial_\tau\nabla
      \bue|^2\right)+ \xi^{2}\left|\frac {\partial \bue}{\partial
        t}\right||\nabla^{2}\bue|\,d\bx \,ds
    \\
    & \leq \int\limits_{0}^{t}\intO c\,\xi^{2} \left(\frac{|\bff|^2
        +\big|\frac {\partial \bue}{\partial
          t}\big|^{2}}{\MC^{N}(|\bD\bue|)}+\MC^{N}(|\bD\bue|)\,|\partial_\tau\nabla
      \bue|^2\right)d\bx \,ds
    \\
    &\quad +\frac{\gamma}{2}
    \int\limits_{0}^{t}\int\limits_{\Omega_{P,\epsilon}}
    \xi^{2}\MC^{N}(|\bD\bue|)|\nabla^{2}\bue|^2\,d\bx\, ds\,.
  \end{align*}
  Now we absorb the last term from the right-hand side in the
  left-hand side. Moreover, we use that $\MC^{N}$ is bounded from
  below by $c\, \para ^{p-2}$ (cf.~Lemma \ref{lem:UAm}), the assumption on $\bff$ and
  \eqref{eq:main-apriori-estimate2} to estimate the first term on the
  right-hand side. To handle the second term we first use that
  $\MC ^N$ is bounded from above by a constant $c$ depending on $p$,
  $\gamma_3, \gamma_4$, $\function ^{N-1}$, and $A_N$ (cf.~Remark
  \ref{rem:VA}, Remark \ref{rem:bal}) and then we use the estimate
  \eqref{eq:est-eps}$_2$. These estimates result in 
  \begin{equation*}
         \int\limits_{0}^{T}\int\limits_{\Omega_{P,\epsilon}} \xi^{2}
         \MC^{N}(|\bD\bue|)|\nabla^{2}\bue|^{2}\,d\bx\, dt\leq
         C(A_N,\function_{N-1},\para ^{2-p},|||\bu_{0},\bff|||)\,.
  \end{equation*}
  By monotone convergence for $\vep \to 0$, this shows that
  \begin{equation*}
    \int\limits_{0}^{T}\int\limits_\Omega \xi^{2}
    \MC^{N}(|\bD\bue|)|\nabla^{2}\bue|^{2}\,d\bx\, dt\leq
    C(A_N,\function_{N-1},\para ^{2-p}, c_P 
    )\,.
      \end{equation*}
  Using Lemma \ref{lem:UAm} we finally get also 
    \begin{equation*}
    \para ^{p-2}\int\limits_{0}^{T}\int\limits_\Omega \xi^{2}
    |\nabla^{2}\bue|^{2}\,d\bx\, dt\leq    C(A_N,\function_{N-1},\para
    ^{2-p}, c_P 
    )\,.
    \end{equation*}
The last two estimates together with Proposition \ref{prop:pFA} and
the definition of the tangential derivatives (cf.~\eqref{eq:1}) 
finish the proof of Proposition \ref{prop:JMAA2017-2}.
\end{proof}
\subsection{Uniform estimates for the second order spatial derivatives}
We now improve the estimate in the normal direction in the sense that
we will show that they are bounded uniformly with respect to the
parameters $A_{n}$, for all $n=1,\dots,N$. The used method is an
adaption to the time evolution problem and the case $p>2$ of the
treatment in~\cite{br-plasticity,SS00} in the case $p<2$
(cf.~\cite{br-parabolic}).  In particular, it involves a technical
steps to justify the treatment of the time derivative, 
which is an adaptation of the method used in \cite{LM72,tem}.
\begin{lemma}
  \label{lem:time-derivative}
    Let $\partial\Omega\in C^{2,1}$ and let
  $\bv\in L^{2}(I;W^{2,2}(\Omega)\cap W^{1,2}_{0}(\Omega))\cap
  W^{1,2}(I;L^{2}(\Omega))$. Then, for all $t\in[0,T]$ it holds
  \begin{equation*}
  -\int\limits_{0}^{t}\int\limits_{\Omega}\frac{\partial \bv}{\partial
    t}\partial_{33}^2\bv\,d\bx\,
  dt=\frac{1}{2}\|\partial_{3}\bv(t)\|_2^{2}-\frac{1}{2}\|\partial_{3}\bv(0)\|_2^{2}. 
\end{equation*}
\end{lemma}
\begin{proof}
  Note that the assumptions on $\bv$ already imply, by parabolic
  interpolation, that $\bv\in C(\overline I;W^{1,2}_{0}(\Omega))$.  We give an
  elementary proof, by heat regularization, since the direct
  integration by parts is not justified under the given
  assumptions. In fact, we have that
  $\frac{\partial \bv}{\partial t}=\bfzero$ on the boundary, but is it
  not clear if this holds also in the sense of traces. Let us define
  $\bphi:=\frac{\partial \bv}{\partial t}-\Delta\bv\in L^{2}(I\times
  \Omega)$ and $\bpsi:=\bv(0)\in W^{1,2}_{0}(\Omega)$ and approximate
  these functions 
  by sequences of smooth and compactly supported
  functions $\bphi_{n}$ and $\bpsi_{n}$, respectively. Let $\bv_{n}$
  be the solution of boundary initial value problem
  \begin{align}
    \label{eq:heat}
    \begin{aligned}
      \frac {\partial \bv_{n}}{\partial t} -\Delta\bv_{n}&=
      \bphi_{n}
      \qquad&&\text{in } I\times \Omega,
      \\
      \bv_{n} &= \bfzero &&\text{on } I\times\partial \Omega\,,
      \\
      \bv_{n}(0)&=\bpsi_{n}&&\text{in }\Omega\,.
    \end{aligned}
  \end{align}
By energy methods, one obtains directly that there exists a unique
solution $\bv_{n}$ belonging to $L^{2}(I;W^{2,2}(\Omega)\cap
W^{1,2}_{0}(\Omega))\cap W^{1,2}(I;L^{2}(\Omega))$.  Moreover, 
we have
\begin{align*}
  &\|\bv_{n}-\bv_k\|_{L^{2}(I;W^{2,2}(\Omega)\cap
    W^{1,2}_{0}(\Omega))\cap W^{1,2}(I;L^{2}(\Omega))}
  \\
  &\leq
    c\, \|\bphi_{n}-\bphi_k\|_{L^{2}(I\times
    \Omega)}+c\,\|\bpsi_{n}-\bpsi_k\|_{W^{1,2}(\Omega)}\qquad  \text{ for }k,n\in\setN \,,
\end{align*}
which implies that $(\bv_n)$ is a Cauchy sequence in the spaces on the
left-hand side. Let $\bu$ be the limit in $L^{2}(I;W^{2,2}(\Omega)\cap
W^{1,2}_{0}(\Omega))\cap W^{1,2}(I;L^{2}(\Omega))$ of the sequence
$\bv_n$. By passing to the limit in \eqref{eq:heat} we see that $\bu
-\bv$ is a solution of \eqref{eq:heat} with vanishing data. Thus, by
uniqueness we proved that 
\begin{equation}\label{eq:ap}
  \bv_n \to \bv \quad \text { in }L^{2}(I;W^{2,2}(\Omega)\cap
W^{1,2}_{0}(\Omega))\cap W^{1,2}(I;L^{2}(\Omega))\cap C(\overline I;W^{1,2}_{0}(\Omega))\,.
\end{equation}
Next, testing by the ``second order time derivative'' of
$\bv_{n}$, which can be justified with the help of an Galerkin
approximation (cf.~\cite{bdr-7-5, dr-7-5}),  one
gets that
\begin{equation*}
 \left\|\frac{\partial\bv_{n}(t)}{\partial t}\right\|_2^{2}+ \int\limits_{0}^{t}
 \left\|\frac{\partial\nabla \bv_{n}(s)}{\partial t}\right\|_2^{2}\,ds\leq
 c \int\limits_{0}^{t}
 \left\|\frac{\partial\bphi_{n}(s)}{\partial
     t}\right\|_2^{2}\,ds+c\,\|\bpsi_{n}\|^{2}_{W^{2,2}(\Omega) }\leq c(n)\,.
\end{equation*}
This proves that
$\frac{\partial\bv_{n}}{\partial t} \in
L^{2}(I;W^{1,2}(\partial\Omega)) \vnor
L^{2}(I;W^{1/2}(\partial\Omega))$. Thus,
$\frac{\partial\bv_{n}}{\partial t} =\bfzero$ holds in the sense of
traces in $L^{2}(I;W^{1/2}(\partial\Omega))$ and we obtain 
\begin{align*}
    -\int\limits_{0}^{t}\int\limits_{\Omega}\frac{\partial \bv_n}{\partial
  t}\partial_{33}^2\bv_n\,d\bx \, dt
  &= \int\limits_{0}^{t}\int\limits_{\Omega}\frac{\partial^2 \bv_n}{\partial
    t\partial_{3}}\partial_{3}\bv_n\,d\bx\, dt
    =\frac{1}{2}\|\partial_{3}\bv_n(t)\|_2^{2}-\frac{1}{2}\|\partial_{3}\bv_n(0)\|_2^{2}\,. 
\end{align*}
Passing with $n\to\infty$, which is justified by \eqref{eq:ap}, we
proved the assertion. 
\end{proof}

\begin{proposition}
  \label{prop:main}
  Let the same hypotheses as in Theorem~\ref{thm:MT} be satisfied with $\delta >0$ and let
  the local description $\Grenze_{P}$ of the boundary and the localization function $\xi_P$
  satisfy $(b1)$--\,$(b3)$ and $(\ell 1)$ (cf.~Section~\ref{sec:bdr}). Then, there exists
  a constant $C_2>0$ such that the time-regular solution $\bue\in L^{\infty}(I;W^{1,2}_0(\Omega))\cap
  L ^2(I;W^{2,2}(\Omega))$ of
  the approximate problem~\eqref{eq:eq-e} satisfies  for every $P\in
  \partial \Omega$ and for a.e.~$t\in I$
  \begin{equation*}
    \int\limits_0^t  \int\limits_\Omega 
      \xi^2_P |\partial_3\bFn{N}(\bD\bue(s))|^2\,d\bx\,ds
     \leq C\,, 
  \end{equation*}
  provided $r_P<C_2$ in $(b3)$, with $C$ depending on the
  characteristics of $\bS$,
  $\para ^{2-p}\|\bff\|_2^2\,$, $|||\bu_0,\bff|||,\norm{\xi_P}_{
    1,\infty},\norm{\Grenze_{P}}_{C^{2,1}}$, and $C_2$.
\end{proposition}
\begin{proof}
  We adapt the strategy in~\cite[Proposition~3.2]{br-plasticity} to
  the time-dependent problem.  Fix an arbitrary point
  $P\in \partial \Omega$ and a local description $\Grenze=\Grenze_{P}$
  of the boundary and the localization function $\xi=\xi_P$ satisfying
  $(b1)$--\,$(b3)$ and $(\ell 1)$. In the following constants $c, C$
  can always depend on the characteristics of $ \bSn{N}$,
  hence on those of $\bS$, i.e., on $\gamma_3,\gamma_4$, and $p$.  First we observe that
  Proposition~\ref{prop:pFA} and Remark \ref{rem:bal} yield that there exists a
  constant $C_0$, depending only on the characteristics of
  $\bS$ 
  such that
  \begin{equation*}
   \frac{1}{C_0}| \partial_3\bFn{N}(\bD\bue)|^2\leq
   \mathbb{P}_3^{N}(\bD\bue)\qquad \text{a.e.  in }I\times\Omega\,.
  \end{equation*}
  Thus, we get, using also the symmetry of both $\bD\bue$ and $\bSn{N}$,
  \begin{align}\label{eq:j123}
      &\frac{1}{C_0}
      \int\limits_{0}^{t}\int\limits_\Omega 
      \xi^2 |\partial_3\bFn{N}(\bD\bue)|^2\,d\bx \,ds
      \\
      & \leq \int\limits_{0}^{t}\int\limits_\Omega\xi^2 
      \partial_3 \Sn{N}_{\alpha\beta} (\bD\bue) \,\partial_3
      D_{\alpha\beta}\bue\,d\bx \,ds +
      \int\limits_{0}^{t}\int\limits_\Omega\xi^2
      \partial_3\Sn{N}_{3\alpha }(\bD\bue)\,\partial _\alpha
      D_{33}\bue\,d\bx \,ds\notag
      \\
      &\quad +\int\limits_{0}^{t}\int\limits_\Omega \sum_{j=1}^3
      \xi^2\partial_3  \Sn{N}_{j3}(\bD\bue)\, \partial_3^2
        u^{N}_{j}\,d\bx \,ds\notag 
      \\
      & =:\mathcal{J}_{1}+\mathcal{J}_{2}+\mathcal{J}_{3}\,.\notag
  \end{align}
  The terms $\mathcal{J}_{1}$ and $\mathcal{J}_{2}$ can be estimated
  exactly as in~\cite{br-plasticity}, if one replaces
  $\phi''(|\bD\bue|)$ used there with the equivalent quantity $\MC^N(|\bD\bue|)$. Let us sketch the main
  steps. All missing details can be found in \cite{br-plasticity}. To
  treat $\mathcal J_2$ we multiply and divide by
  $\sqrt{\MC^{N}(|\bD\bue|)}$, use Proposition \ref{prop:pFA} and
  Young inequality, to show that, for any given $\lambda>0$, it holds
\begin{equation*}
  \begin{aligned}
    |\mathcal{J}_{2}| \leq &\,\param
    \int\limits_{0}^{t}\int\limits_\Omega\xi^2
    |\partial_3\bFn{N}(\bD\bue)|^2 \,d\bx \,ds +c_{\param^{-1}}\,
    \sum_{\beta=1}^2 \int\limits_{0}^{t}\int\limits_\Omega\xi^2
    |\partial_\beta\bFn{N}(\bD\bue)|^2 \,d\bx \,ds\,,
  \end{aligned}
\end{equation*}
for some constant $c_{\param^{-1}} $ depending  on
$\param^{-1}$. To treat the term $\mathcal J_1$ we first use the
algebraic identity
$ \partial _3 D_{\alpha \beta}\bue= \partial _\alpha D_{3\beta}\bue
+\partial _\beta D_{3\alpha}\bue -\partial_\beta\partial_\alpha u^N_3
$. The first two terms in the resulting equation are treated as
$\mathcal J_2$, while in the term with
$\partial_\beta\partial_\alpha u^N_3 $ we use the definition of
tangential derivatives \eqref{eq:1}. This results in three terms,
where one is again treated as $\mathcal J_2$. This procedure leads
to\footnote{The estimated terms correspond to the terms $A$ and  $B_3$ in \cite{br-plasticity}.}
\allowdisplaybreaks
\begin{align}\label{eq:j1}
  \begin{aligned}
    |\mathcal{J}_{1}| \leq &\, \param
    \int\limits_{0}^{t}\int\limits_\Omega\xi^2
    |\partial_3\bFn{N}(\bD\bue)|^2 \,d\bx \,ds
    \\
    &\quad +c_{\param^{-1}}\, \big (1+\|\nabla \Grenze\|_\infty^2\big )
    \sum_{\beta=1}^2 \int\limits_{0}^{t}\int\limits_\Omega\xi^2
    |\partial_\beta\bFn{N}(\bD\bue)|^2 \,d\bx \,ds
    \\
    &\quad+ \int\limits_{0}^{t}\int\limits_\Omega\xi^2 
      |\partial_3 \bS^{N} (\bD\bue) |\, |\nabla
      ^2\Grenze| \, |\bD\bue|\,d\bx \,ds 
    \\
    &\quad + \bigg |\int\limits_{0}^{t}\int\limits_\Omega\xi^2 
      \partial_3 \Sn{N}_{\alpha\beta} (\bD\bue) \,\partial_\alpha
     \partial _{\tau_\beta} u_3^N\,d\bx \,ds \,\bigg|\,.
  \end{aligned}
\end{align}
In the last but one term in \eqref{eq:j1} we multiply and divide by
$\sqrt{\MC^{N}(|\bD\bue|)}$, use Proposition \ref{prop:pFA}, Young
inequality and
$\MC^N (|\bD\bue|) |\bD\bue|^2 \sim |\bF^N(|\bD\bue|)|^2$
(cf.~Pro\-position \ref{prop:SA-ham}), yielding that it is estimated by
\begin{equation*}
  \param
    \int\limits_{0}^{t}\int\limits_\Omega\xi^2
    |\partial_3\bFn{N}(\bD\bue)|^2 \,d\bx \,ds+c_{\param^{-1}} \, \|\nabla^2 \Grenze\|^2_\infty
  \int\limits_{0}^{t}\int\limits_\Omega|\bF^N(\bD\bue)|^2 
\, d\bx \,ds\,.
\end{equation*}
To handle the last term in \eqref{eq:j1} we want to perform the crucial partial
integration, which reads (neglecting the localization $\xi$)
\begin{align*}
  \int\limits_{0}^{t}\int\limits_\Omega
      \partial_3 \Sn{N}_{\alpha\beta} (\bD\bue) \,\partial_\alpha
     \partial _{\tau_\beta} u_3^N\,d\bx \,ds &=\int\limits_{0}^{t}\int\limits_\Omega 
      \partial_\alpha \Sn{N}_{\alpha\beta} (\bD\bue) \,\partial_3
                                               \partial _{\tau_\beta} u_3^N\,d\bx \,ds
  \\
  &=\int\limits_{0}^{t}\int\limits_\Omega 
      \partial_\alpha \Sn{N}_{\alpha\beta} (\bD\bue) \,
      \partial _{\tau_\beta} D_{33} \bue\,d\bx \,ds\,.
\end{align*}
Again, we multiply and divide by $\sqrt{\MC^{N}(|\bD\bue|)}$, use
Proposition \ref{prop:pFA}, Young inequality, Corollary
\ref{cor:Ftau}, and the definition of the tangential derivatives, yielding that the last term is estimated by
\begin{align*}
  \begin{aligned}
    &c \sum_{\alpha=1}^2\int\limits_{0}^{t}\int\limits_\Omega
    |\partial_\alpha\bFn{N}(\bD\bue)|^2 \,d\bx \,ds+c \sum_{\beta=1}^2
    \int\limits_{0}^{t}\int\limits_\Omega|\partial _{\tau
      _\beta}\bF^N(\bD\bue)|^2 \, d\bx \,ds
    \\
    &\le c \sum_{\alpha=1}^2
    \int\limits_{0}^{t}\int\limits_\Omega|\partial _{\tau
      _\alpha}\bF^N(\bD\bue)|^2 \, d\bx \,ds+c \,\|\nabla
    \Grenze\|_\infty ^2\int\limits_{0}^{t}\int\limits_\Omega 
    |\partial_3\bFn{N}(\bD\bue)|^2 \,d\bx \,ds\,. \hspace*{-5mm}
  \end{aligned}
\end{align*}
The presence of the localization $\xi $ leads to several additional
terms, which all can be handled as in \cite{br-plasticity}. All
together we arrive at the following estimate
\pagebreak
\begin{equation}
  \begin{aligned}\label{eq:j1-2}
    |\mathcal{J}_{1}|+|\mathcal{J}_{2}| \leq &\,\big (\param  
    +c_{\param^{-1}}\,\|\nabla \Grenze\|_\infty^2\big ) \int\limits_{0}^{t}\int\limits_\Omega\xi^2
    |\partial_3\bFn{N}(\bD\bue)|^2 \,d\bx \,ds
    \\
    &+c_{\param^{-1}}\sum_{\beta=1}^2 \int\limits_{0}^{t}\int\limits_\Omega\xi^2
    |\partial_{\tau_\beta}\bFn{N}(\bD\bue)|^2 \,d\bx \,ds
    \\
    &+c_{\param^{-1}} \big(1+\|\nabla\xi\|_\infty^2\big
    ) \int\limits_{0}^{t}\int\limits_\Omega|\bF^N(|\bD\bue|)|^2\, d\bx \,ds\,.
  \end{aligned}
\end{equation}
In this estimate we used for the terms with
$\partial _\beta \bFn{N}(\bD\bue)$ 
the definition of the tangential derivative in \eqref{eq:1} to get
\begin{align}\label{eq:tau}
  \begin{aligned}
    \int\limits_\Omega\xi^2
    |\partial_\beta\bFn{N}(\bD\bue)|^2
    \,d\bx &\le \int\limits_\Omega\xi^2
    |\partial_{\tau_\beta}\bFn{N}(\bD\bue)|^2 
    \,d\bx
    \\
    &\quad + \norm{\nabla \Grenze}^2_\infty\int\limits_\Omega\xi^2
    |\partial_3\bFn{N}(\bD\bue)|^2 
    \,d\bx\,.
  \end{aligned}
\end{align}
Also the term $\mathcal{J}_3$ is treated essentially as in
\cite{br-plasticity}. Since in this step the equation \eqref{eq:eq-e}
is used, in addition we have to handle the term with the time
derivative. More precisely, we  re-write the
equations~\eqref{eq:eq-e} as follows
  \begin{equation*}
    \partial_3 \Sn{N}_{j3}(\bD\bue)=\frac {\partial u^{N}_{j}}{\partial t}-f^j-\partial_\beta
    \Sn{N}_{j\beta}(\bD\bue)\qquad \text{a.e. in     }I\times\Omega\,,
  \end{equation*}
  multiply it by $\partial_{33}\bue$, use the algebraic identity
  \begin{align}
   \partial_{j} \partial_{k}u^{N}_{i}=\partial_{j}D_{i
     k}\bue+\partial_{k}D_{i j} \bue-\partial_{i}D_{j k}
    \bue,\label{eq:alg}
  \end{align}
  treat all terms without the time derivative as $I_3$
  in~\cite[p.~186]{br-plasticity} and integrate by parts the term
  involving $\frac{\partial \bue}{\partial t}$, use 
  Lemma~\ref{lem:time-derivative}, to get the following
  \allowdisplaybreaks
\begin{align}
    \mathcal{J}_{3}
    &=\sum_{j=1}^3 \int\limits_{0}^{t}\intO \xi^2\frac{\partial u^{N}_{j}}{\partial t}\,\partial^2_{33}u^{N}_{j}-\xi^2\big(f^j
    +\partial_\beta \Sn{N}_{j\beta}(\bD\bue)\big)\big
    (2\partial _3D_{j3}\bue - \partial _j D_{33}\bue\big )\,d\bx \,ds\, \notag 
    \\
    &=-\frac{1}{2}\intO \xi^2  |\partial_3\bue (t)|^2\,d\bx\!+\!\frac{1}{2}\intO \xi^2  |\partial_3\bue (0)|^2\,d\bx
    \!-\!2\sum_{j=1}^3 \int\limits_{0}^{t} \intO \xi \partial _3\xi
    \frac{\partial u^{N}_{j}}{\partial t} \partial _3 u^{N}_{j}\, d\bx \,ds \notag 
    \\
    &\qquad +\sum _{j=1}^3 \int\limits_{0}^{t}\intO\xi^2\big(f_{j}
    +\partial_\beta S^N_{j\beta}(\bD\bue)
\big)\big
    (2\partial _3D_{j3}\bue - \partial _j D_{33}\bue\big )\,d\bx \,ds\, \notag 
    \\
    &\leq -\frac{1}{2}\intO
    \xi^2|\partial_3\bue(t)|^2\,d\bx+\frac{1}{2}\intO
    \xi^2|\partial_3\bue(0)|^2\,d\bx\label{eq:j3}
  \\
  &\quad + \big (\param + c_{\lambda^{-1}}\,\|\nabla \Grenze \|_\infty
^2\big )  \int\limits_{0}^{t}\int\limits_\Omega\xi^2|\partial_3\bFn{N}(\bD\bue)|^2\,d\bx \,ds \notag 
    \\
    &\quad
      +c_{\param^{-1}}\sum_{\beta=1}^2\int\limits_{0}^{t}\int
      \limits_\Omega\xi^2|\partial_{\tau_\beta}\bFn{N}(\bD\bue)|^2\,d\bx
      \,ds +c\int\limits_0^t\!\intO
\xi^3
|\partial_{3}\bue|^2\,d\bx\,ds \notag 
    \\
    &\quad
       + c\,\norm{\nabla\xi }_\infty ^2 \int\limits_{0}^{t}\int\limits_\Omega
       \Bigabs{\frac{\partial \bue }{\partial
       t}}^2\, d\bx\,ds+c_{\param^{-1}}\int\limits_{0}^{t}
       \int\limits_{\Omega}\frac{|\bff|^{2}}{\MC^{N}(|\bD\bue|)}\,d\bx
       \,ds\,, \notag 
\end{align}
where we used again \eqref{eq:tau}.
Now we choose in the estimates \eqref{eq:j1-2}, \eqref{eq:j3} first
$\lambda>0$ small enough and then the covering of the boundary
$\partial \Omega$ such that
$\|\nabla \Grenze\|_\infty$ is small enough in order to absorb in the
left-hand side of \eqref{eq:j123} the term involving
$\partial_3\bFn{N}(\bD\bue)$. This way we obtain the following estimate
\begin{equation*}
  \begin{aligned}
    &    \intO\xi^3|\partial_3\bue(t)|^2\,d\bx+ \frac {1}{C_0}\int\limits_0^t\int\limits_\Omega 
\xi^2|\partial_3\bFn{N}(\bD\bue(s))|^2\,d\bx\,ds
    \\
    &\leq     \intO\xi^3|\partial_3\bu_0|^2\,d\bx +c 
    \sum_{\beta=1}^2 \int\limits_0^T
    \int\limits_\Omega\xi^2|\partial_{\tau_{\beta}}\bFn{N}(\bD\bue(s))|^{2}\,d\bx\,ds
    \\
    &\quad+c     \int\limits_0^T\!\intO
    \frac{|\bff(s)|^{2}}{\MC^{N}(|\bD\bue(s)|)}  +|\bF^N(\bD\bue(s))|^2
    +\Bigabs{\frac{\partial \bue(s)}{\partial t}}^2d\bx\,ds
    \\
    &\quad +c\int\limits_0^T\!\intO \xi^3
    |\partial_{3}\bue(s)|^2\,d\bx\,ds
    \end{aligned}
\end{equation*}
with constants depending only on the characteristics of $\bS$,
$\norm{g}_{C^{2,1}}$, and $\norm{\xi}_{1,\infty}$.

Using the uniform estimates \eqref{eq:main-apriori-estimate2}, 
\eqref{eq:est-eps} and the lower bound in
Lemma~\ref{lem:UAm}, which yields 
\begin{equation*}
  \int\limits_0^T\!\intO \frac{|\bff|^{2}}{\MC^{N}(|\bD\bue|)}\,d\bx \,ds \leq
  C\delta^{2-p}\int\limits_0^T\intO    {|\bff|^{2}}\,d\bx\, ds,
\end{equation*}
we get from the last estimate the assertion of Proposition \ref{prop:main}.  
%
%
\end{proof}
Choosing now an appropriate finite covering of the boundary (for the
details see also~\cite{br-reg-shearthin}),
Proposition~\ref{prop:JMAA2017-1} and Proposition~\ref{prop:main} yield the
following result:
\begin{proposition}
\label{thm:estimate_for_ue}
Let the same hypotheses as in Theorem~\ref{thm:MT} with $\delta>0$ be
satisfied. Then, it holds for a.e. $t\in I$
  \begin{equation*}
    \int\limits_0^t 
    \| \nabla \bFn{N}(\bD\bue(s))\|^2_2\,ds\leq 
    C\, 
  \end{equation*}
  with $C$ depending only on 
the characteristics of $\bS$, $|||\bu_0,\bff|||$,
$\delta^{2-p}\|\ff\|_2^2$, $\para^{p-2}$ and
$\partial \Omega$. In particular is $C$ independent of $A_{n}$, $n=1,\dots,N$.
\end{proposition}
\subsection{Multiple Passage to the limit}
From Proposition~\ref{thm:existence_perturbation} and Proposition~\ref{thm:estimate_for_ue}
we obtain the following estimate, uniform with respect to
$A_{n}\geq1$, $n=1,\dots,N$, and valid for a.e.~$t\in I$.
\begin{equation}
  \label{eq:reg-N}
  \begin{aligned}
    \hspace*{-1mm}\|\bue(t)\|_{{1,2}}^{2}\!+\!\|\bFn{N}(\bD\bue(t))\|_{2}^{2}+\!\int\limits_{0}^{t}\left\|\frac{\partial\bue(s)}{\partial
      t}\right\|_{2}^{2}
    \!+\! \| \nabla \bFn{N}(\bD\bue(s))\|^2_2\,ds\leq C\,
  \end{aligned}
\end{equation}
with $C$ depending only on the data of the problem \eqref{eq:pfluid}. 

Note that the functions $\buen{N}$ and $\bFn{N}$ depend (implicitly)
on the parameters $A_{n}$. Since these parameters are relevant for the
various limiting processes, we  now start to write these
dependencies in an explicit way. The uniform estimates for $\bue(t,\bx,A_{1},\dots, A_{N})$ and
$\bF^{N}(\bD\bue(t, \bx,A_{1},\dots, A_{N}))$ are inherited by
taking appropriate limits of the various $A_{n}$. In particular, we will
define (when the limit exists in appropriate spaces)
\begin{equation*}
\buen{N-1}(t,\bx,A_{1},\dots,
A_{N-1}):=\lim_{A_N\to\infty}\bue(t,\bx,A_{1},\dots,  A_{N-1}  ,A_{N}),\,
\end{equation*}
and then inductively
\begin{equation*}
  \buen{n-1}(t,\bx,A_{1},\dots, A_{n-1})=\lim_{A_{n}\to\infty}\buen{n}(t,\bx,A_{1},\dots,A_{n-1},  A_{n})\qquad n=1,\dots,N,
\end{equation*}
 in such a way that the function $\bu:=\buen{0}$ will be shown to be the unique
regular solution to the initial boundary value problem~\eqref{eq:pfluid}.
\begin{proof}[Proof of Theorem~\ref{thm:MT}]
  From estimate \eqref{eq:reg-N} we obtain that $\buen{N}$ is
  uniformly bounded in
  $W^{1,2}(I;L^2(\Omega))\cap L^\infty(I;W^{1,2}(\Omega))$ and that 
  $\bFn{N}(\bD\buen{N})$ is uniformly bounded in
  $L^\infty(I;L^2(\Omega))\cap L^{2}(I;W^{1,2}(\Omega))$.

These bounds directly imply that there exists a sequence $A_{N_{k}}\to\infty$
(which we call again $A_{N}$), a vector field
$\buen{N-1}(t,\bx,A_{1},\dots,A_{N-1})$, and a tensor field $\widehat{\bFn{N-1}}$
\begin{equation}\label{eq:conv}
  \begin{aligned}
    \lim_{A_{N}\to\infty} \buen{N}&=\buen{N-1}&\qquad&\text{weakly in
    }W^{1,2}(I;L^{2}(\Omega))\,, 
    \\
    \lim_{A_{N}\to\infty}    \buen{N}&=
    \buen{N-1}&&\text{
weakly* in }L^{\infty}(I;W^{1,2}(\Omega))\,, 
          \\
 \lim_{A_{N}\to\infty}   \bFn{N}(\bD\buen{N})
 &=\widehat{\bFn{N-1}}&&\text{
   weakly in
      }L^{2}(I;W^{1,2}(\Omega))\,,
\\
 \lim_{A_{N}\to\infty}   \bFn{N}(\bD\buen{N})
&=\widehat{\bFn{N-1}}&&\text{weakly* in
     }L^{\infty}(I;L^{2}(\Omega))\,.
  \end{aligned}
\end{equation}
%
From $\norm{\bFn{N}(\bD\bue)}_{L^{2}(I;W^{1,2}(\Omega)}\leq C$ 
it follows, using Proposition \ref{prop:pFA}, the lower bound on
$\MC^{N}$ proved in Lemma~\ref{lem:UAm}, and the identity
\eqref{eq:alg}, that
\begin{equation}
  \begin{aligned}\label{eq:w22}
    \delta^{p-2}\norm{\nabla^{2}\bue}_{L^2(I;W^{2,2}(\Omega))}\leq C\,
  \end{aligned}
\end{equation}
with $C$ depending only on the data of the problem \eqref{eq:pfluid}, but independent of
$A_{N}$. The estimates \eqref{eq:reg-N}, \eqref{eq:w22} and the
Aubin-Lions compactness lemma imply that (up to a further
sub-sequence)
\begin{equation*}
 \lim_{A_{N}\to\infty} \bD\buen{N}=\bD\buen{N-1}\qquad \text{
    a.e. in $I\times \Omega$ and  strongly in } L^{2}(I\times \Omega) , 
\end{equation*}
for all fixed $A_{n}$, with $n=1,\dots,N-1$.
Next, we observe that since 
\begin{equation*}
  \begin{aligned}
    \lim_{A_{N}\to\infty}\MC^{N}(t)&=\MC^{N-1}(t) \,,
  \end{aligned}
\end{equation*}
 uniformly for $t$ belonging to  compact sets in $\setR^{\geq0}$ (but
in reality even more since $\MC^{N}(t)=\MC^{N-1}(t)$ for all $0\leq t\leq
A_{N}$), it follows by the definition of $\bFn{N}$ and $\bSn{N}$ that a.e.~in $
I\times \Omega$ and for all fixed $A_{n}$, $n=1,\dots,N-1$, there holds
\begin{equation}\label{eq:FS}
  \begin{aligned}
\lim_{A_{N}\to\infty}
\bFn{N}(\bD\buen{N}(A_{1},\dots,A_{N-1},A_{N}))=\bFn{N-1}(\bD\buen{N-1}(A_{1},\dots,A_{N-1}))\,, 
   \\
\lim_{A_{N}\to\infty}
\bSn{N}(\bD\buen{N}(A_{1},\dots,A_{N-1},A_{N}))=\bSn{N-1}(\bD\buen{N-1}(A_{1},\dots,A_{N-1}))\,.
  \end{aligned}
\end{equation}
In fact, by the definition of multiple approximation it follows that
for all given $\bP\in \setR^{3\times 3}$ and for all fixed
$A_{1},\dots,A_{N-1}$ it holds
\begin{equation*}
  \begin{aligned}
    \lim_{A_{N}\to\infty}
    \bFn{N}(\bP,A_{1},\dots,A_{N-1},A_{N})=\bFn{N-1}(\bP,A_{1},\dots,A_{N-1}),
\\
    \lim_{A_{N}\to\infty}
    \bSn{N}(\bP,A_{1},\dots,A_{N-1},A_{N})=\bSn{N-1}(\bP,A_{1},\dots,A_{N-1}),
  \end{aligned}
\end{equation*}
hence 
\begin{equation*}
  \begin{aligned}
&\bFn{N}(\bD\buen{N}(t,\bx,A_{1},\dots,A_{N-1},A_{N}))-\bFn{N-1}(\bD\buen{N-1}(t,\bx,A_{1},\dots,A_{N-1})),
\\
&=
  \bFn{N}(\bD\buen{N}(t,\bx,A_{1},\dots,A_{N-1},A_{N}))-\bFn{N-1}(\bD\buen{N}(t,\bx,A_{1},\dots,A_{N-1},A_N)),
\\
&\quad +  \bFn{N-1}(\bD\buen{N}(t,\bx,A_{1},\dots,A_{N-1},A_N))-\bFn{N-1}(\bD\buen{N-1}(t,\bx,A_{1},\dots,A_{N-1})),
  \end{aligned}
\end{equation*}
and the first line on the right-hand side vanishes for large enough
$A_{N}$,  by the properties of the multiple 
approximation; while the second one converges to zero due to the continuity of
$\bFn{N-1}$ and the point-wise convergence of $\bD\bue$. The same argument
applies also to $\bSn{N}$.

The classical result stating that the weak limit in Lebesgue spaces and the a.e.~limit 
coincide (cf.~\cite{GGZ}) and \eqref{eq:conv} imply that
\begin{equation*}
    \begin{aligned}
      \widehat{\bFn{N-1}} =\bFn{N-1}(\bD\buen{N-1}
      (A_{1},\dots,A_{N-1}))
      \quad\text{in }L^{2}(0,T;W^{1,2}(\Omega))\,.
    \end{aligned}
  \end{equation*}
This identification, the convergences in \eqref{eq:conv}, and the lower
semicontinuity of norms proves that, for a.e.~$t\in I$, it holds
\begin{equation}\label{eq:reg-N-1}
  \begin{aligned}
   & \|\buen{N-1}(t)\|_{W^{1,2}}^{2}+\|\bFn{N-1}(\bD\buen{N-1}(t))\|_{2}^{2}
   \\
   &\quad +\int\limits_{0}^{t}\left\|\frac{\partial\buen{N-1}(s)}{\partial
      t}\right\|_{2}^{2}+  \| \nabla
    \bFn{N-1}(\bD\buen{N-1}(s))\|^2_2\,ds\leq C\,
  \end{aligned}
\end{equation}
with a constant $C$ depending on the data of the problem
\eqref{eq:pfluid}, but  independent of $A_{n}$, for $n=1,\dots,N-1$.

We have now to pass to the limit in the weak
formulation~\eqref{eq:weak-eps} of the approximate
problem~\eqref{eq:eq-e}. Since, in view of \eqref{eq:conv}, we easily deal with the time derivative
and the right-hand side $\ff$, the crucial point is the justification
of the limit
\begin{equation}
\label{eq:relevant-limit}
  \int\limits_0^T\hskp{\bSn{N}(\bD\buen{N}(t))}{\bD\bw}\,\psi (t)\,dt\to   
  \int\limits_0^T\hskp{\bSn{N-1}(\bD\buen{N-1}(t))}{\bD\bw}\,\psi (t)\,dt,
\end{equation}
for all $\psi \in C_0^\infty (I)$ and all $\bw\in C^\infty_0(\Omega)$.
At the moment we already know that
$\lim_{A_{N}\to\infty}\bSn{N}(\bD\buen{N})=\bSn{N-1}(\bD\buen{N-1})$
holds a.e.~in $I\times \Omega$. Thus, to conclude it is sufficient to show that
$\bSn{N}(\bD\bue)$ is bounded uniformly with respect to $N$ in
$L^q(I\times \Omega)$, for some $q>1$. To this end we observe that 
Corollary~\ref{cor:UAm}, Proposition \ref{prop:hammer-phi}, the
definition of $\bF_{\function_{q_N,\delta}}$ in \eqref{def:F1}, and
$q_N \ge 2$ imply that for all $\bP\in\setR^{3\times 3}$ there holds
\begin{align}\label{eq:un}
  \begin{aligned}
    | \bFn{N}(\bP)|^2&\ge c\,\para^{p-q_N}
    |\bF_{\function_{q_N,\delta}}(\bP)|^2 = c\,\para^{p-q_N}
    (\delta+|\bP^{\sym}|)^{q_{N}-2}|\bP^{\sym}|^2
    \\
    &\ge c\,\para^{p-q_N} |\bP^{\sym}|^{q_{N}}\,.
  \end{aligned}
\end{align}
The a priori bound \eqref{eq:reg-N} and parabolic embedding imply that
$\bF^N(\bD \bue)$ is bounded in $L^{\frac {10}3}(I\times \Omega)$ by a
constant depending only on the data of problem \eqref{eq:pfluid}. This
together with \eqref{eq:un} and $\frac 53q \ge q+\frac 43$, valid for all
$q\ge 2$, implies  
\begin{equation*}
  \begin{aligned}
    \|\bD\buen{N}\|_{L^{q_N +\frac 43}(I\times \Omega)}\leq C\,
  \end{aligned}
\end{equation*}
with a constant independent of $A_{N}$.
Corollary \ref{cor:UAm} also implies that 
\begin{equation*}
  \begin{aligned}
      |\bSn{N}(\bD\bue(t,\bx))|&\leq
     c\, A_{N-1}^{p-q_{N-1}}\,
    (\function_{q_{N-1},\para})' (\abs{\bD\bue(t,\bx)}) \,
\\
&    \leq C \, A_{N-1}^{p-q_{N-1}}\,
\big(\delta^{q_{N-1}-1}+|\bD\buen{N}(t,\bx)|^{q_{N-1}-1}\big)\,. 
  \end{aligned}
\end{equation*}
Hence, the latter estimates prove that
\begin{equation*}
  \begin{aligned}
    &\norm{\bSn{N}(\bD\bue)}_ {L^{(4/3+q_{N})/
        (q_{N-1}-1)}(I\times \Omega)}\leq C (A_{N-1}),
  \end{aligned}
\end{equation*}
where the constant $C$ depends on the data of the problem
\eqref{eq:pfluid}, on $A_{N-1}$, but is independent of $A_{N}$.
Thus, we can infer that there exists
$\widehat{\bSn{N-1}}$ 
such that (up possibly to a further relabelled sub-sequence)
\begin{equation}\label{eq:S*}
  \lim_{A_{N}\to\infty}
  \bSn{N}(\bD\buen{N})=\widehat{\bSn{N-1}}\qquad\text{ weakly in
  }{L^{(4/3+q_{N})/
      (q_{N-1}-1)}(I\times \Omega)}\,,
\end{equation}
provided  ${(4/3+q_{N})/ (q_{N-1}-1)}>1$, which is equivalent to 
\begin{equation*}
  q_{N-1}-q_{N}<\frac{7}{3}\,,
\end{equation*}
which motivated the choice of $q_n$ in Definition \ref{def:spec}.
Using again the classical result stating that the weak limit in Lebesgue spaces and the a.e.~limit 
coincide (cf.~\cite{GGZ}) we infer from \eqref{eq:FS} and \eqref{eq:S*} that
\begin{equation*}
  \begin{aligned}
\widehat{\bSn{N-1}}=\bSn{N-1}(\bD\buen{N-1} (A_{1},\dots,A_{N-1}))
    \quad\text{in }L^{(\frac 43+q_{N})/ (q_{N-1}-1)}(I\times \Omega)
    \,,
  \end{aligned}
\end{equation*}
which in turn implies \eqref{eq:relevant-limit}. Thus we proved that
$\buen{N-1}$ satisfies \eqref{eq:reg-N-1} and 
\begin{align*}
    &\int\limits _0^T\!\Bighskp{\frac {\partial\buen{N-1}(t)}{\partial t}
    }{\bfw}\,\psi (t) \, dt
    \!+\!\int\limits_0^T\!\hskp{\bSn{N-1}(\bD\buen{N-1}(t))}{\bD\bfw}
    \,\psi(t)\,dt =\int\limits_0^T\!\hskp{\bff(t)}{\bfw}\, \psi(t)\,
    dt\,,
\end{align*}
for all $\psi \in C_0^\infty (I)$ and all $\bw\in
C^\infty_0(\Omega) $.


\medskip

At this point we can repeat exactly the same argument by replacing
$N$ with $N-1$. Thus,  one obtains inductively that for all
$n=1,\dots,N-1$ there holds 
\begin{equation*}
\int\limits _0^T\Bighskp{\frac {\partial\buen{n-1}(t)}{\partial t}
  }{\bfw}\,\psi (t) \, dt
  +\int\limits_0^T\hskp{\bSn{n-1}(\bD\buen{n-1}(t))}{\bD\bfw}\,\psi
  (t)\,dt 
  =\int\limits_0^T\hskp{\bff(t)}{\bfw}\, \psi(t)\, dt\,,
\end{equation*}
for all $\psi \in C_0^\infty (I)$ and all $\bw\in
C^\infty_0(\Omega) $. After $N$ iterations we find, using also the
density of $C^\infty_0(\Omega) $ in $W^{1,p}_0(\Omega)$ in the last step,  that
the vector field $\bu^{0}=:\bu$ is a  regular solution of the original
problem problem~\eqref{eq:pfluid}. This finishes the proof of Theorem
\ref{thm:MT}. 
\end{proof}

Let us finish with stating the corresponding result to Theorem
\ref{thm:MT} in the steady case. This result can be proved, with many simplifications due to the
absence of the time derivative and the better embedding results in the
steady case (cf.~Section~\ref{sec:ske}), exactly in the same way as the
unsteady result Theorem
\ref{thm:MT}. Thus, we have the following result:
 \begin{theorem}
  \label{thm:MTs}
  Let $\Omega\subset\setR^3$ be a bounded domain with $C^{2,1}$
  boundary, and assume that $\bff \in L^{2}(\Omega)$.  Let the
  operator $\bS$, derived from a potential $\pot $, have
  $(p,\delta)$-structure for some $p\in(2,\infty)$, and
  $\delta\in(0, \infty)$ fixed but arbitrary.
  
  Then, there exists a unique regular solution of the steady version
  of the system~\eqref{eq:pfluid}, i.e., $\bfu \in W^{1,p}_0
  (\Omega)$ fulfils for all $ \bw\in C^{\infty}_{0}(\Omega)$
  \begin{equation*}
      \int\limits_{\Omega}\mathbf{S}(\bD\bu)\cdot\bD\bw\,d\bx=\int\limits_{\Omega}\mathbf{f}\cdot
    \bw\,d\bx\,,
  \end{equation*}
  and satisfies $\bF(\bD\bu) \in W^{1,2}(\Omega)$ 
  with norm  depending  only on the characteristics of
  $\bS$, $\delta^{-1}$, $\Omega$, and $\|\bff\|_2$.
\end{theorem}

\section*{Acknowledgement}
Luigi C. Berselli was partially supported by a grant of the group
GNAMPA of INdAM. 

\def\cprime{$'$} \def\cprime{$'$} \def\cprime{$'$}

\end{document}